\pdfoutput=1
\documentclass[review]{elsarticle}
\usepackage{setspace}
\usepackage{longtable}
\usepackage{lineno,hyperref}
\usepackage{amssymb,amsfonts,amsthm}
\usepackage{amsmath}
\usepackage{graphicx}
\usepackage{float}
\usepackage{caption}
\usepackage{subcaption}
\usepackage{parselines}
\usepackage[margin=1.0in]{geometry}
\begin{document}

\title{Neural Networks with Local Converging Inputs (NNLCI) for Solving Conservation Laws, Part II: 2D Problems}
\author{Haoxiang Huang\footnotemark[1],\; Yingjie Liu\footnotemark[2] \; and Vigor Yang\footnotemark[3]}

\begin{abstract}
In our prior work ~\cite{2CGNN_1D_RiemannProb_21}, neural network methods with inputs based on domain of dependence and a converging sequence were introduced for solving one dimensional conservation laws, in particular the Euler systems. To predict a high-fidelity solution at a given space-time location, two solutions of a conservation law from a converging sequence, computed from low-cost numerical schemes, and in a local domain of dependence of the space-time location, serve as the input of a neural network. In the present work, we extend the methods to two dimensional Euler systems and introduce variations. Numerical results demonstrate that the methods not only work very well in one dimension ~\cite{2CGNN_1D_RiemannProb_21}, but also perform well in two dimensions. Despite smeared local input data, the neural network methods are able to predict shocks, contacts, and smooth regions of the solution accurately. The neural network methods are efficient and relatively easy to train because they are local solvers.                 
\end{abstract}

\maketitle%\thispagestyle{empty}

\renewcommand{\thefootnote}{\fnsymbol{footnote}}
%\footnotetext[1]{{\tt E-mail: hcwong@gatech.edu}. Woodruff School of Mechanical Engineering, Georgia Institute of Technology, Atlanta, GA 30332, USA}
\footnotetext{1. {\tt E-mail: hcwong@gatech.edu}. Woodruff School of Mechanical Engineering, Georgia Institute of Technology, Atlanta, GA 30332, USA.}
%\footnotetext[2]{{\tt E-mail: yingjie@math.gatech.edu}. School of Mathematics, Georgia Institute of Technology, Atlanta, GA 30332, USA. Research is supported in part by NSF grant DMS-1522585.}
\footnotetext{2. {\tt E-mail: yingjie@math.gatech.edu}. School of Mathematics, Georgia Institute of Technology, Atlanta, GA 30332, USA.}
%\footnotetext[3]{{\tt Email: vigor.yang@aerospace.gatech.edu}. Daniel Guggenheim School of Aerospace Engineering, Georgia Institute of Technology, Atlanta, GA 30332, USA}
\footnotetext{3. {\tt Email: vigor.yang@aerospace.gatech.edu}. Daniel Guggenheim School of Aerospace Engineering, Georgia Institute of Technology, Atlanta, GA 30332, USA.}
\renewcommand{\thefootnote}{\arabic{footnote}}

\makeatletter
\def\ps@pprintTitle{%
  \let\@oddhead\@empty
  \let\@evenhead\@empty
  \let\@oddfoot\@empty
  \let\@evenfoot\@oddfoot
}
\makeatother

%\linenumbers
\section{Introduction}
Artificial neural networks~\cite{NN89} are an important tool for computations in science and engineering.
Indeed, many approaches have recently been developed that incorporate artificial neural networks for solving hyperbolic equations. For example, Raissi \textit{et al.}~\cite{RAISSI2019686} employed physics-informed neural networks (PINN) by constraining neural networks with physics informed conditions, such as initial conditions, boundary conditions and functional forms of partial differential equations (PDEs), and utilizing automatic differentiation~\cite{BaydinPearlmutterRadulSiskind17}. This method has achieved many successes in data driven methods for predicting turbulent mixing, vortex induced vibration (VIV) with given governing PDEs, including the Navier-Stokes equations~\cite{ RaissiBabaeeGivi19,RAISSI2019686, raissiJFM19, Raissi1026}, hypersonic flow~\cite{lou2020physicsinformed},  electro-convection~\cite{MFNN_CompxFluids21} and so on. 
In~\cite{Constraint_NN20}, the Rankine-Hugoniot jump conditions are added as a constraint to the loss function of the neural network for solving Riemann problems. 
In~\cite{ChenT95, DeepOnet21}, finite expansions of neural networks that can be trained off-line are introduced to form a mapping, which can map the initial values and a spatial location to the high-fidelity solution at the location in a later time.  
There are also approaches using neural networks trained off-line to predict key parameters of a numerical scheme. In \cite{Control_Oscillation_high-order_Galerkin20, Control_Oscillation_Spectral_Method21}, neural networks are used to detect discontinuities and determine the size of artificial viscosity needed in the presence of discontinuities for a scheme, using a local solution as the input. In ~\cite{NN_Based_LImiter20}, neural networks are used to detect discontinuities and tune the slope limiter of a scheme, using local solution as the input.
Another approach is to use a numerical solution in a region computed on a coarse grid as input to predict the high fidelity solution in the region.
In \cite{MFNNliu19}, a low-cost low-fidelity solution is used in a neural network to predict its difference from the high-fidelity solution. In \cite{Nguyen2020ASP}, the gradient of the numerical solution of a second order wave equation on a coarse grid is used to predict the solution on a fine grid.

In our prior work~\cite{2CGNN_1D_RiemannProb_21}, a novel neural network method is introduced to solve conservation laws whose solutions may contain discontinuities. We wish to use a local low-cost solution as the input of a neural network to predict a high-fidelity solution at a given space-time location. In order for the neural network to distinguish a numerically smeared discontinuity from a steep smooth solution in its input, two approximate solutions of the conservation laws from a sequence (converging to the solution), computed from low-cost numerical schemes, and in a local domain of dependence of the space-time location, serve as the input instead. The neural network is now expected to distinguish the numerical discontinuity from a smooth solution in its input and make the correct prediction because the former becomes increasingly steeper in a converging sequence in the input, and the latter does not. This works very well, not only for discontinuities, but also for smooth areas of the solution. The cost is low because it is a local post-processing-type solver, and we can use low-cost (first-order) schemes on coarse grids to compute inputs. There are many ways to compute a converging sequence to the solution. We first compute first-order numerical solutions on two coarse grids (one coarser than the other), and then use the results in a carefully selected local domain of dependence of the space-time location as input. This is referred to as the 2-coarse-grid neural network (2CGNN). Alternatively, instead of using two coarse grids, we can use a low-cost scheme on one coarse grid to solve the conservation laws perturbed with two diffusion coefficients, and carefully apply the numerical solutions to a local domain of dependence of the space-time location as input. This is referred to as the 2-Diffusion-Coefficient neural network (2DCNN).

In this study, we introduce some variations of \cite{2CGNN_1D_RiemannProb_21} and
extend it to two dimensions. Several configurations of 2D Riemann problems from \cite{PDLax98} are studied, and the prediction results demonstrate the effectiveness of this method in capturing both the discontinuities and the smooth areas of the solution. TensorFLow~\cite{abadi2016tensorflow} is used for the neural networks involved in the numerical experiments. Compared to widely used traditional numerical methods, such as MUSCL~\cite{vanLeer79}, ENO~\cite{ENO87,ShuOsh88}, and WENO~\cite{WENO,JiShu96}, the proposed neural network methods seem to be particularly useful for applications that need repetitive computations with varying parameters. 

The paper is structured as follows.  Sec.~\ref{Sec: 1D-variants} introduces some new variants of the neural network method. 2-Coarse-Grid neural networks for 2D are introduced in Sec.~\ref{2CGNN-4-2D}. 2-Diffusion-Coefficient neural networks for 2D are introduced in Sec.~\ref{Sec: 2DCNN-2D}. Sec.~\ref{Sec: summary} summarizes numerical errors. Conclusions are presented in Sec.~\ref{Sec: conclusion}.

\section{New Variants of the Neural Network Method for 1D  Problems}
\label{Sec: 1D-variants}
Consider a scalar conservation law 
\begin{equation}
\label{1D-cons-law}
   \frac{\partial U}{\partial t} + \frac{\partial f(U)}{\partial x} = 0,  x\in\Omega\subset \mathcal{R},  t\in[0, T],
\end{equation}
and the 1-D Euler equations for an ideal gas    
\begin{equation}
\label{1D-Euler}
    \frac{\partial}{\partial t}\left(\begin{array}{c} \rho \\ \rho v \\E \end{array}\right) +\frac{\partial }{\partial x} \left(\begin{array}{c} \rho v\\ \rho v^2 \\ v( E + p)\end{array}\right) = 0,\; x\in\Omega\subset \mathcal{R}, \; t\in[0, T]~,
\end{equation}
where $\rho$, $u$ and $p$ are density, velocity,
and pressure, respectively, $\Omega=[a,b]$ is an interval,
\begin{equation}
    E = \frac{p}{\gamma - 1}+\frac{1}{2}\rho v^2~,
\end{equation}
$\gamma$ is the specific heat ratio and its value is $1.4$ throughout all the cases studied in the paper.  

We first describe 2CGNN introduced in \cite{2CGNN_1D_RiemannProb_21} as follows. Let $[a,b]$ be partitioned with a coarse uniform grid $a=x_0<x_1<\dots <x_M=b$ having spatial grid size $\Delta x=x_1-x_0$ and let the time step size be $\Delta t$. Refine the grid to obtain a finer uniform grid with spatial grid size $\frac12 \Delta x$ and time step size $\frac12 \Delta t$.  Let $L$ be a low-cost scheme used to compute (\ref{1D-cons-law}) on both grids. 
Given a grid point $x_{i'}$ at time $t_{n'}$ (on the coarsest uniform grid) where the solution is to be predicted by a neural network, we choose the coarsest grid solution (computed by $L$) at $3$ points $x_{i'-1}$, $x_{i'}$ and $x_{i'+1}$ at  time level $t_{n'-1}$ and also at point $(x_{i'}, t_{n'})$ as the first part of the input, and the finer grid solution (also computed by $L$) at the same space-time locations as the second part of input. Note that the chosen $4$ space-time locations of either grid enclose a local (space-time) domain of dependence of the exact solution at $(x_{i'}, t_{n'})$ (with $\Delta t$ satisfying the CFL restriction.)
Since the two parts of the input solution have different levels of approximation to the exact or reference solution, the neural network utilizes the information to extrapolate a prediction of the exact or reference solution. Denote the first part of the input as
$$u^{n'-1}_{i'-1}, u^{n'-1}_{i'}, u^{n'-1}_{i'+1}, u^{n'}_{i'}~,$$
and 
the second part of input as
$$u^{n''-2}_{i''-2}, u^{n''-2}_{i''}, u^{n''-2}_{i''+2}, u^{n''}_{i''}~.$$
Note that the space-time index $(i',n')$ in the coarsest grid refers to the same location as $(i'',n'')$ does in the finer grid, $(i'-1, n'-1)$ refers to the same location as $(i''-2, n''-2)$ does, and so on. Suppose we are interested in the predicted solution at $(x,t)$ which is referred to as $(i',n')$ in the coarsest grid,
the input of 2CGNN is 
\begin{equation}
\label{standard-2CGNN-input}
    \{u^{n'-1}_{i'-1}, u^{n'-1}_{i'}, u^{n'-1}_{i'+1}, u^{n'}_{i'}, u^{n''-2}_{i''-2}, u^{n''-2}_{i''}, u^{n''-2}_{i''+2}, u^{n''}_{i''}\}~,
\end{equation}
called ``input of $u$,'' and the corresponding output of 2CGNN
is the predicted solution at $(x,t)$.
For the Euler system (\ref{1D-Euler}), the input and output of 2CGNN are made up of the corresponding ones for each prime variable. For example, if the input is
the vector 
$$\{ {\rm input}\; {\rm of}\; \rho,\;
{\rm input}\; {\rm of}\; v,\;
{\rm input}\; {\rm of}\; p\}
$$
with $8\times 3=24$ elements, the corresponding output will be 
$\{\rho, v, p\}$
at $(x,t)$ with $3$ elements.

The loss function measures the difference between the output and the reference solution corresponding to the input, and is defined as follows.
$$
\begin{array}{ccc}
{\rm Loss}&=&\sum_{k}\|({\rm output} \; {\rm corresponding}\; {\rm to}\; k^{th}\; {\rm set}\; {\rm of}\; {\rm input})- \\
&& ({\rm reference}\; {\rm solution}\; {\rm corresponding}\; {\rm to}\; k^{th}\; {\rm set}\; {\rm of}\; {\rm input}) \|_2^2~,
\end{array}
$$
where $\|\cdot\|_2$ is the $2$-norm, and the summation goes through every set of input in the training data.

\subsection{Cross-training and Incorporating the Time Step Size into the Input}
\label{Cross-training}
The inputs in \cite{2CGNN_1D_RiemannProb_21} imply that $\Delta t$ is fixed as the size of time step of the low-cost scheme used in computing inputs on the coarsest grid. To account for training processes that require different time step sizes, we can treat the time step size $\Delta t$ as part of the input.
The modified input for problem (\ref{1D-cons-law}) is
\begin{equation}
\label{cross-training-input1}
    \{u^{n'-1}_{i'-1}, u^{n'-1}_{i'}, u^{n'-1}_{i'+1}, u^{n'}_{i'}, u^{n''-2}_{i''-2}, u^{n''-2}_{i''}, u^{n''-2}_{i''+2}, u^{n''}_{i''}, \Delta t\}~,
\end{equation}
and the output is unchanged. Similarly, the input for the Euler system~(\ref{1D-Euler}) can include $\Delta t$
and the output is unchanged.
Incorporating the time step size into the input allows the low-cost scheme to use various time step sizes
in computing inputs for the training data, and in computing inputs for the neural network to make predictions.

The associated neural network for 2CGNN typically consists of $6$ hidden layers, and each layer has $180$ neurons. When performing ``cross-training'' (i.e., the training data covers several different problems, such as the Lax and Sod problems~\cite{ShuLec98}), the neural network can consist of $6$ hidden layers with $255$ neurons per layer, or even $300$ neurons per layer. When $\Delta t$ is chosen to be the same in the Lax and Sod problems, there is no need to treat $\Delta t$ as part of the input. The
neural network structure to predict final-time of Lax and Sod problems consists of $5$ hidden layers with $66$ neurons per layer.

During the training process, the neural network minimizes the difference between the outputs of the neural network and a reference solution by using an Adam optimizer first and an L-BFGS optimizer thereafter in TensorFLow (with the number of iterations of optimization procedure under $50000$ each.) After the training is done, the neural network is used to predict a solution (different from the training data), given an input computed by the same low-cost scheme(s) and grids, which are used to compute inputs of the training data. We use the trained neural network to predict final solutions for $14$ initial values of the Euler system which include the original initial values of the Lax and Sod problems, $\pm3\%$, $\pm5\%$, and $\pm7\%$ perturbations of their initial values. 

In order to generate the training data, we use a first order scheme on
the coarsest uniform grid ($50$ cells) and the finer uniform grid ($100$ cells) to compute the input data from several initial values of the Euler system, which include $\pm2\%$, $\pm4\%$, $\pm6\%$, $\pm8\%$ and $\pm10\%$ perturbations of the initial values of the Lax and Sod problems.
The high resolution reference solutions of the training data are computed on a uniform grid
with $200$ cells by a 3rd order finite volume scheme using non-oscillatory hierarchical reconstruction (HR) limiting ~\cite{LiuShu07b} and partial neighboring cells~\cite{XuLiuShu09} in HR (note that solution values at grid points need to be interpolated from cell averages.) This scheme is used to compute reference solutions for all 1D examples in the paper. The first order scheme used for generating the inputs is  (\ref{leapfrog-diffusion-splitting}) as in \cite{2CGNN_1D_RiemannProb_21}. The time step size for the first order scheme is fixed during the evolution in time, e.g., $\Delta t$ for the $50$-cell grid and $\frac12 \Delta t$ for the $100$-cell grid, satisfying the CFL condition.
For (\ref{1D-cons-law}), the first order leapfrog and diffusion splitting scheme is
\begin{equation}
\label{leapfrog-diffusion-splitting}
\left \{
\begin{array}{l}
    \frac{\tilde{U}_i - U^{n-1}_i}{2\Delta t} +\frac{{f(U)}|^{n}_{i+1}-{f(U)}|^{n}_{i-1}}{2\Delta x} =0~, \\
    \frac{U^{n+1}_i - \tilde{U}_i}{\Delta t}- \alpha\frac{\tilde{U}_{i+1} - 2\cdot \tilde{U}_i + \tilde{U}_{i-1}}{\Delta x^2}=0~,
\end{array}
\right .
\end{equation}
\\
where $\alpha = \Delta x$.

 Figures~\ref{2CGNN: cross-training, Lax prob.} and \ref{2CGNN: cross-training, Sod prob.} depict the predicted density profiles of the Lax and Sod problems using 3 different types of training approach.

The first training strategy is to use input type (\ref{standard-2CGNN-input}) with training data generated from perturbed Lax and Sod problems described above. Note that the same spatial and temporal grid sizes must be used
when computing inputs for the two problems. The predicted results are shown in the right graphs in Fig.~\ref{2CGNN: cross-training, Lax prob.} and \ref{2CGNN: cross-training, Sod prob.}. The second training strategy is to use input type (\ref{cross-training-input1}) with training data generated from
perturbed Lax and Sod problems as described above. Since the time step size is treated as part of the input, 
we can use the same spatial grid sizes but different time step sizes 
when computing inputs for the two problems. The two different time step sizes help the neural network switch between the two problems and the prediction results are improved as shown in the left graphs of Figs.~\ref{2CGNN: cross-training, Lax prob.} and \ref{2CGNN: cross-training, Sod prob.}. 
In order to make a smoother transition in the training data, we also introduce two evenly distributed intermediate time step sizes between the above two time step sizes, and compute the corresponding input data for the two problems. The prediction results are shown in the middle graphs of Figs.~\ref{2CGNN: cross-training, Lax prob.} and \ref{2CGNN: cross-training, Sod prob.}.

\begin{figure}[H]\centering
\begin{subfigure}[b]{0.33\textwidth}
    \centering
    \includegraphics[width=1.0\linewidth]{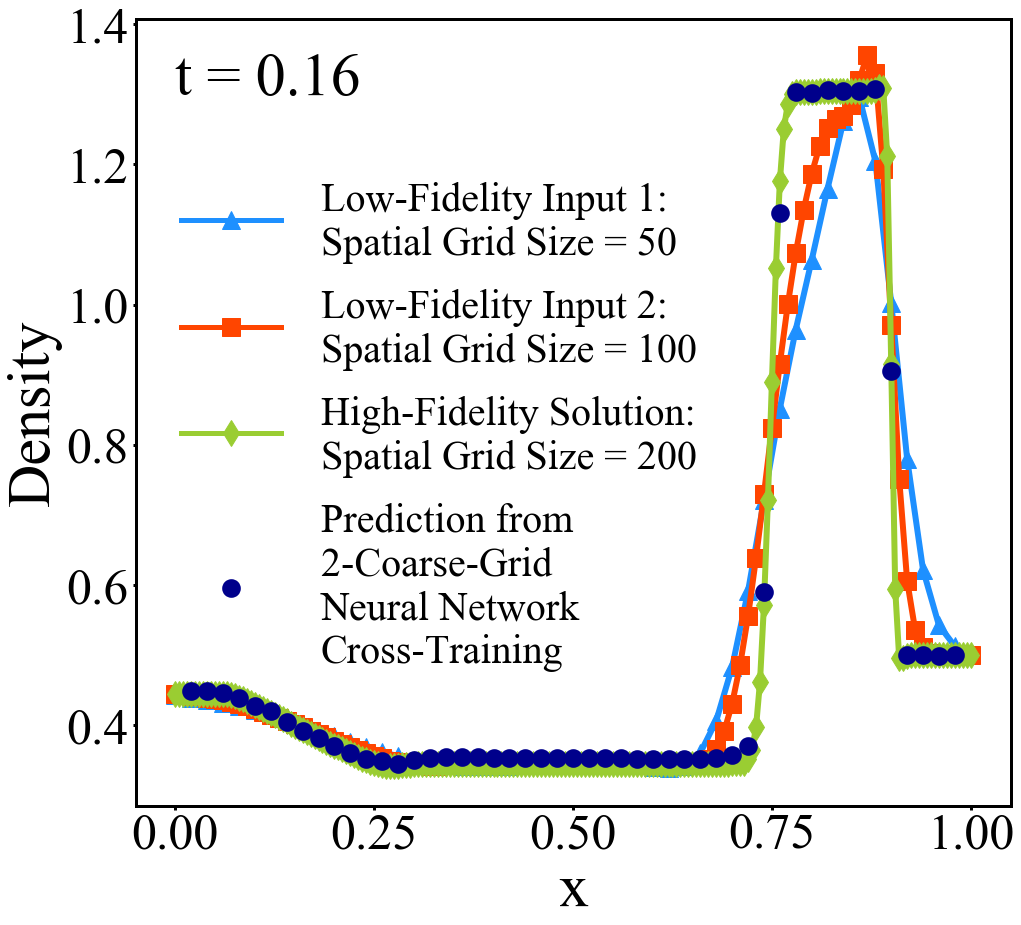}
    \caption{}
\end{subfigure}
\hfill
\begin{subfigure}[b]{0.33\textwidth}
    \centering
    \includegraphics[width=1.0\linewidth]{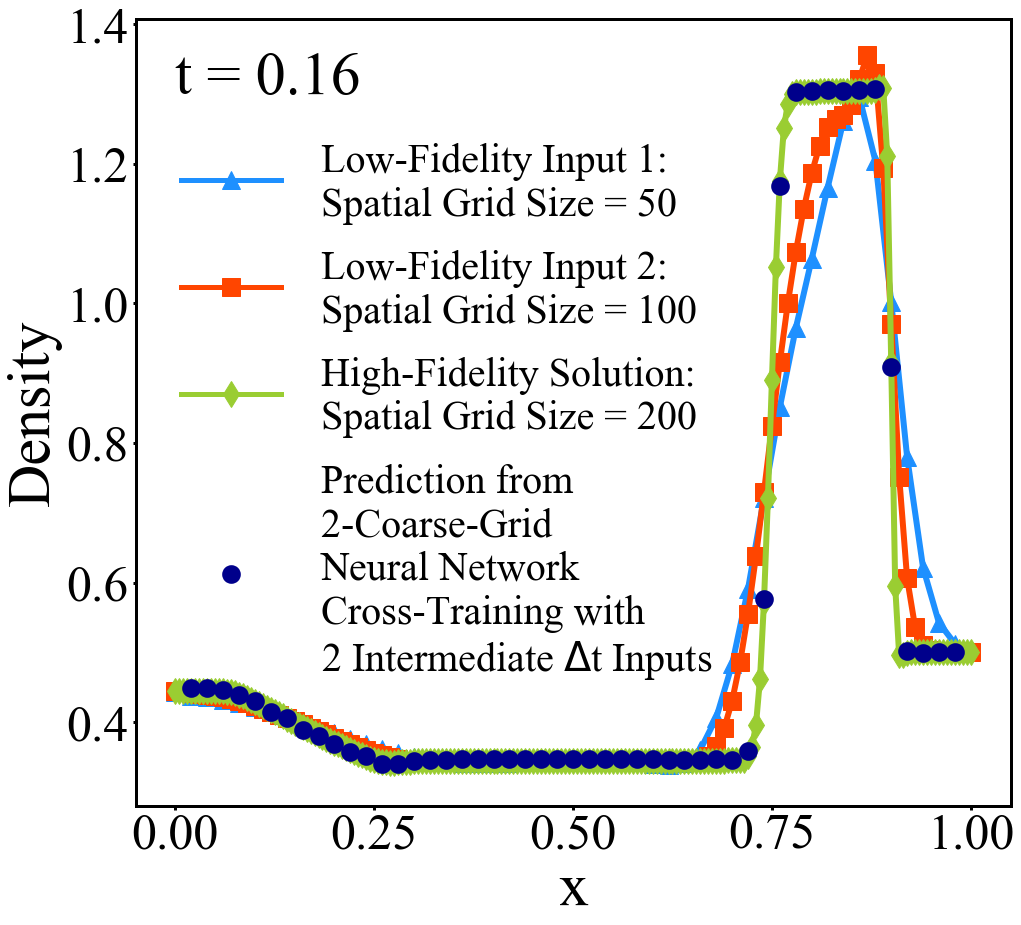}
    \caption{}
\end{subfigure}
\hfill
\begin{subfigure}[b]{0.327\textwidth}
    \centering
    \includegraphics[width=1.0\linewidth]{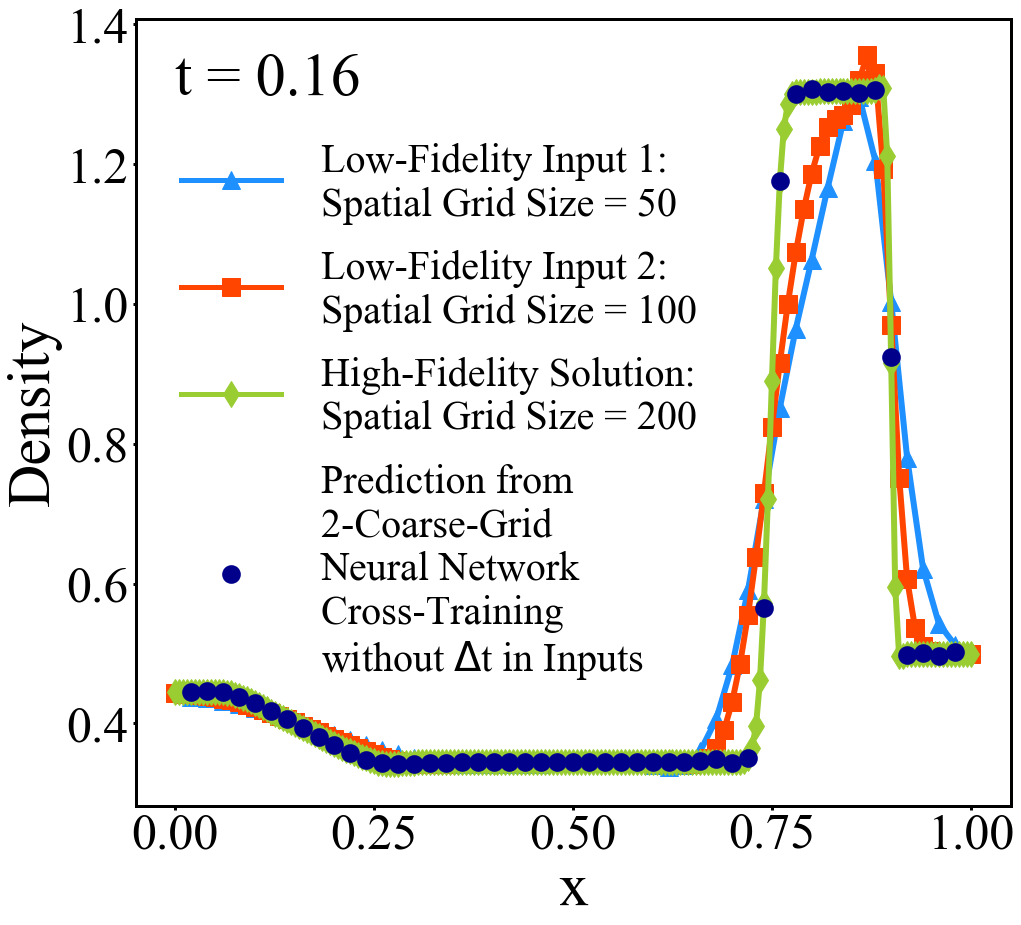}
    \caption{}
\end{subfigure}
\caption{2CGNN Cross-Training predictions of density profiles at final-time $(t=0.16)$ solution of {\bf Lax problem} (dark blue), low-fidelity input solutions (blue and red) by leapfrog and diffusion splitting scheme (\ref{leapfrog-diffusion-splitting}) on $2$ different grids (with $50$ and $100$ cells resp.), and ``exact'' (reference) solution (green): (a) input type~(\ref{cross-training-input1}), (b) input with 2 intermediate time steps, (c) input type~(\ref{standard-2CGNN-input}).}
\label{2CGNN: cross-training, Lax prob.}
\end{figure}

\begin{figure}[H]\centering
\begin{subfigure}[b]{0.33\textwidth}
    \centering
    \includegraphics[width=1.0\linewidth]{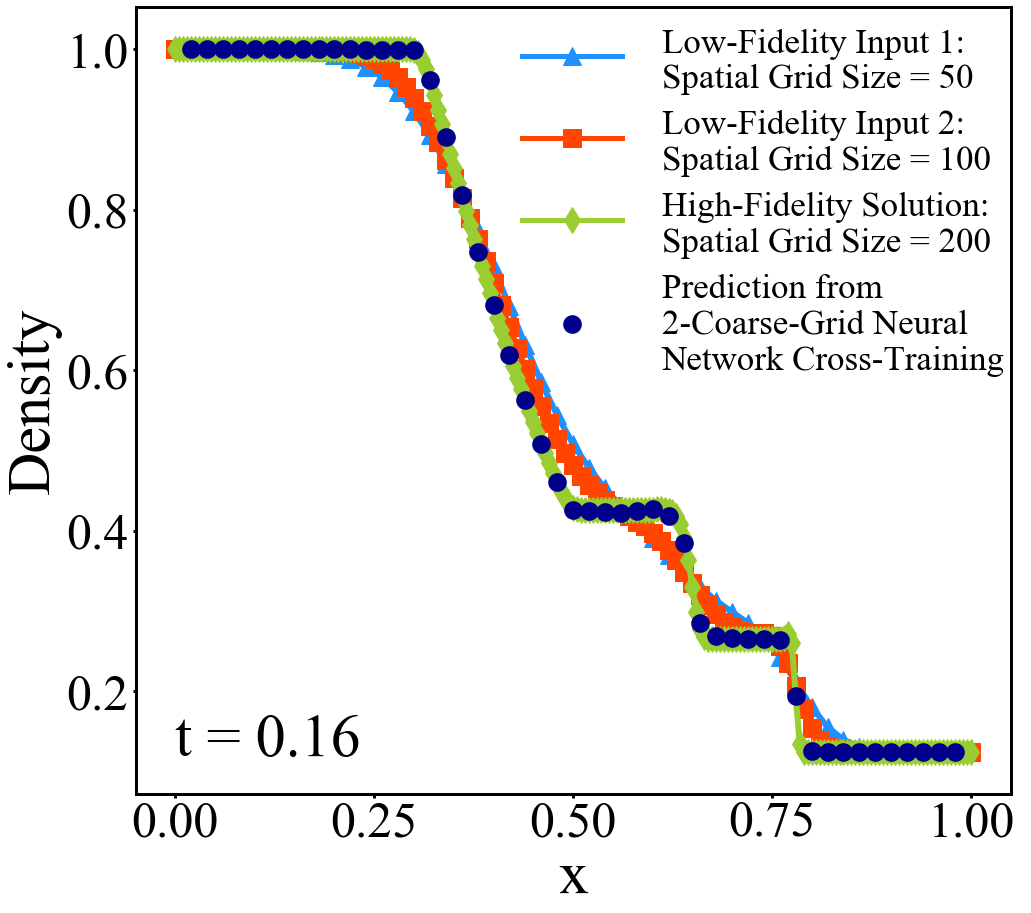}
    \caption{}
\end{subfigure}
\begin{subfigure}[b]{0.33\textwidth}
    \centering
    \includegraphics[width=1.0\linewidth]{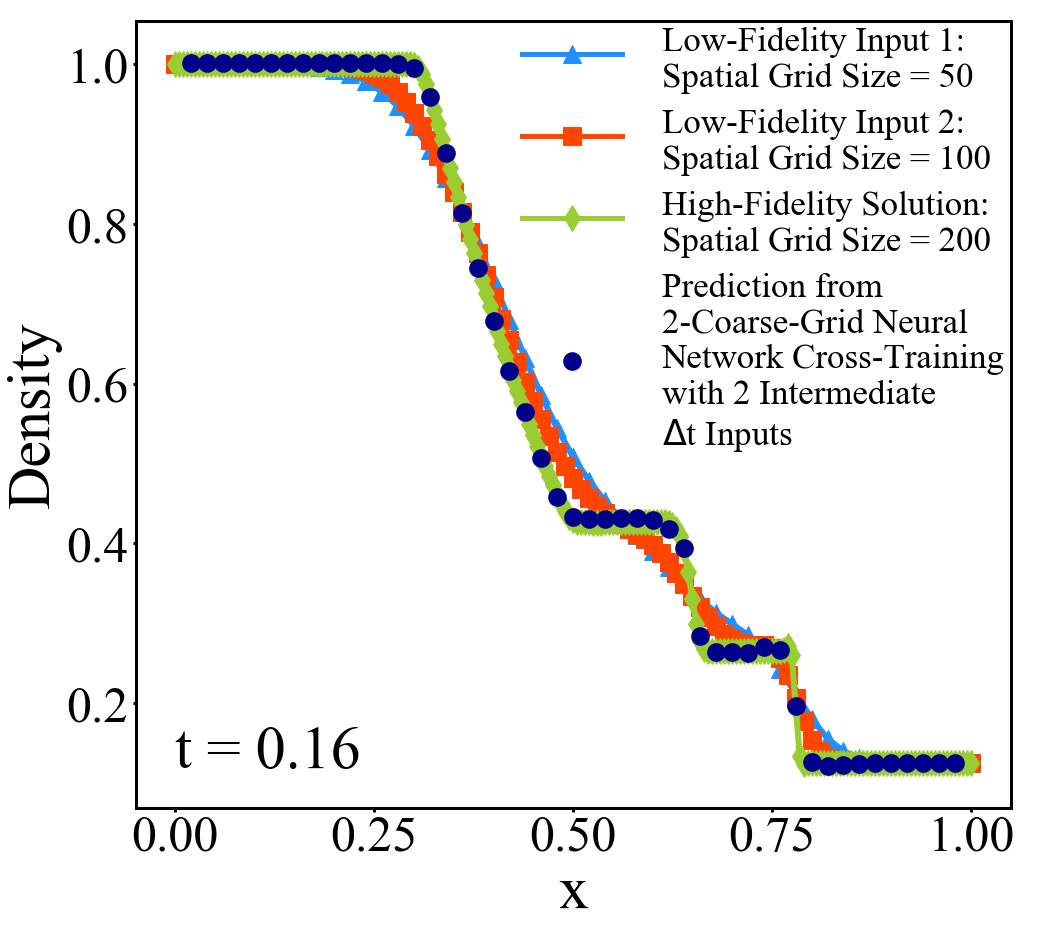}
    \caption{}
\end{subfigure}
\begin{subfigure}[b]{0.327\textwidth}
    \centering
    \includegraphics[width=1.0\linewidth]{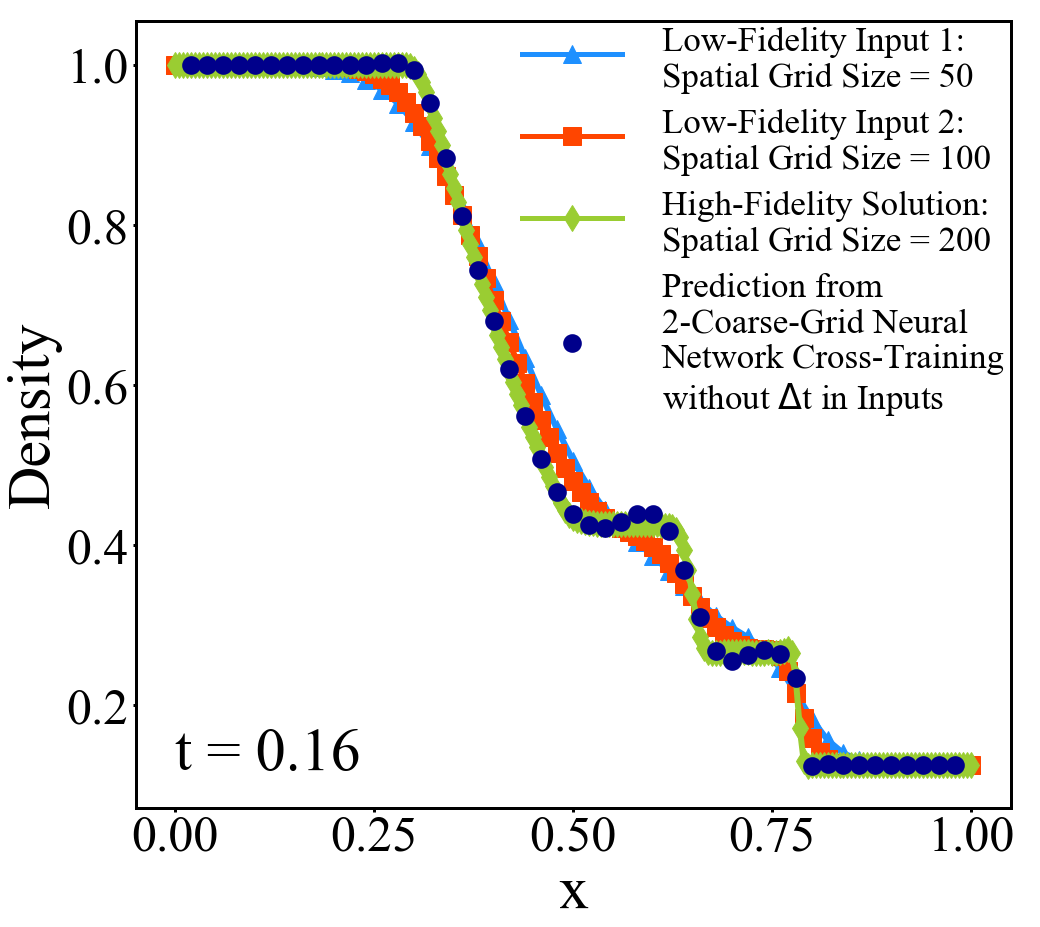}
    \caption{}
\end{subfigure}
\caption{2CGNN Cross-Training predictions of density profiles at final-time $(t=0.16)$ solution of {\bf Sod problem} (dark blue), low-fidelity input solutions (blue and red) by leapfrog and diffusion splitting scheme (\ref{leapfrog-diffusion-splitting}) on $2$ different grids (with $50$ and $100$ cells resp.), and ``exact'' (reference) solution (green): (a) input type~(\ref{cross-training-input1}), (b) input with 2 intermediate time steps, (c) input type~(\ref{standard-2CGNN-input}).}
\label{2CGNN: cross-training, Sod prob.}
\end{figure}

\subsection{2CGNN for the Woodward-Colella Problem}
\label{2CGNN-4-WC}

\begin{figure}[H]\centering
\begin{subfigure}[b]{0.32\textwidth}
    \centering
    \includegraphics[width=1.0\linewidth]{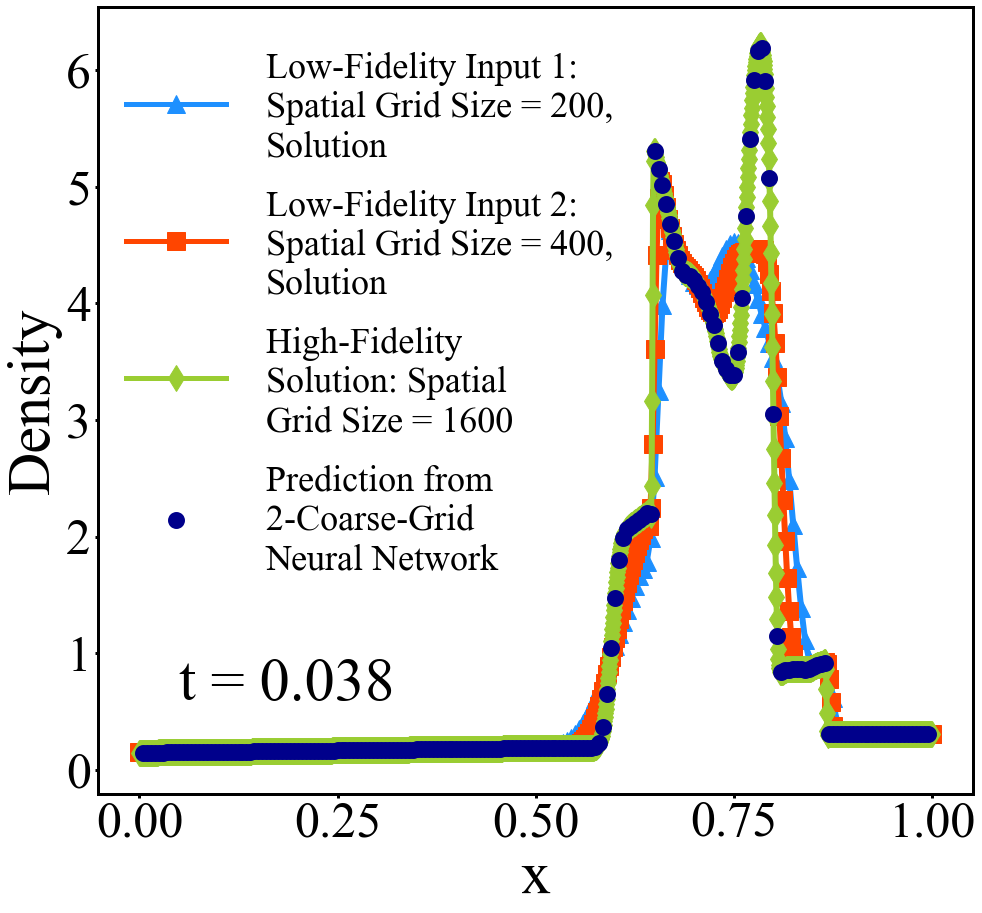}
    \caption{}
\end{subfigure}
\begin{subfigure}[b]{0.335\textwidth}
    \centering
    \includegraphics[width=1.0\linewidth]{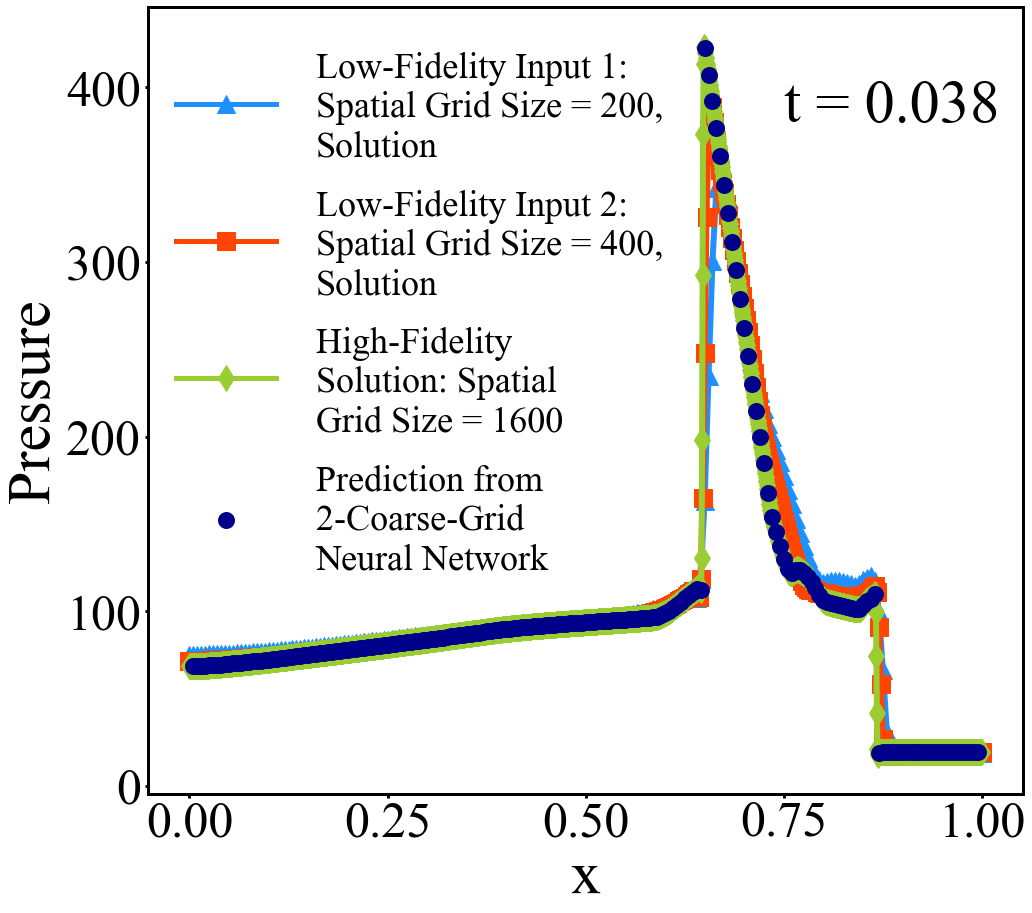}
    \caption{}
\end{subfigure}
\begin{subfigure}[b]{0.327\textwidth}
    \centering
    \includegraphics[width=1.0\linewidth]{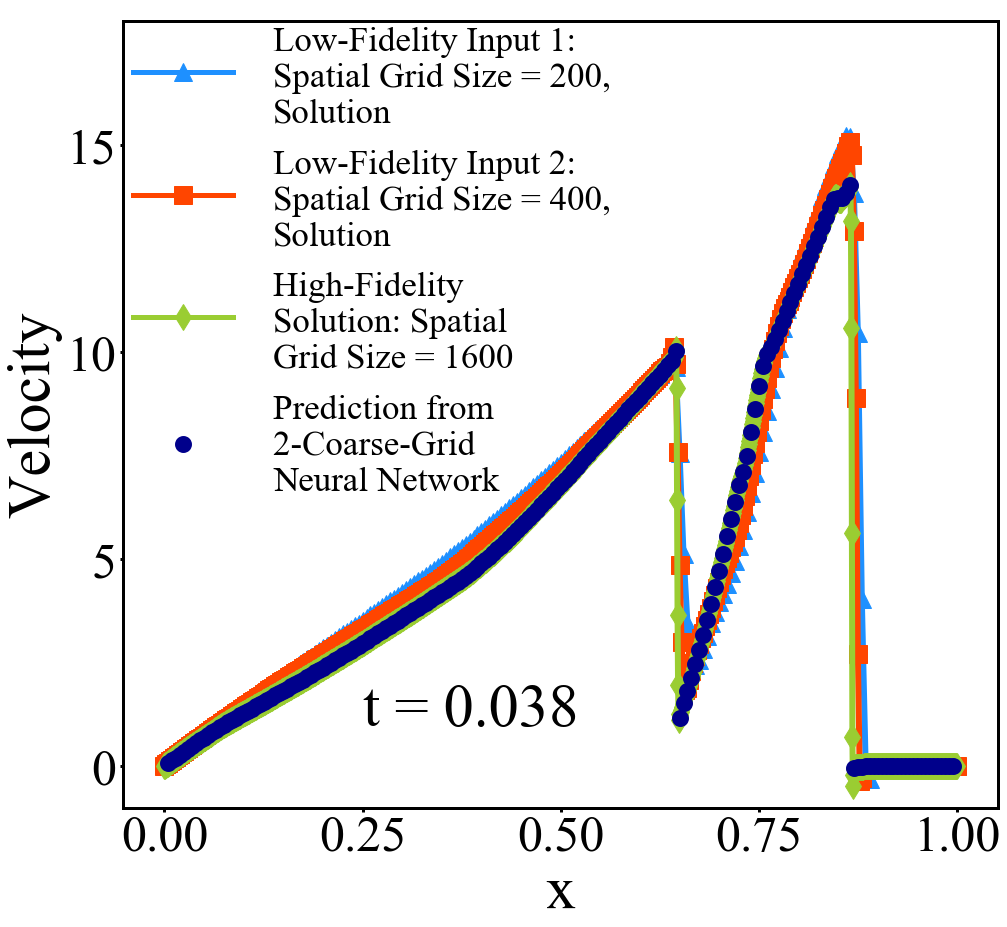}
    \caption{}
\end{subfigure}
\caption{2CGNN prediction of final-time $(t=0.038)$ solution of {\bf Woodward-Colella problem}  (dark blue), low-fidelity input solutions (blue and red) by a 3rd order finite volume method on $2$ different grids (with $200$ and $400$ cells resp.), and ``exact'' (reference) solution (green): (a) density, (b) pressure, (c) velocity.}
\label{2CGNN_src_input_WC_prob_original}
\end{figure}

\begin{figure}[H]\centering
\begin{subfigure}[b]{0.32\textwidth}
    \centering
    \includegraphics[width=1.0\linewidth]{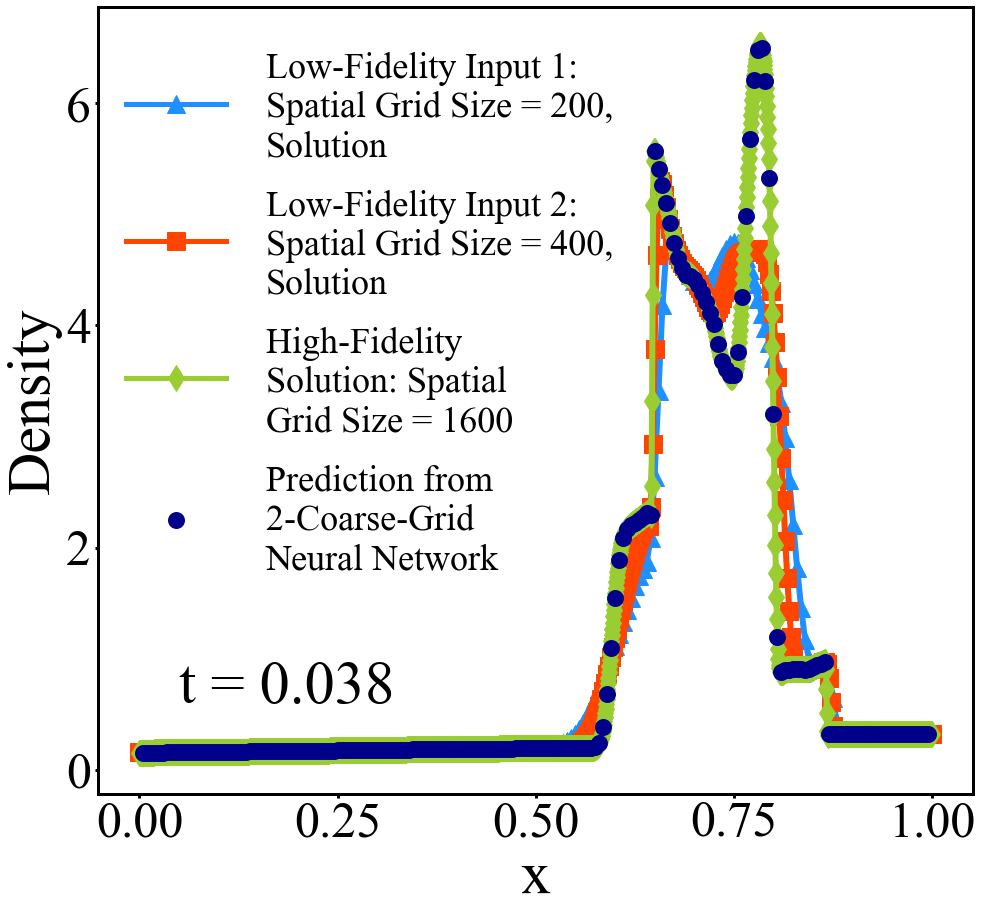}
    \caption{}
\end{subfigure}
\begin{subfigure}[b]{0.335\textwidth}
    \centering
    \includegraphics[width=1.0\linewidth]{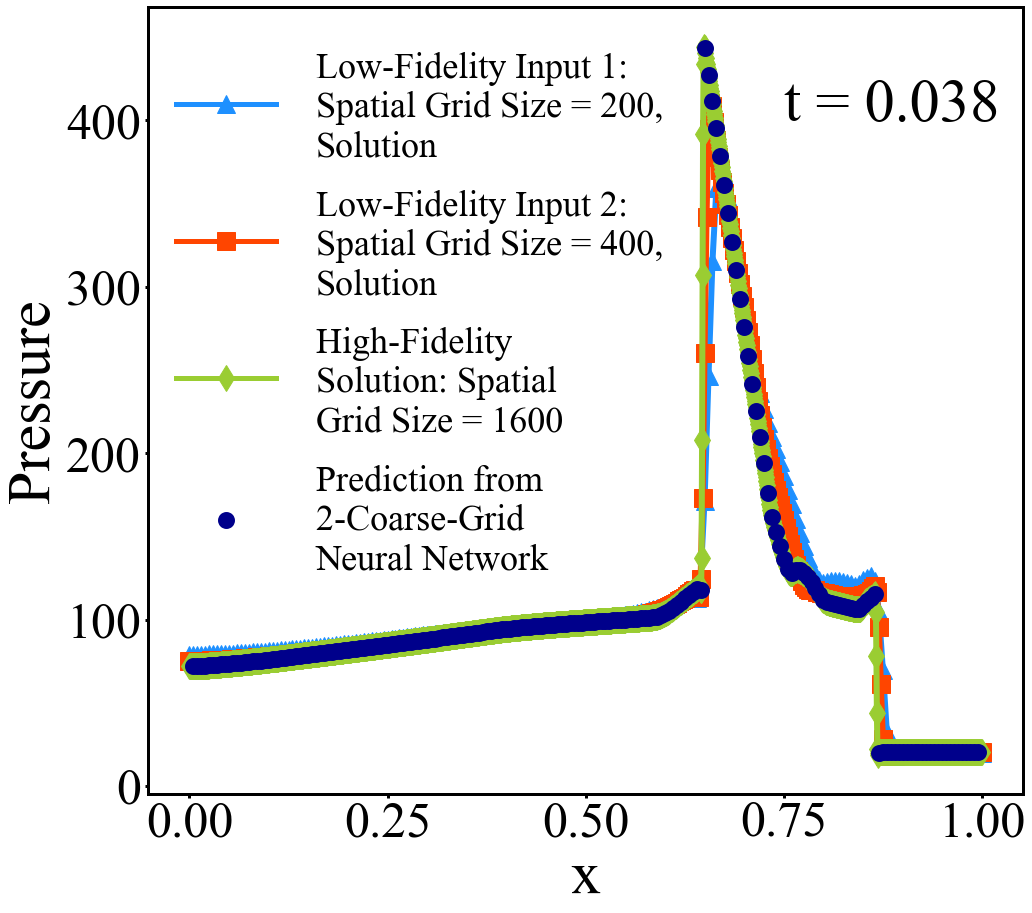}
    \caption{}
\end{subfigure}
\begin{subfigure}[b]{0.327\textwidth}
    \centering
    \includegraphics[width=1.0\linewidth]{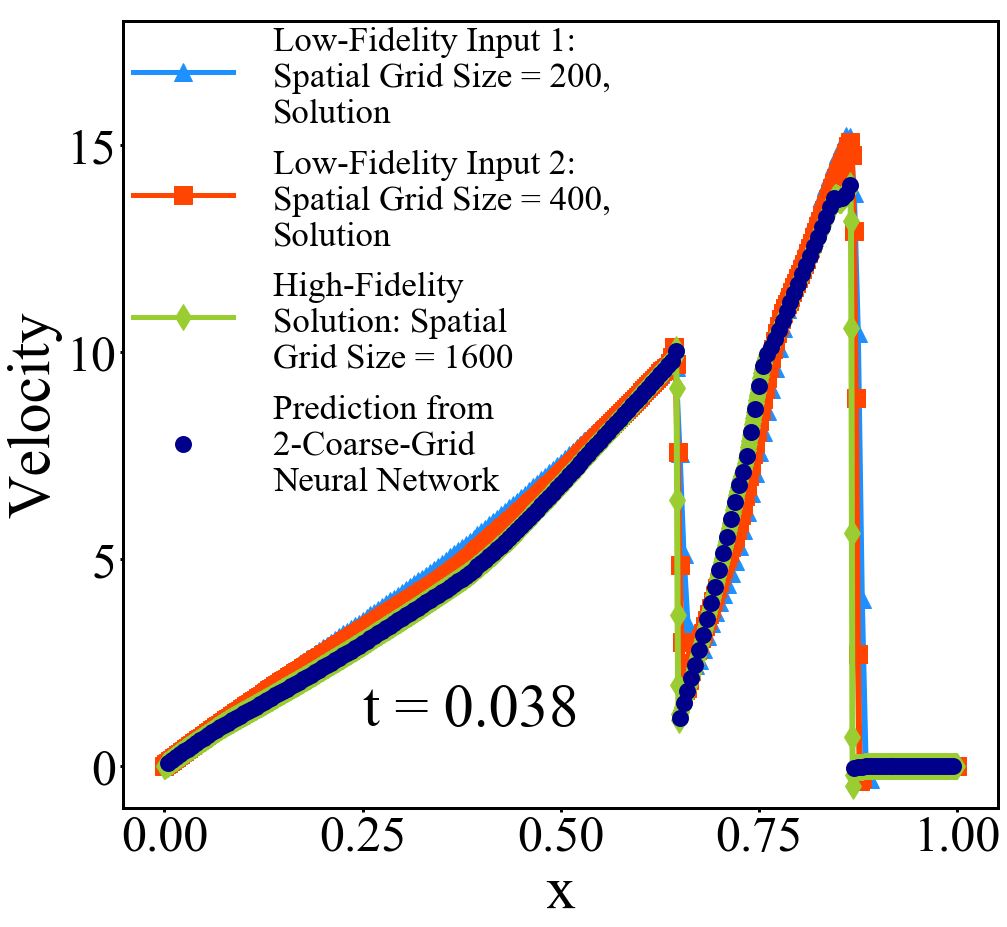}
    \caption{}
\end{subfigure}
\caption{2CGNN prediction of final-time $(t=0.038)$ solution (dark blue) of the blast waves with its {\bf initial value being $+5\%$ perturbation of that of the Woodward-Colella problem}, low-fidelity input solutions (blue and red) by a 3rd order finite volume method on $2$ different grids (with $200$ and $400$ cells resp.), and ``exact'' (reference) solution (green): (a) density, (b) pressure, (c) velocity.}
\label{2CGNN_src_input_WC_prob_p5}
\end{figure}

In \cite{2CGNN_1D_RiemannProb_21}, the low-cost schemes used
are all first order schemes. These schemes do not work well for computing inputs for 2CGNN for the Woodward-Colella (W-C) problem~\cite{WCProb84}, because on reasonable grid sizes these inputs are too qualitatively different from the solution. We therefore use the 3rd order finite volume scheme for computing reference solutions to compute inputs on grids with $200$ and $400$ cells (for 2CGNN with 
input type (\ref{standard-2CGNN-input}).)
Figures \ref{2CGNN_src_input_WC_prob_original} and \ref{2CGNN_src_input_WC_prob_p5} show the interactive blast waves of the W-C problem predicted by 2CGNN. It is clear that the predicted solution improves greatly on the input ones. The training of the neural network is similar to that in the previous subsection. 

\subsection{A Variant of 2DCNN for the Woodward-Colella Problem}
\label{low-high-2DCNN-1D}
\begin{figure}[H]\centering
\begin{subfigure}[b]{0.32\textwidth}
    \centering
    \includegraphics[width=1.0\linewidth]{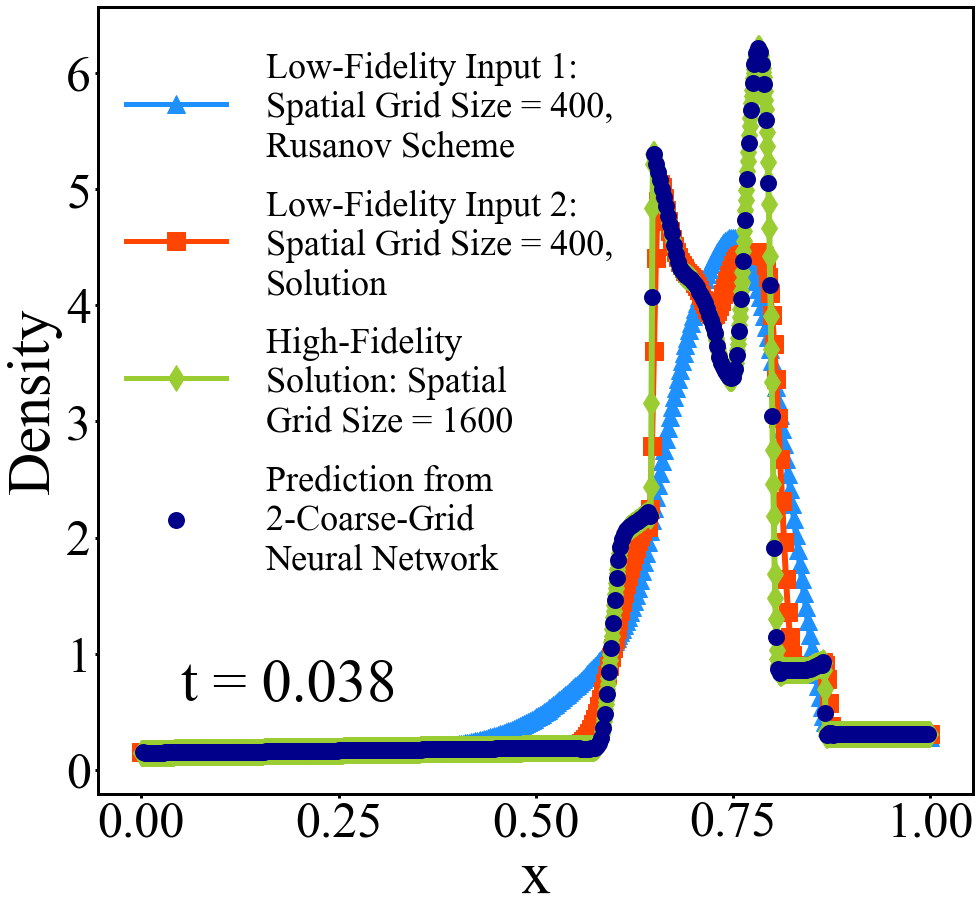}
    \caption{}
\end{subfigure}
\begin{subfigure}[b]{0.335\textwidth}
    \centering
    \includegraphics[width=1.0\linewidth]{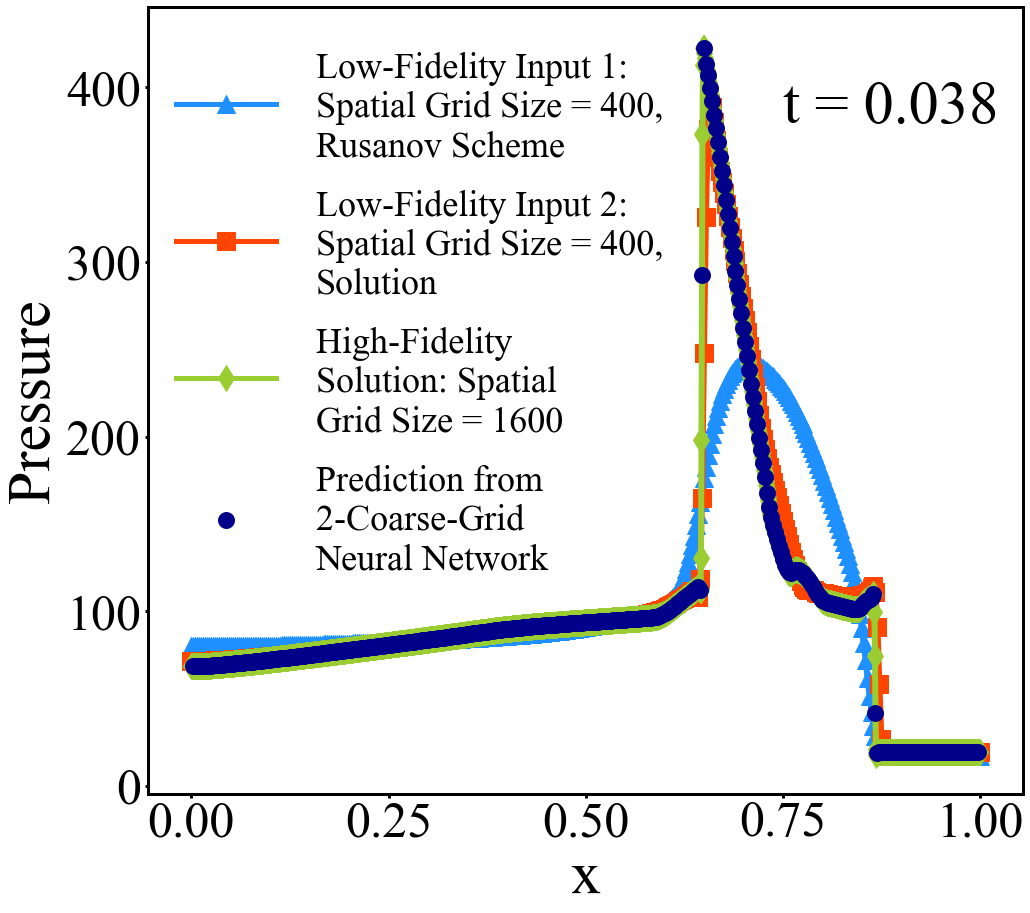}
    \caption{}
\end{subfigure}
\begin{subfigure}[b]{0.327\textwidth}
    \centering
    \includegraphics[width=1.0\linewidth]{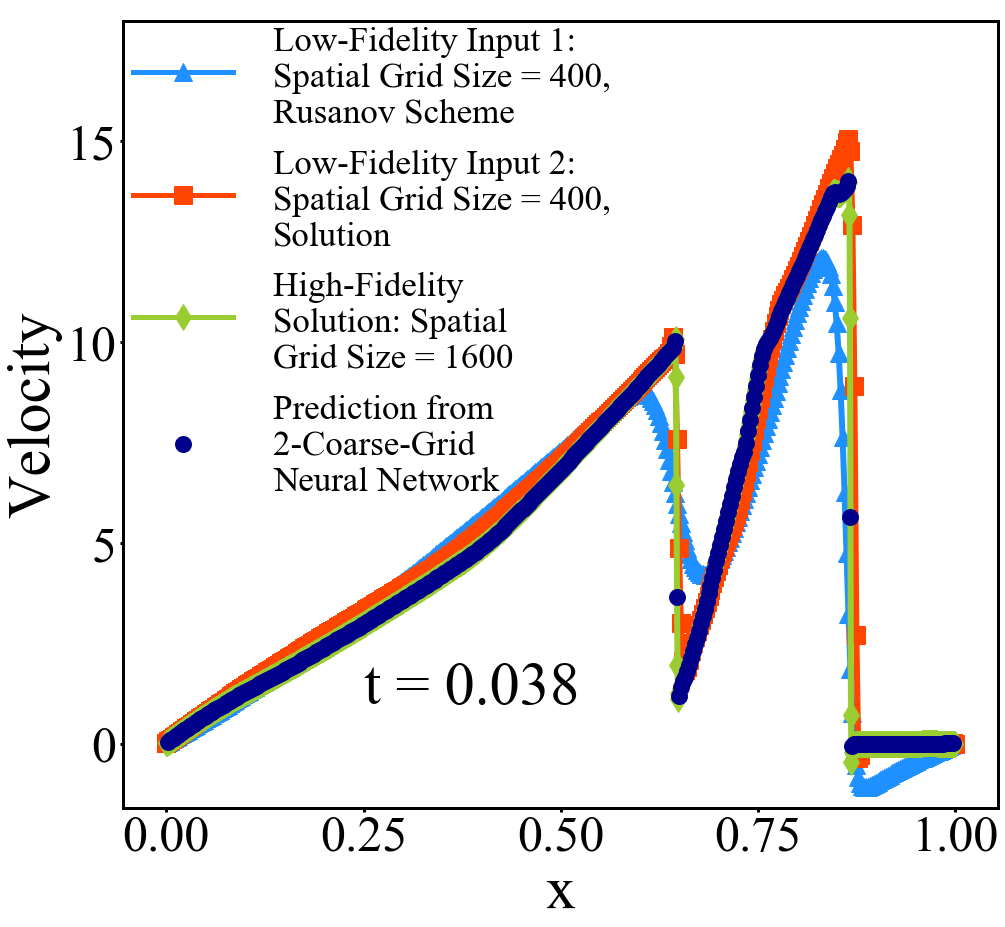}
    \caption{}
\end{subfigure}
\caption{2DCNN prediction of final-time $(t=0.038)$ solution of {\bf Woodward-Colella problem}  (dark blue), low-fidelity input solutions (blue and red) by the Rusanov scheme and a 3rd order finite volume method respectively on the same grid (with $400$ cells), and ``exact'' (reference) solution (green): (a) density, (b) pressure, (c) velocity.}
\label{2DCNN_src_input_rusanov_WC_prob_original}
\end{figure}

\begin{figure}[H]\centering
\begin{subfigure}[b]{0.32\textwidth}
    \centering
    \includegraphics[width=1.0\linewidth]{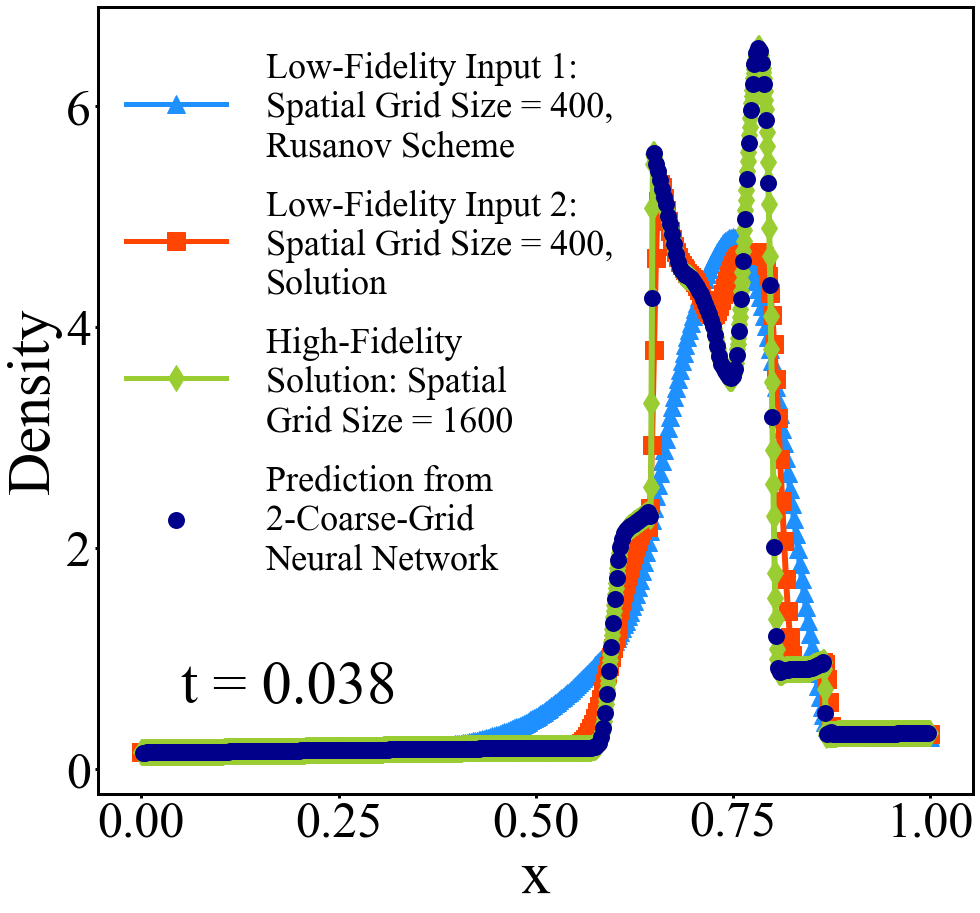}
    \caption{}
\end{subfigure}
\begin{subfigure}[b]{0.335\textwidth}
    \centering
    \includegraphics[width=1.0\linewidth]{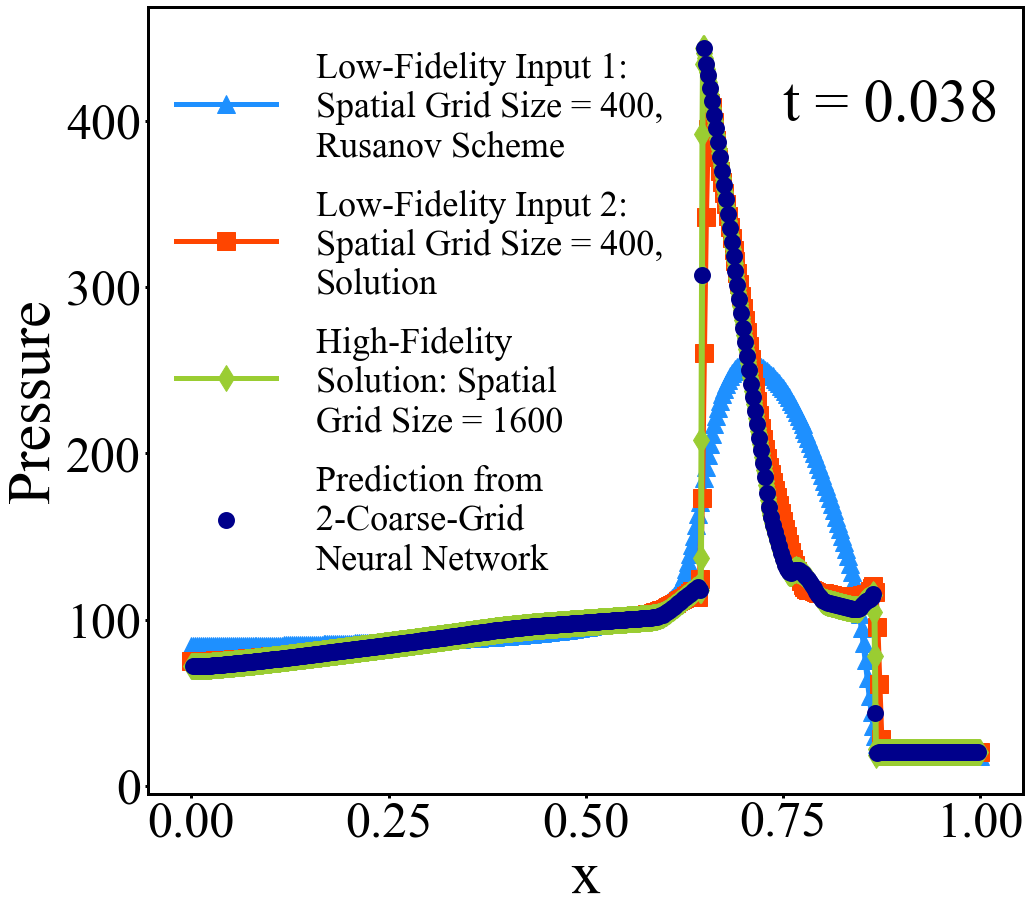}
    \caption{}
\end{subfigure}
\begin{subfigure}[b]{0.327\textwidth}
    \centering
    \includegraphics[width=1.0\linewidth]{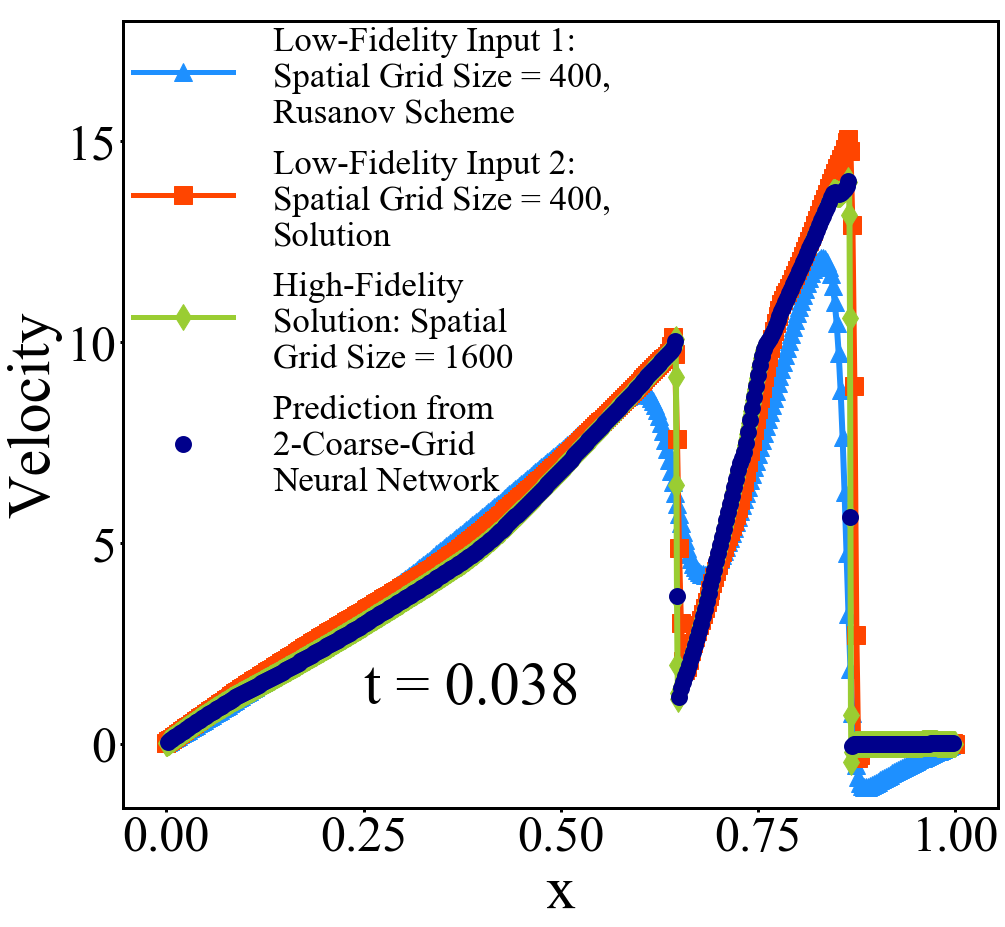}
    \caption{}
\end{subfigure}
\caption{2DCNN prediction of final-time $(t=0.038)$ solution (dark blue) of the blast waves with its {\bf initial value being $+5\%$ perturbation of that of the Woodward-Colella problem}, low-fidelity input solutions (blue and red) by the Rusanov scheme and a 3rd order finite volume method respectively on the same grid (with $400$ cells), and ``exact'' (reference) solution (green): (a) density, (b) pressure, (c) velocity.}
\label{2DCNN_src_input_rusanov_WC_prob_p5}
\end{figure}

In~\cite{2CGNN_1D_RiemannProb_21}, in addition to 2CGNN, we proposed the 2-Diffusion-Coefficient Neural Network (2DCNN) for solving 1-D Riemann problems. The low-cost scheme used is a first order scheme (leapfrog and diffusion splitting scheme), with different diffusion coefficients for computing low-fidelity solutions as inputs for the Lax and Sod problems. However, first order schemes do not work well for computing inputs for 2CGNN for the Woodward-Colella problem, as mentioned in Sec.~\ref{2CGNN-4-WC}. Therefore, we use the first order Rusanov scheme and the 3rd order finite volume scheme for reference solutions to compute inputs on the same coarse grid. This process can be viewed as a variant of the 2DCNN from \cite{2CGNN_1D_RiemannProb_21} because only one grid is used in computing inputs. The input format is
similar to (\ref{standard-2CGNN-input}) with the first part of the input computed by the Rusanov scheme~\cite{Rusanov61} at corresponding space-time locations, and the second part of the input computed by the higher order scheme at the same space-time locations. 
Figures \ref{2DCNN_src_input_rusanov_WC_prob_original} and \ref{2DCNN_src_input_rusanov_WC_prob_p5} show the interactive blast waves of the W-C problem predicted by 2DCNN. It is clear that the predicted solution improves greatly over the input ones. The training of the neural network is similar to that in the previous subsection. 

\section{2CGNN for 2D Riemann Problems}
\label{2CGNN-4-2D}
Consider the 2D scalar conservation law
\begin{equation}
\label{cons-law}
   \frac{\partial U}{\partial t} + \frac{\partial f(U)}{\partial x} + \frac{\partial g(U)}{\partial y} = 0,  (x, y)\in\Omega\subset \mathcal{R},  t\in[0, T],
\end{equation}
\\
and the 2D Euler equations for the ideal gas 

\begin{equation}
\label{2D Riemann Problem formula}
    \frac{\partial}{\partial t}\left(\begin{array}{c} \rho \\ \rho u \\ \rho v \\E \end{array}\right) +\frac{\partial }{\partial x} \left(\begin{array}{c} \rho u\\ \rho u^2 + p \\ \rho u v \\ u( E + p)\end{array}\right) +\frac{\partial }{\partial y} \left(\begin{array}{c} \rho v\\ \rho v u \\ \rho v^2 + p \\ v( E + p)\end{array}\right) = 0,\; (x, y)\in\Omega\subset \mathcal{R}, \; t\in[0, T]~,
\end{equation}
where $\rho$, $u$, $v$ and $p$ are density, $x$ and $y$ components of velocity, and pressure, respectively, 
\begin{equation}
\label{energy}
    E = \frac{p}{\gamma - 1}+\frac{1}{2}\rho (u^2 + v^2)~,
\end{equation}
and $\gamma=1.4$.

 In this section, we extend the application of 2CGNN from one dimensional problems to two dimensional problems. 

\subsection{Input, Output and Loss Function}
Let  $\Omega =[a,b]\times [c,d]$ be partitioned with the coarsest uniform rectangular grid $a=x_0<x_1<\dots <x_M=b$ and $c=y_0<y_1<\dots <y_N=d$. The spatial grid size is $\Delta x=x_1-x_0$ and $\Delta y=y_1-y_0$, and the time step size is $\Delta t$. We refine the grid to obtain a finer uniform grid with spatial grid size $\frac12 \Delta x$ and $\frac12 \Delta y$, and time step size $\frac12 \Delta t$. Let $L$ be
a low-cost scheme used to compute (\ref{cons-law}) on both grids.
Suppose we are interested in predicting the solution at $(x,y,t)$, that is referred as $(i',j',n')$ in the coarsest grid, we choose the coarsest grid solution (computed by $L$) at $9$ points $(x_{i'-1},y_{j'-1})$, $(x_{i'-1},y_{j'})$, $(x_{i'-1},y_{j'+1})$, $(x_{i'},y_{j'-1})$, $(x_{i'},y_{j'})$, $(x_{i'},y_{j'+1})$, $(x_{i'+1},y_{j'-1})$, $(x_{i'+1},y_{j'})$, and $(x_{i'+1},y_{j'+1})$ at time level $t_{n'-1}$, and also at point $(x_{i'}, y_{i'}, t_{n'})$ as the first part of the input of the neural network, and the finer grid solution (also computed by $L$) at the same space-time locations as the second part of input. Note that, the chosen $10$ space-time locations on both grids enclose a local (space-time) domain of dependence of the exact solution at $(x_{i'}, y_{i'}, t_{n'})$ (with $\Delta t$ satisfying the CFL condition.)

\begin{figure}[H] \centering
\begin{subfigure}[b]{.8\textwidth}
    \centering
    \includegraphics[width=1.0\linewidth]{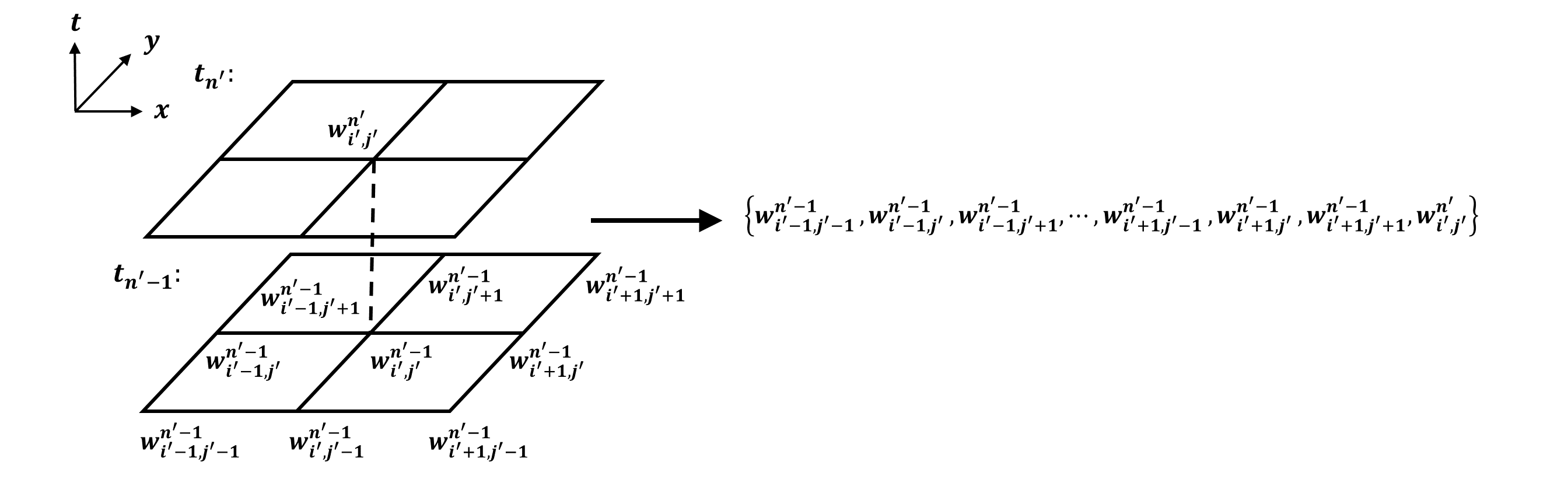}
    \caption{}
\end{subfigure}
\vskip10pt

\begin{subfigure}[b]{.8\textwidth}
    \centering
    \includegraphics[width=1.0\linewidth]{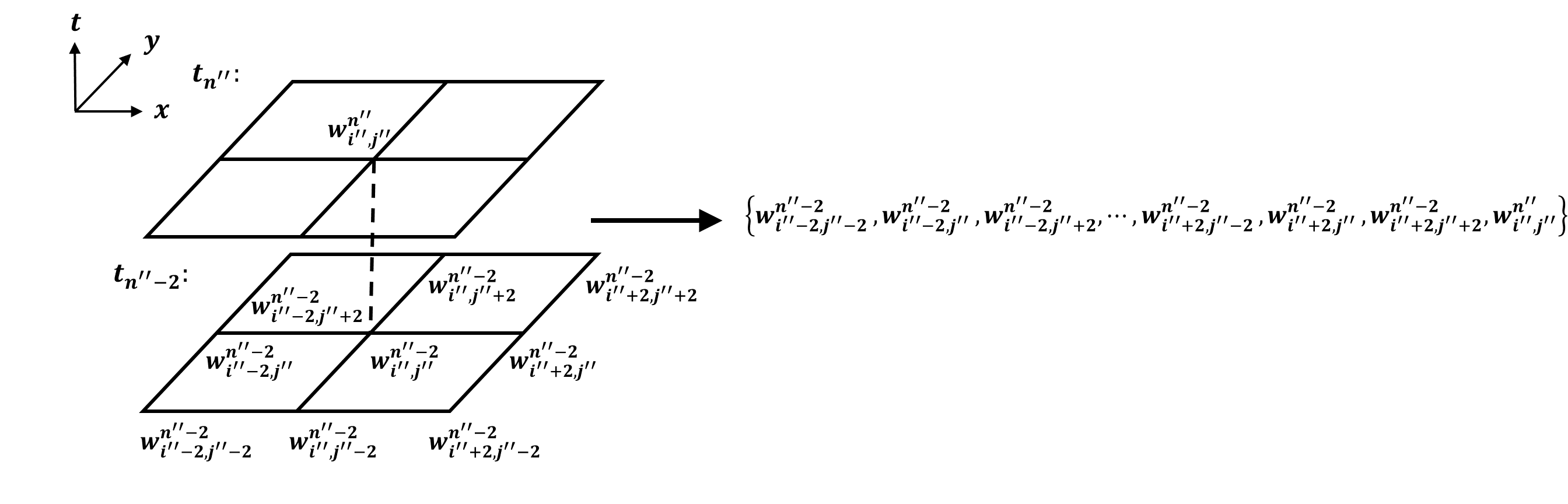}
    \caption{}
\end{subfigure}
\vskip10pt

\begin{subfigure}[b]{.8\textwidth}
    \centering
    \includegraphics[width=1.0\linewidth]{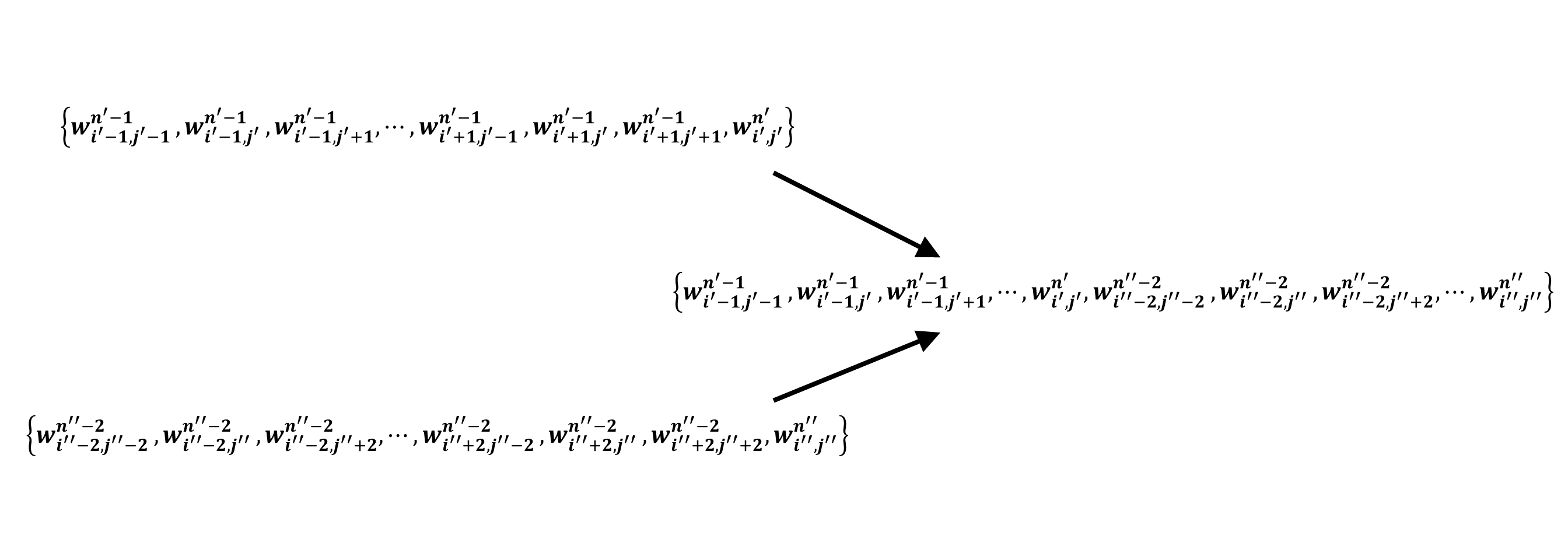}
    \caption{}
\end{subfigure}
\caption{Procedure for formatting the input for 2CGNN.}
\label{Input_2CGNN_2D}
\end{figure}

\begin{figure}[H]
\centering
\begin{subfigure}[b]{.8\textwidth}
    \centering
    \includegraphics[width=1.0\linewidth]{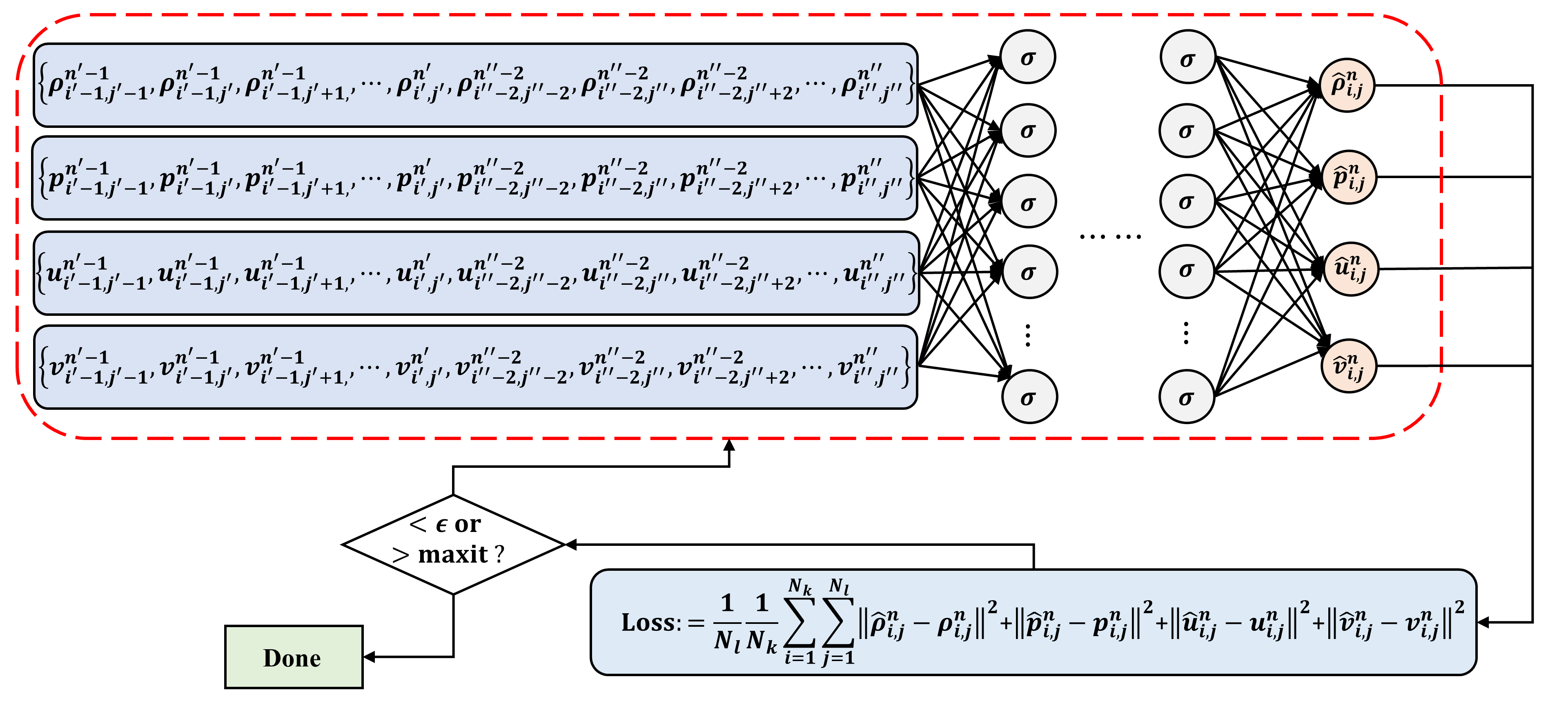}
\end{subfigure}
\caption{Training procedure for 2CGNN.}
\label{2CGNN_2D}
\end{figure}

Denote the first part of the 2D input as
$$w^{n'-1}_{i'-1,j'-1}, w^{n'-1}_{i'-1,j'}, w^{n'-1}_{i'-1,j'+1}, w^{n'-1}_{i',j'-1}, w^{n'-1}_{i',j'}, w^{n'-1}_{i',j'+1}, w^{n'-1}_{i'+1,j'-1}, w^{n'-1}_{i'+1,j'}, w^{n'-1}_{i'+1,j'+1}, w^{n'}_{i',j'}~,$$
and 
the second part of the 2D input as
$$w^{n''-2}_{i''-2,j''-2}, w^{n''-2}_{i''-2,j''}, w^{n''-2}_{i''-2,j''+2}, w^{n''-2}_{i'',j''-2}, w^{n''-2}_{i'',j''}, w^{n''-2}_{i'',j''+2}, w^{n''-2}_{i''+2,j''-2}, w^{n''-2}_{i''+2,j''}, w^{n''-2}_{i''+2,j''+2}, w^{n''}_{i'',j''}~.$$
Note that the space-time indices $(i',j',n')$ on the coarsest grid refers to the same location as $(i'',j'',n'')$ does on the finer grid, $(i'-1, j'-1, n'-1)$ refers to the same location as $(i''-2, j''-2, n''-2)$ does, and so on.
The input of 2CGNN is now
\begin{equation}
\begin{array}{cc}
    \{w^{n'-1}_{i'-1,j'-1}, w^{n'-1}_{i'-1,j'}, w^{n'-1}_{i'-1,j'+1}, w^{n'-1}_{i',j'-1}, w^{n'-1}_{i',j'}, w^{n'-1}_{i',j'+1}, w^{n'-1}_{i'+1,j'-1}, w^{n'-1}_{i'+1,j'}, w^{n'-1}_{i'+1,j'+1}, w^{n'}_{i',j'},\\
    w^{n''-2}_{i''-2,j''-2}, w^{n''-2}_{i''-2,j''}, w^{n''-2}_{i''-2,j''+2}, w^{n''-2}_{i'',j''-2}, w^{n''-2}_{i'',j''}, w^{n''-2}_{i'',j''+2},\\
     w^{n''-2}_{i''+2,j''-2}, w^{n''-2}_{i''+2,j''}, w^{n''-2}_{i''+2,j''+2}, w^{n''}_{i'',j''}\}~,
\end{array}
\label{standard-input}
\end{equation}
called the ``input of $w$,'' and the corresponding output of 2CGNN
is the predicted solution of (\ref{cons-law}) at $(x,y,t)$ (or $(i',j',n')$ on the coarsest grid.)
See Fig.~\ref{Input_2CGNN_2D} and \ref{2CGNN_2D} for an illustration of 2CGNN.
For the Euler system, the input and output of 2CGNN are made up of corresponding ones for each prime variable. For example, the input can be
the vector 
$$\{ {\rm input}\; {\rm of}\; \rho,\;
{\rm input}\; {\rm of}\; u,\;
{\rm input}\; {\rm of}\; v,\;
{\rm input}\; {\rm of}\; p\}
$$
with $20\times 4=80$ elements,
and the corresponding output will be 
$\{\rho, u, v, p\}$ at $(x,y,t)$ (or $(i',j',n')$ on the coarsest grid) with $4$ elements.

The loss function measures the difference between the output and the reference solution corresponding to the input, and is defined as follows.
$$
\begin{array}{ccc}
{\rm Loss}&=&\sum_{k}\|({\rm output} \; {\rm corresponding}\; {\rm to}\; {\rm the}\; k^{th}\;  {\rm set}\; {\rm of}\; {\rm input})- \\
&& ({\rm reference}\; {\rm solution}\; {\rm corresponding}\; {\rm to}\; {\rm the}\; k^{th}\; {\rm set}\; {\rm of}\; {\rm input}) \|_2^2~,
\end{array}
$$

where the summation goes through every set of input in the training data.

\subsection{Generation of Input and Training Data}
We will use 2CGNN to predict solutions of Riemann problems from \cite{PDLax98}.
The associated neural network for 2CGNN varies for different configurations. For Configurations 1, 2, and 3, the corresponding neural network consists of $8$ hidden layers, each with $320$ neurons. That for Configurations 6 and 8 consists of $8$ hidden layers, and Configuration 4 needs a neural network with $9$ hidden layers; each of these layers has $360$ neurons. During the training process, the neural network minimizes the difference between outputs and a reference solution by using first an Adam optimizer then an L-BFGS optimizer in TensorFLow (with the number of iterations per optimization procedure under $50000$ each.) After the training is done, the neural network is used to predict a solution (different from the training data), given an input computed by the same low-cost scheme(s) and grids that are used to compute the inputs of the training data. We use the trained neural network to predict final solutions for $5$ initial values of the Euler system for each configuration, including the original initial value of the configuration, and $\pm3\%$ and $\pm5\%$ perturbations of the initial value. 

In order to generate the training data, we use a first order scheme on the coarsest uniform grid ($200$ cells in each spatial dimension) and the finer uniform grid ($400$ cells each spatial dimension) to compute the input data from several initial values of the Euler system for each configuration, including $\pm2\%$, $\pm4\%$, $\pm6\%$, $\pm8\%$ and $\pm10\%$ perturbations of the initial value of the configuration.
The high resolution reference solutions of the training data are computed on a uniform grid
with $400$ cells in each spatial dimension (note that solution values at grid points need to be interpolated from cell averages.) A 4th order central scheme on overlapping cells with HR limiting \cite{LiuShu07b} is used for computing reference solutions for all 2D Riemann problems. The low-cost scheme for (\ref{cons-law}) used for computing inputs is the first order leapfrog and diffusion splitting scheme 
\begin{equation}
\label{leapfrog-diffusion-splitting, 2D}
\left \{
\begin{array}{l}
    \frac{\tilde{U}_{i,j} - U^{n-1}_{i,j}}{2\Delta t} +\frac{{f(U)}|^{n}_{i+1,j}-{f(U)}|^{n}_{i-1,j}}{2\Delta x} +\frac{{g(U)}|^{n}_{i,j+1}-{g(U)}|^{n}_{i,j-1}}{2\Delta y} =0~, \\
    \frac{U^{n+1}_{i,j} - \tilde{U}_{i,j}}{\Delta t}- \alpha [\frac{\tilde{U}_{i+1,j} - 2\cdot \tilde{U}_{i,j} + \tilde{U}_{i-1,j}}{\Delta x^2} + \frac{\tilde{U}_{i, j+1} - 2\cdot \tilde{U}_{i,j} + \tilde{U}_{i,j-1}}{\Delta y^2} ] =0~,
\end{array}
\right .
\end{equation}
where $\alpha = \Delta x$, and $\Delta x=\Delta y$ throughout the paper.
The time step size for the first order scheme is fixed during the evolution in time, e.g., $\Delta t$ for the $200\times 200$ grid and $\frac12 \Delta t$ for the $400\times 400$ grid, satisfying the CFL condition.

\subsection{Results and Discussion for 2D Riemann Problems}
\label{2d-2CGNN-leapfrog}
We use 2CGNN with inputs computed by the first order 
leapfrog and diffusion splitting scheme (\ref{leapfrog-diffusion-splitting, 2D}) (applied component-by-component for systems) to predict most of the problems. The predictions of the final-time solutions of the 2D Riemann problems are shown in contour plots below. The predicted solution profiles of cross sections perpendicular to the $y$-axis at the final time capture shocks and contacts, as well as the smooth regions of the solution, very well.  The spatial computational domain is $[0,1]\times [0,1]$ unless otherwise specified.  

Fig.~\ref{2CGNN: Final time of config. 6, original} shows the predicted final-time solution of the 2D Euler system in Configuration 6 from \cite{PDLax98}. The corresponding density cross sections perpendicular to the $y$-axis of Configuration 6 at $y = 0.34$, $y = 0.50$, $y = 0.60$, $y = 0.70$ and $y = 0.80$ are shown in Fig.~\ref{2CGNN: Final time cross-section of config. 6, original}.  
Fig.~\ref{2CGNN: Final time of config. 6, +5} shows the prediction of the final-time solution of the 2-D Euler system, with the initial value being $+5\%$ perturbation of that of Configuration 6, and Fig.~\ref{2CGNN: Final time cross-section of config. 6, +5} shows the density cross section profiles perpendicular to the y-axis at the locations of the original initial values case.

Since the coarsest grid for input has $200$ grid cells, the spatial grid for the predicted solution also has $200$ grid cells. The predictions do not depict smeared solutions like their low-cost input solutions. Note that the configuration cases for predictions are not included in the training data.

\begin{figure}[H] \centering
\begin{subfigure}[b]{.48\textwidth}
    \centering
    \includegraphics[width=1.0\linewidth]{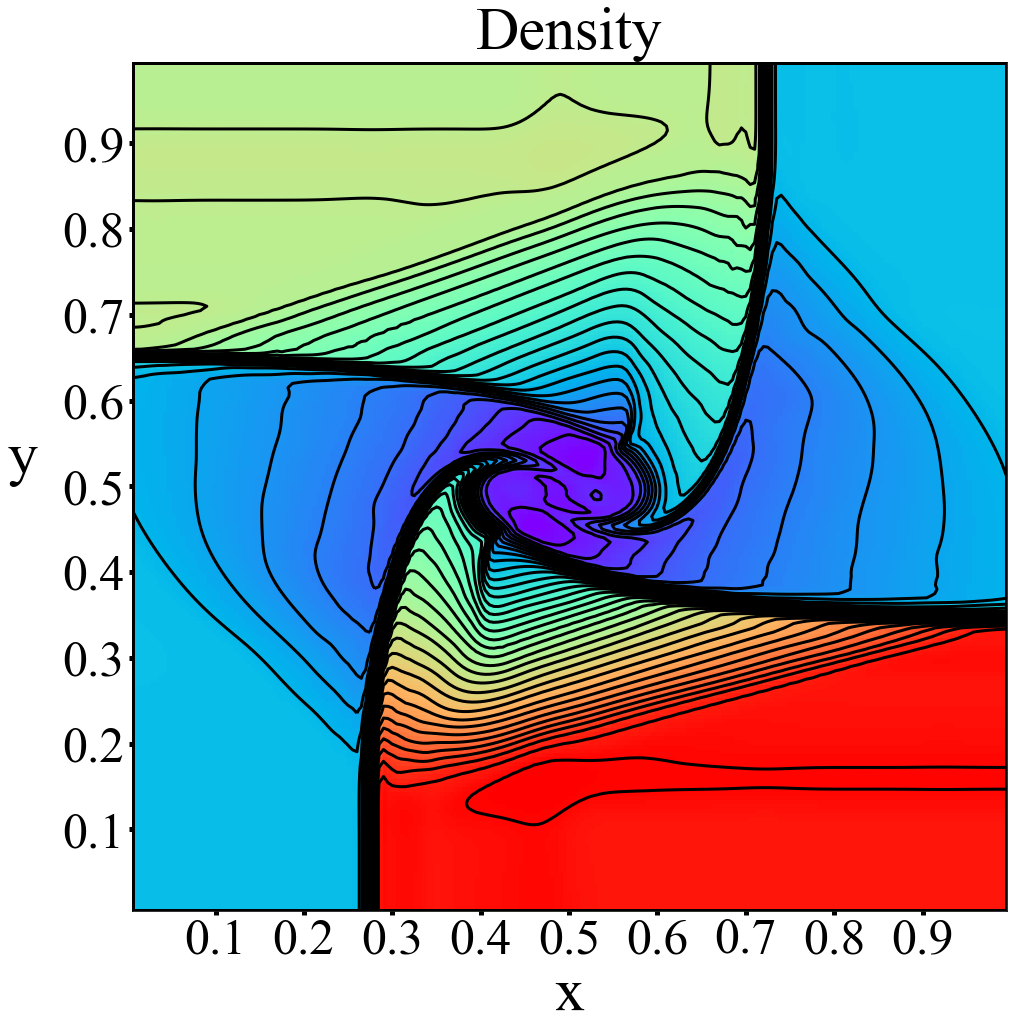}
    \caption{}
\end{subfigure}
\begin{subfigure}[b]{.48\textwidth}
    \centering
    \includegraphics[width=1.0\linewidth]{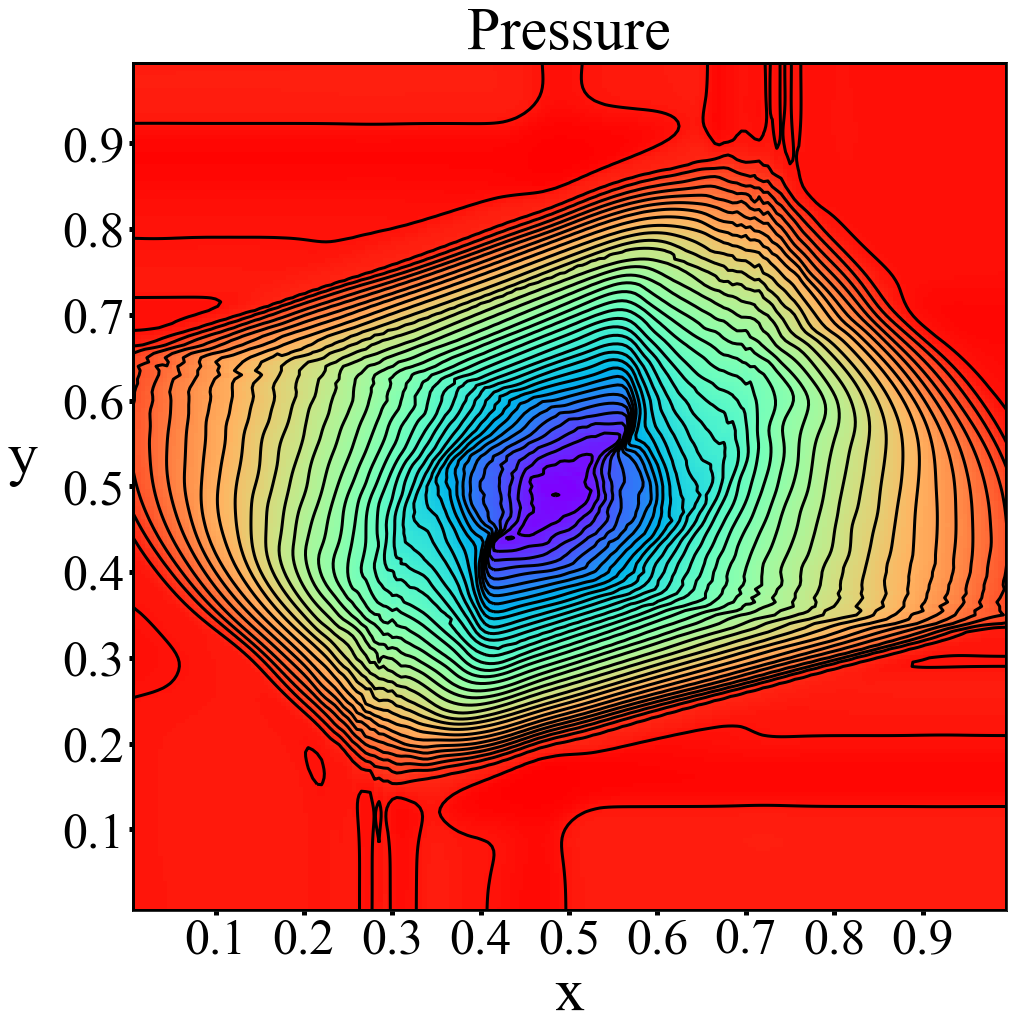}
    \caption{}
\end{subfigure}
\end{figure}

\begin{figure}[H]\ContinuedFloat
\centering
\begin{subfigure}[b]{.48\textwidth}
    \centering
    \includegraphics[width=1.0\linewidth]{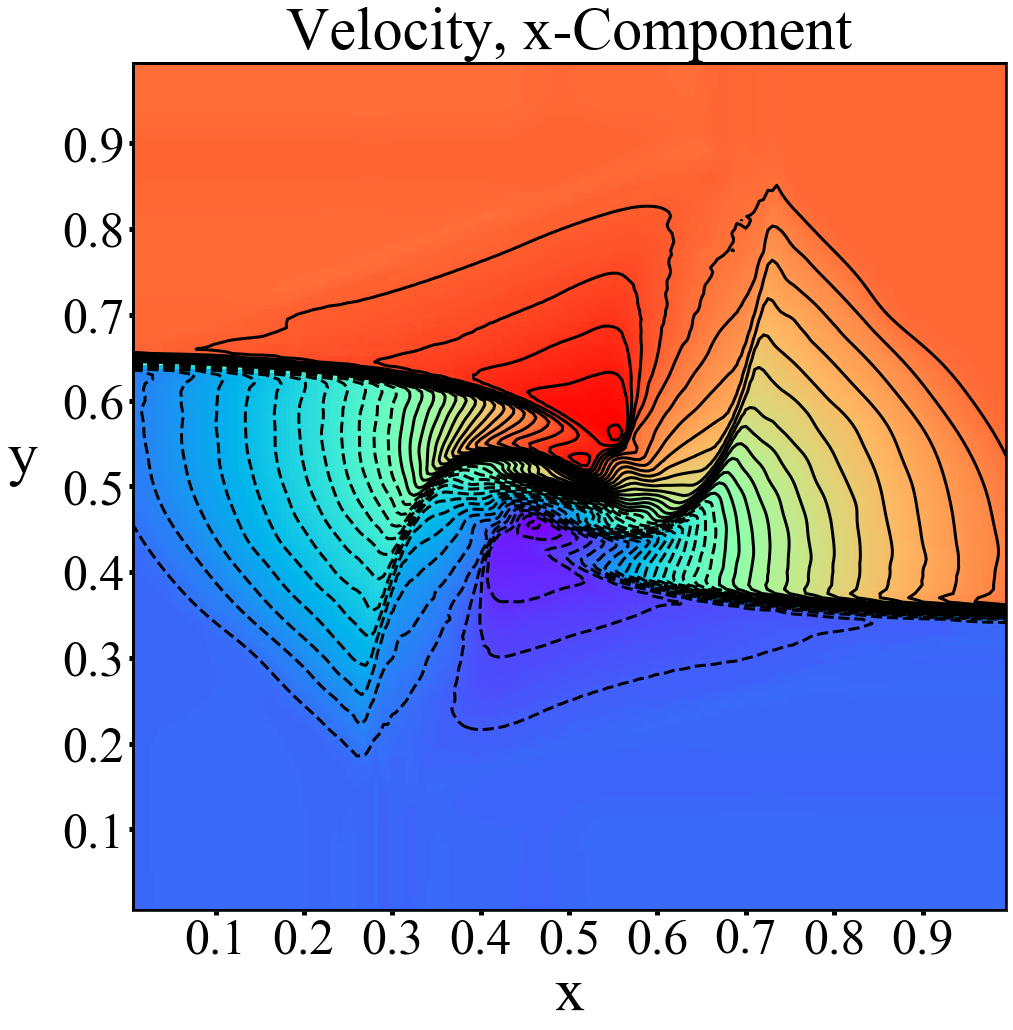}
    \caption{}
\end{subfigure}
\begin{subfigure}[b]{.48\textwidth}
    \centering
    \includegraphics[width=1.0\linewidth]{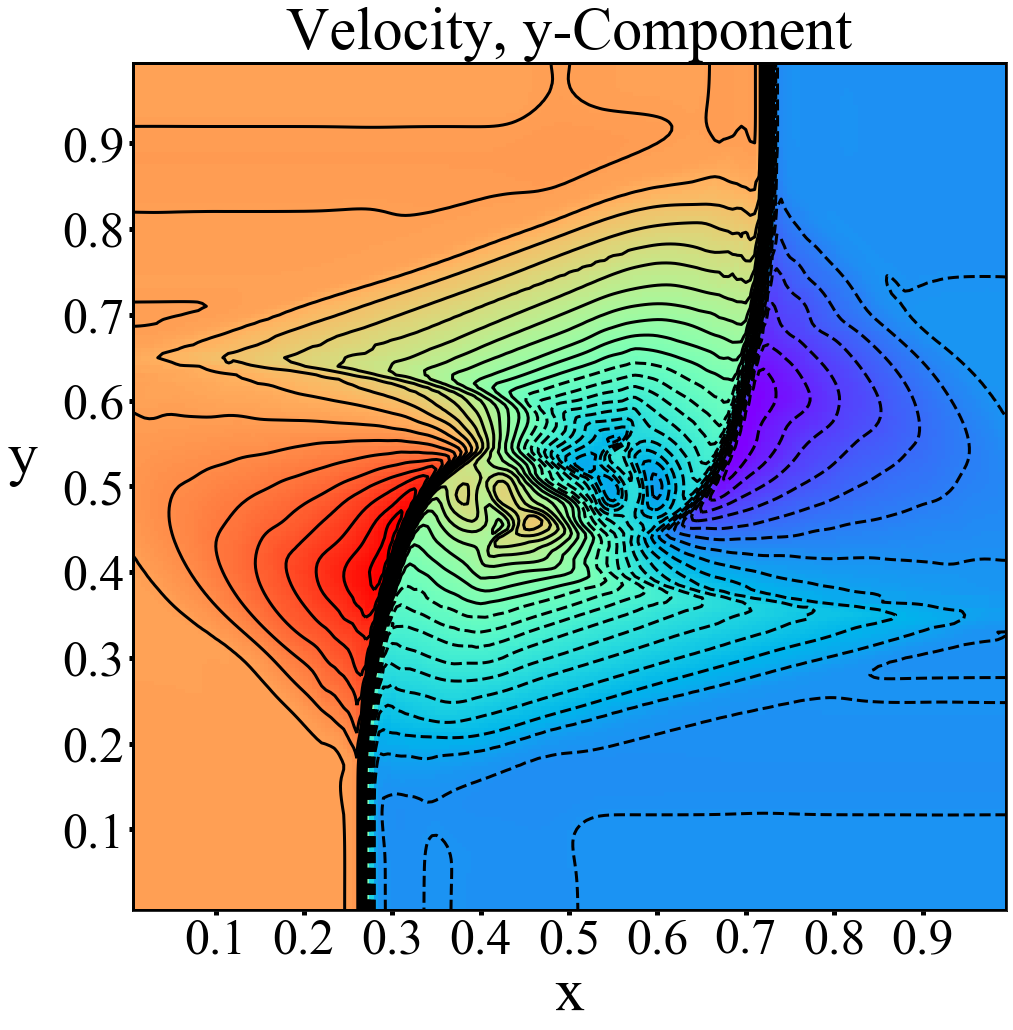}
    \caption{}
\end{subfigure}
\caption{2CGNN prediction of the final-time $(t=0.3)$ solution of {\bf Configuration 6} in \cite{PDLax98}: (a) density, (b) pressure, (c) velocity, x-component, (d) velocity, y-component.}
\label{2CGNN: Final time of config. 6, original}
\end{figure}

\begin{figure}[H]\centering
\begin{subfigure}[b]{0.327\textwidth}
    \centering
    \includegraphics[width=1.0\linewidth]{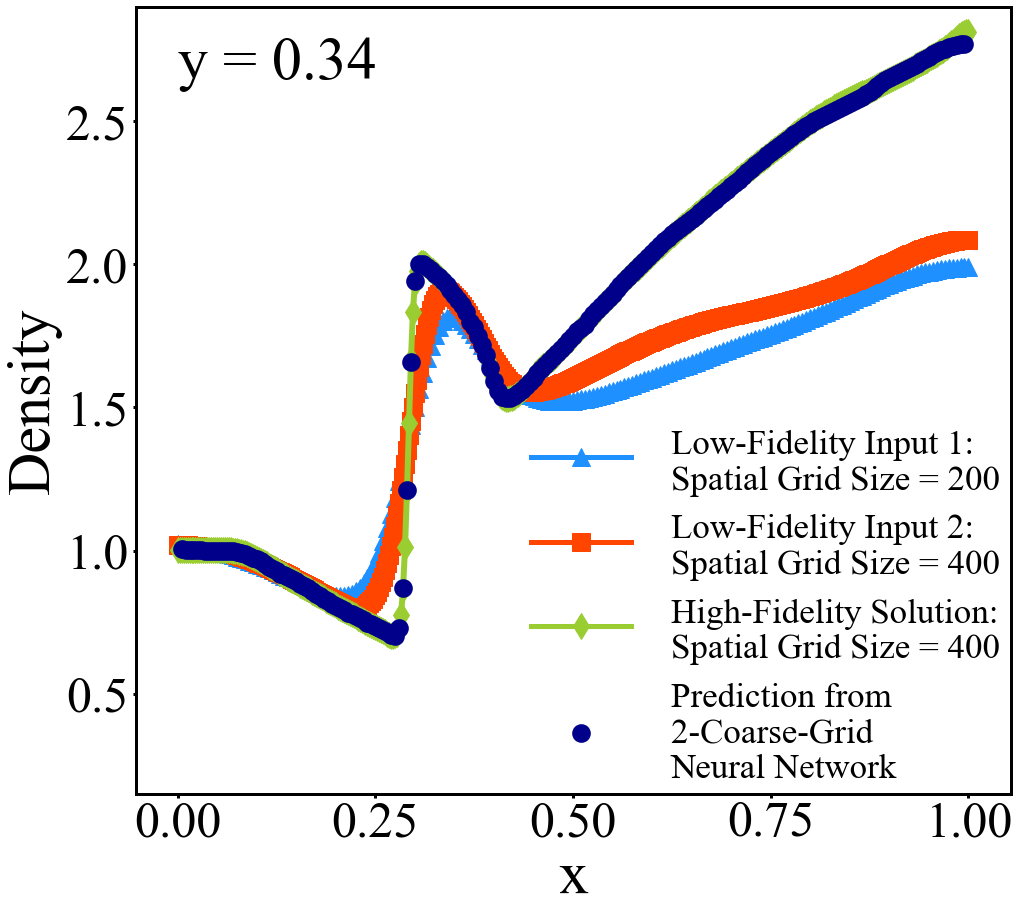}
    \caption{}
\end{subfigure}
\begin{subfigure}[b]{0.334\textwidth}
    \centering
    \includegraphics[width=1.0\linewidth]{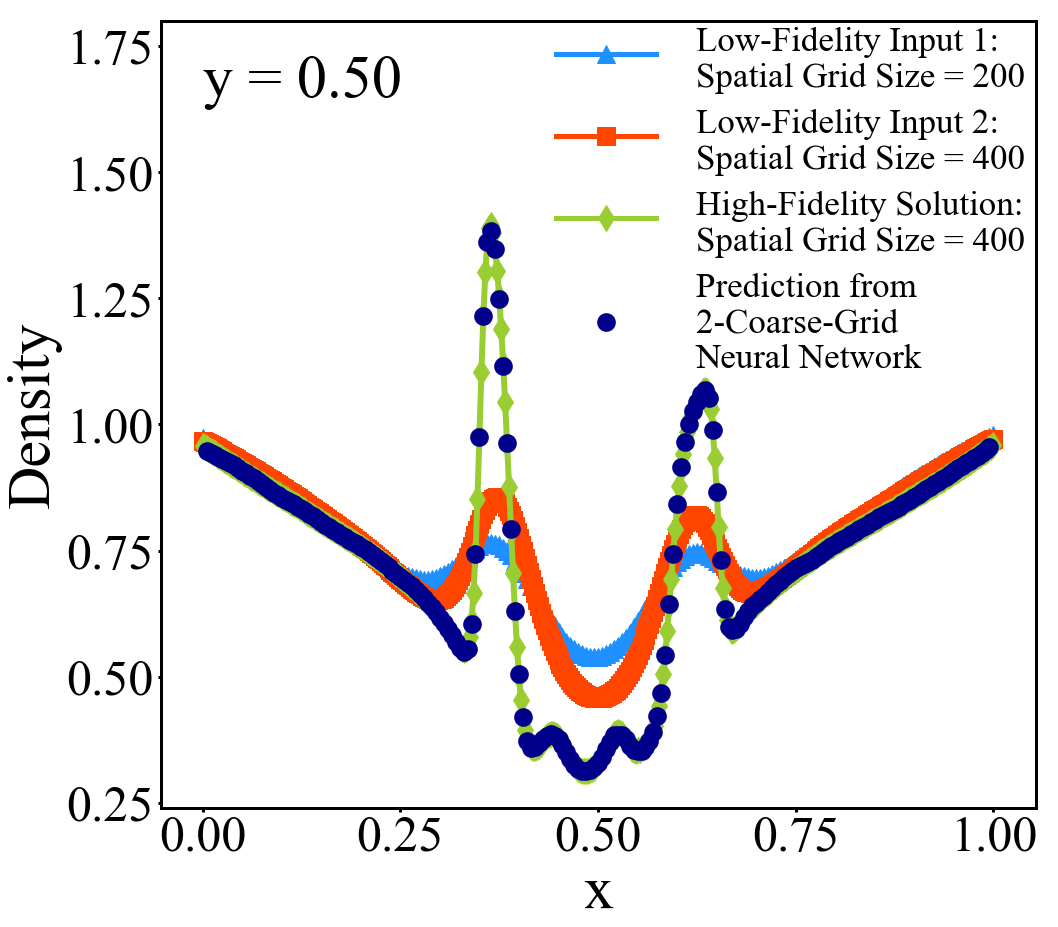}
    \caption{}
\end{subfigure}
\begin{subfigure}[b]{0.325\textwidth}
    \centering
    \includegraphics[width=1.0\linewidth]{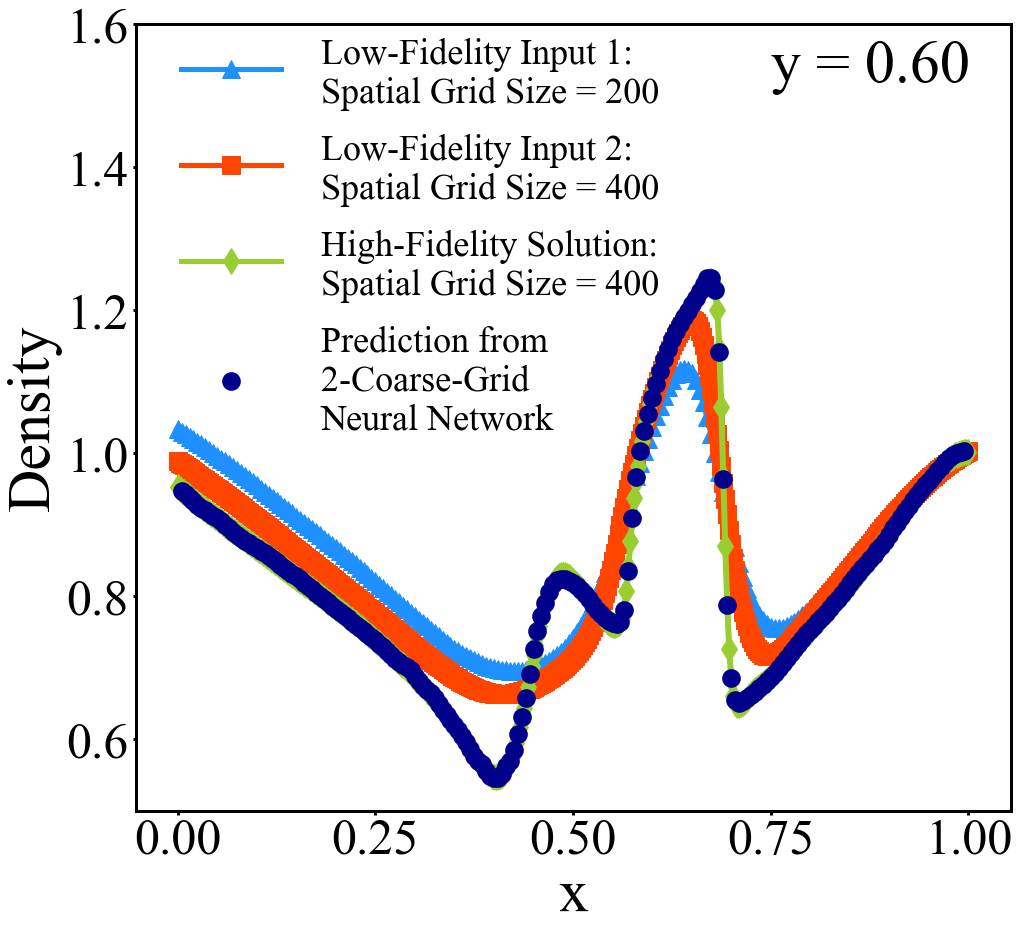}
    \caption{}
\end{subfigure}
\begin{subfigure}[b]{0.33\textwidth}
    \centering
    \includegraphics[width=1.0\linewidth]{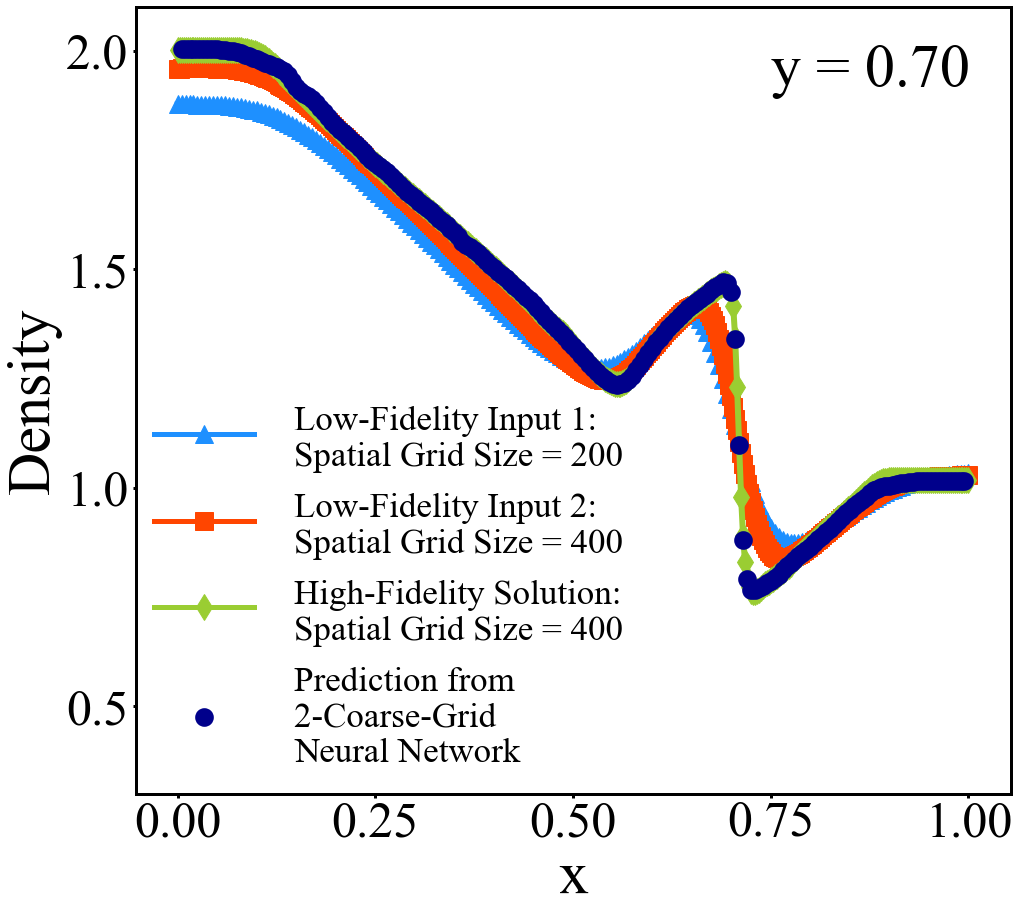}
    \caption{}
\end{subfigure}
\begin{subfigure}[b]{0.33\textwidth}
    \centering
    \includegraphics[width=1.0\linewidth]{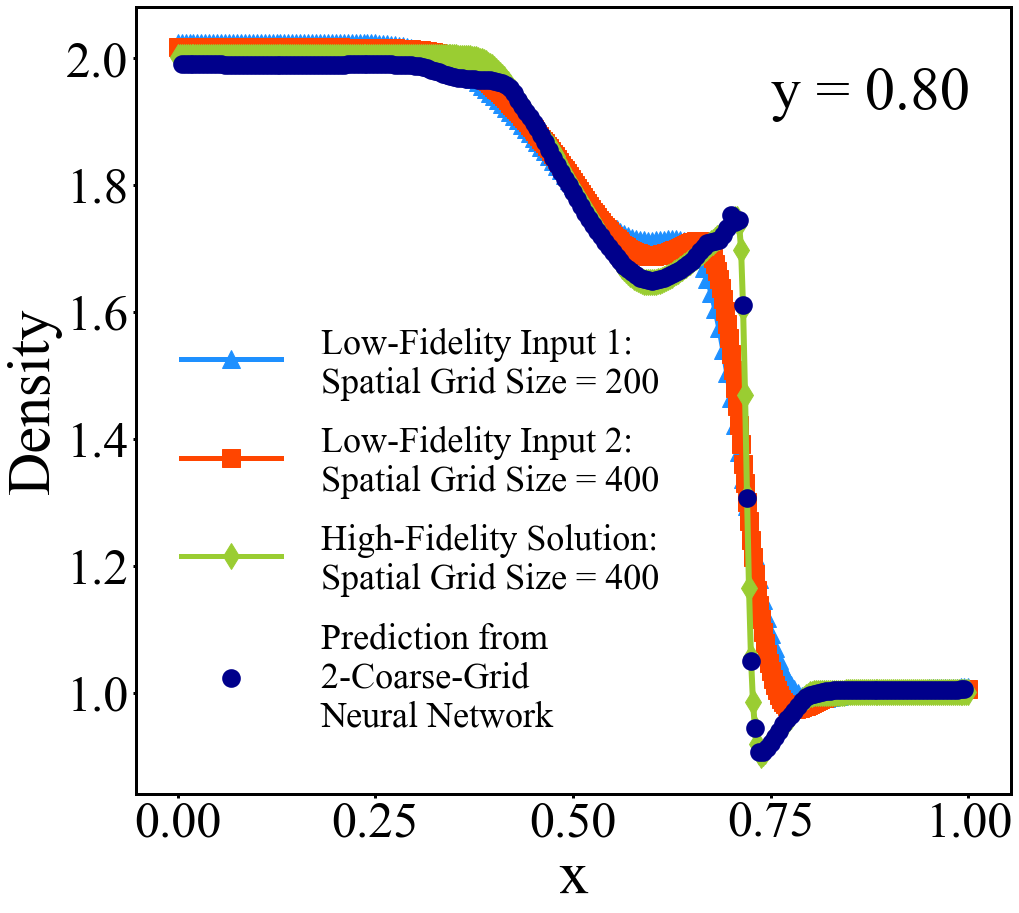}
    \caption{}
\end{subfigure}

\caption{Density cross-section profiles (along (a) \textit{y =} 0.34, (b) \textit{y =} 0.50, (c) \textit{y =} 0.60, (d) \textit{y =} 0.70, (e) \textit{y =} 0.80) of the 2CGNN prediction of the final-time $(t = 0.3)$ solution of {\bf Configuration 6} \cite{PDLax98}. Predicted density (dark blue), low-fidelity input solutions (blue and red) by leapfrog and diffusion splitting scheme (\ref{leapfrog-diffusion-splitting, 2D}) on  $200\times 200$ and $400\times 400$ grids respectively, and ``exact'' (reference) solution (green).} 
\label{2CGNN: Final time cross-section of config. 6, original}
\end{figure}

\begin{figure}[H] \centering
\begin{subfigure}[b]{.48\textwidth}
    \centering
    \includegraphics[width=1.0\linewidth]{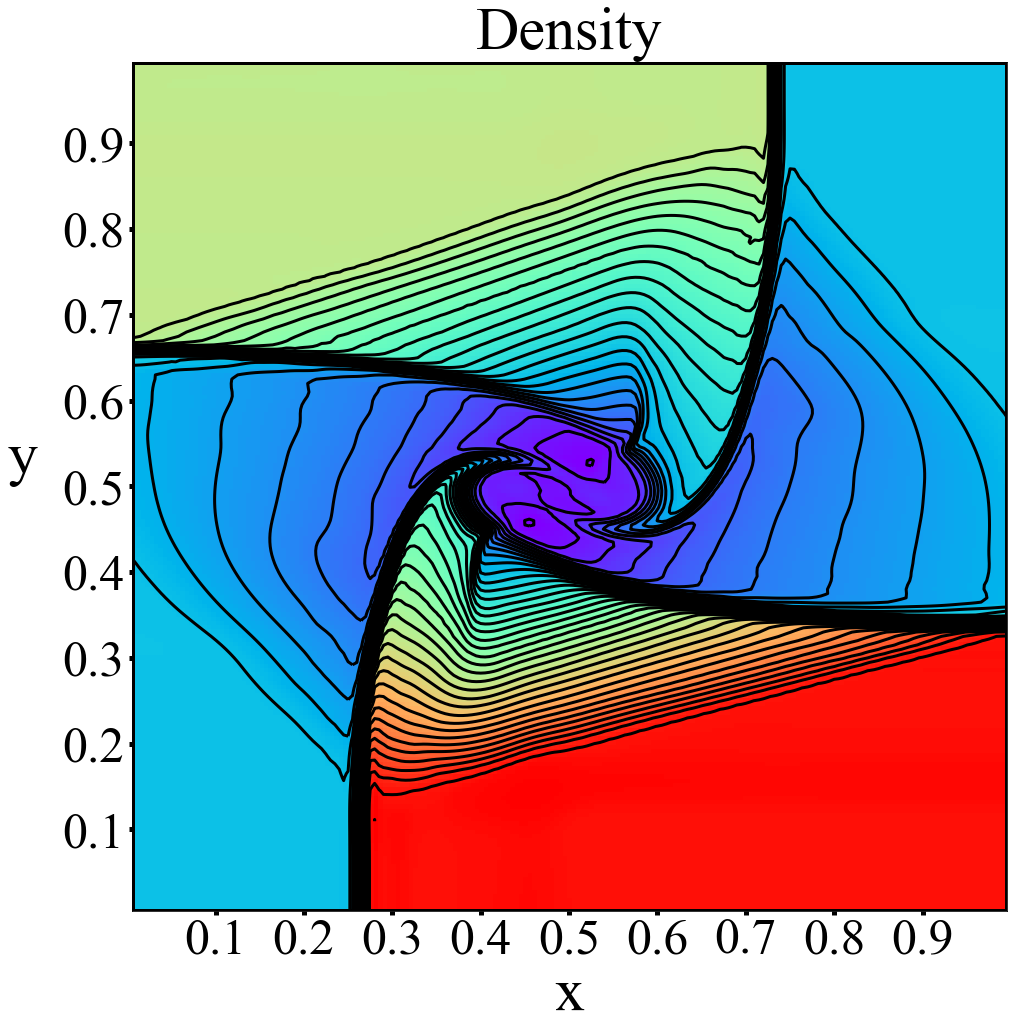}
    \caption{}
\end{subfigure}
\begin{subfigure}[b]{.48\textwidth}
    \centering
    \includegraphics[width=1.0\linewidth]{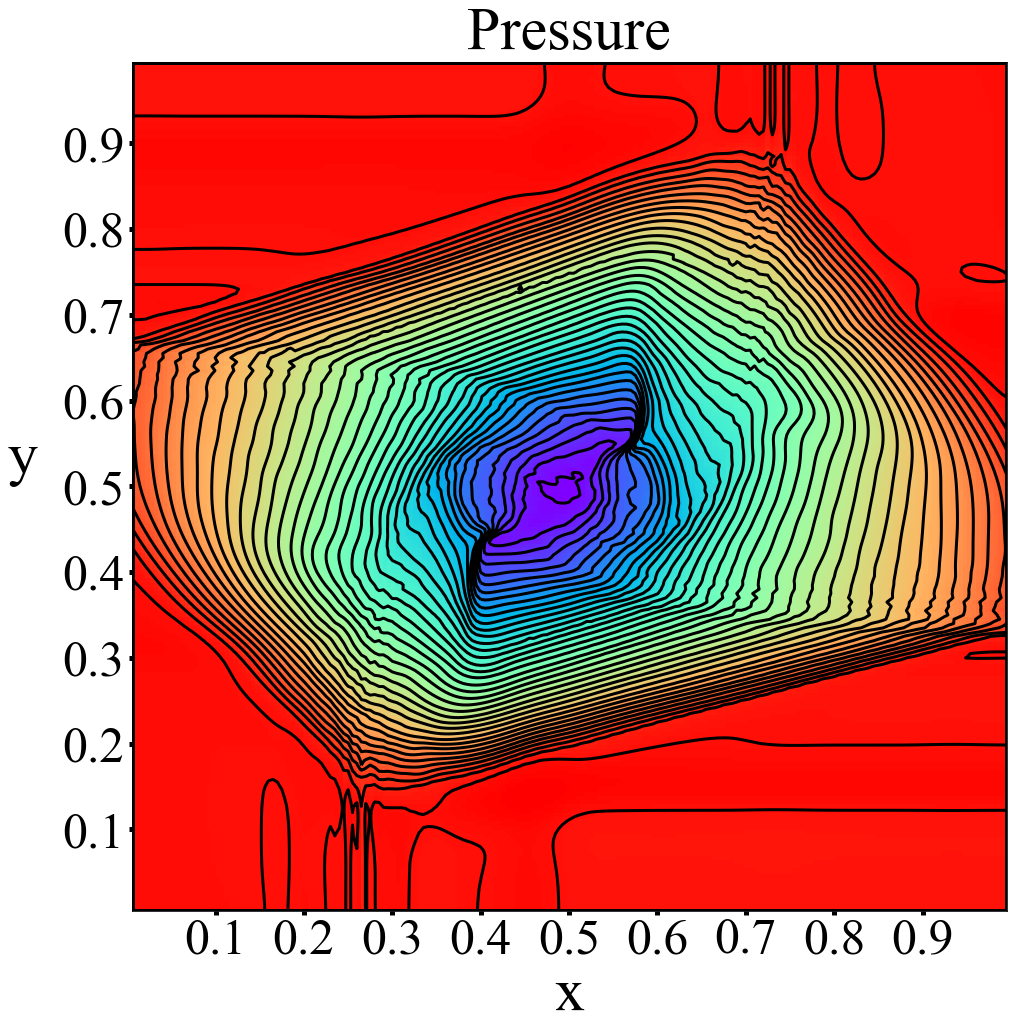}
    \caption{}
\end{subfigure}
\end{figure}

\begin{figure}[H]\ContinuedFloat
\centering
\begin{subfigure}[b]{.48\textwidth}
    \centering
    \includegraphics[width=1.0\linewidth]{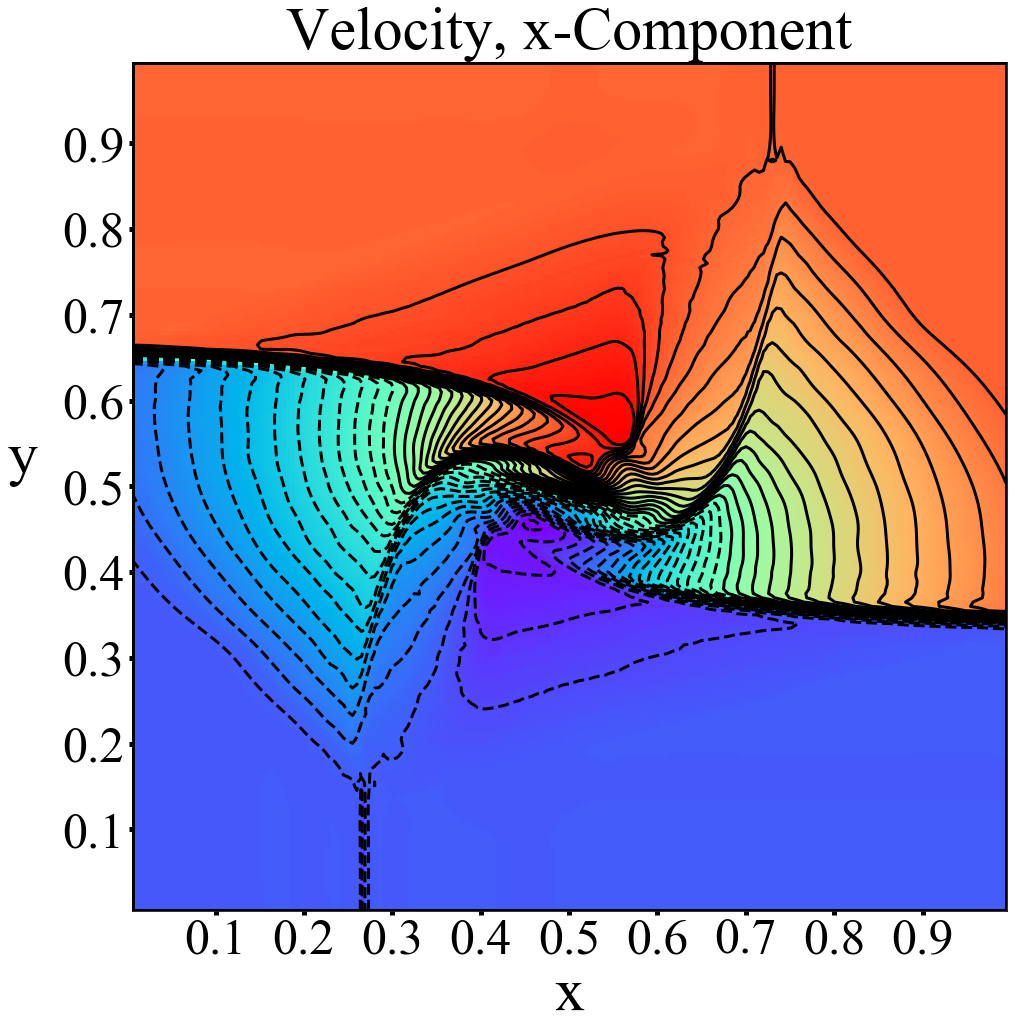}
    \caption{}
\end{subfigure}
\begin{subfigure}[b]{.48\textwidth}
    \centering
    \includegraphics[width=1.0\linewidth]{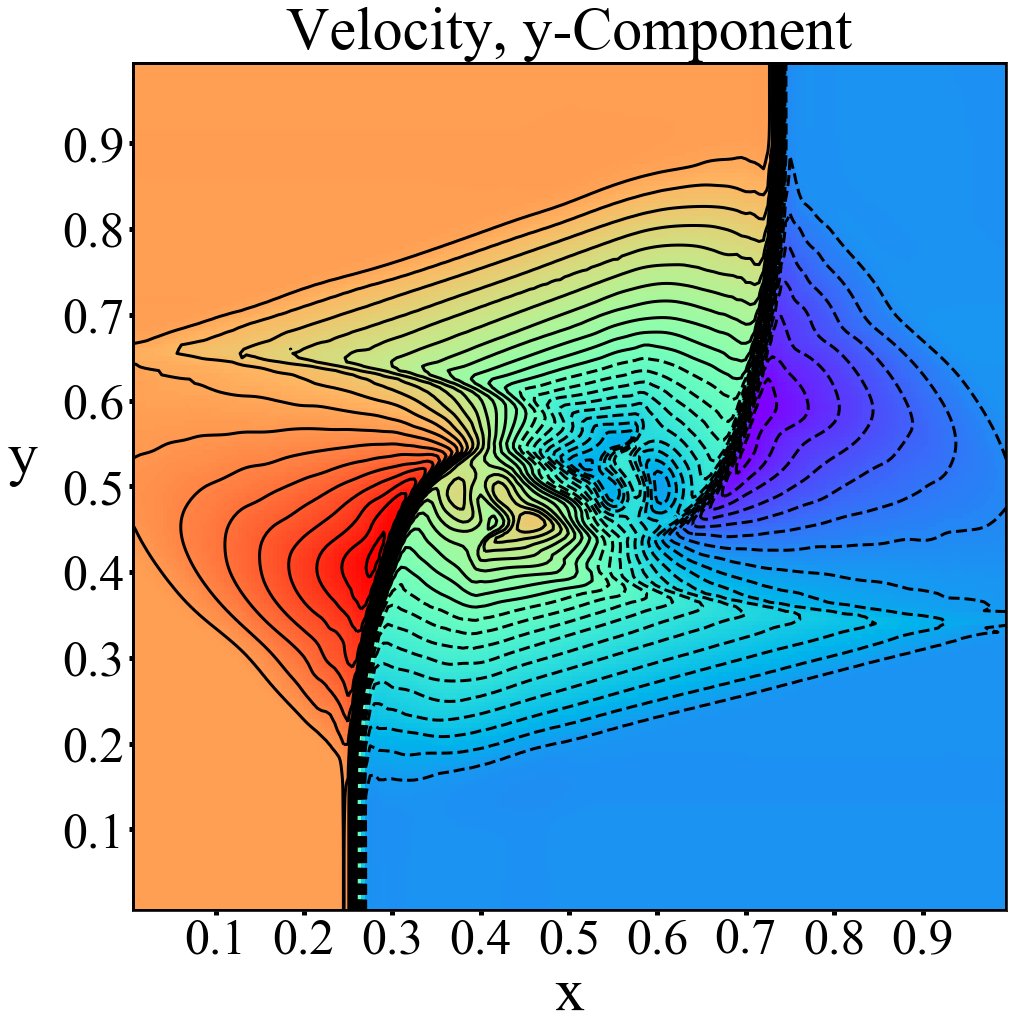}
    \caption{}
\end{subfigure}
\caption{2CGNN prediction of the final-time $(t=0.3)$ solution of the 2D Euler system, with {\bf initial value $+5\%$ perturbation of that of Configuration 6 in \cite{PDLax98}}: (a) density, (b) pressure, (c) velocity, x-component, (d) velocity, y-component.}
\label{2CGNN: Final time of config. 6, +5}
\end{figure}

\begin{figure}[H]\centering
\begin{subfigure}[b]{0.329\textwidth}
    \centering
    \includegraphics[width=1.0\linewidth]{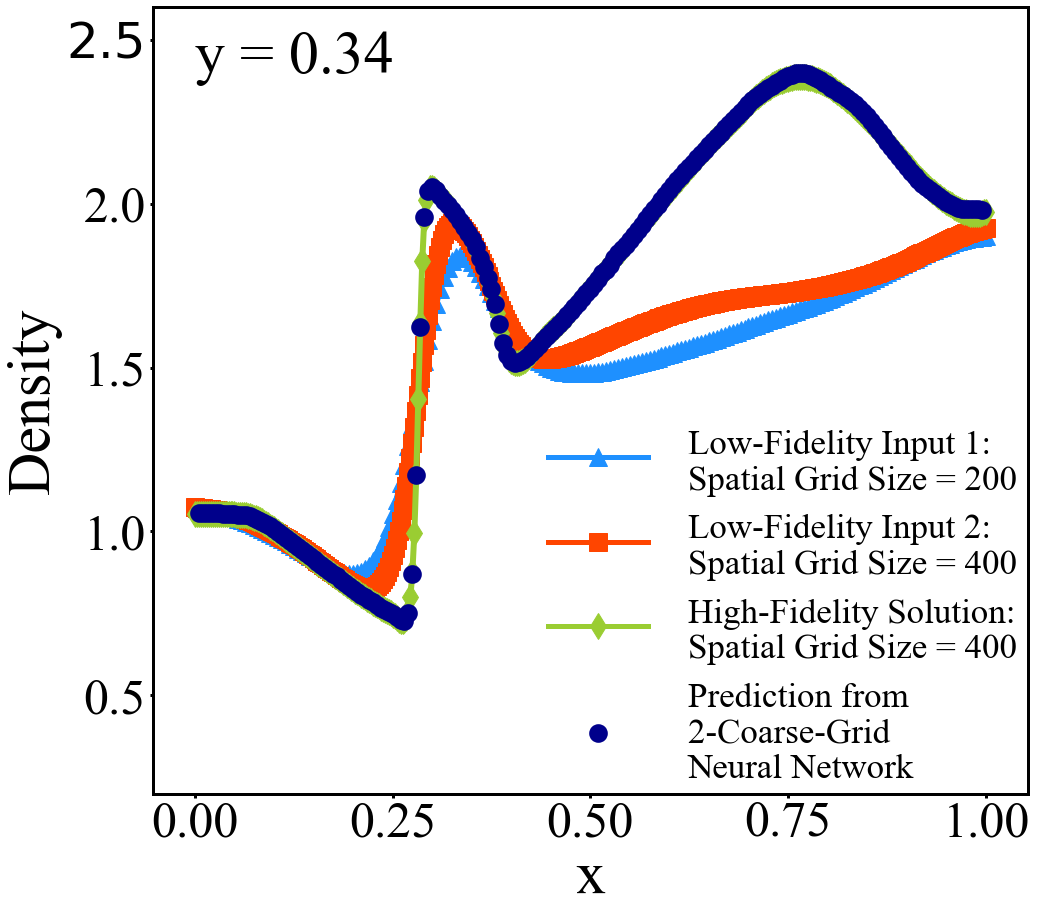}
    \caption{}
\end{subfigure}
\begin{subfigure}[b]{0.334\textwidth}
    \centering
    \includegraphics[width=1.0\linewidth]{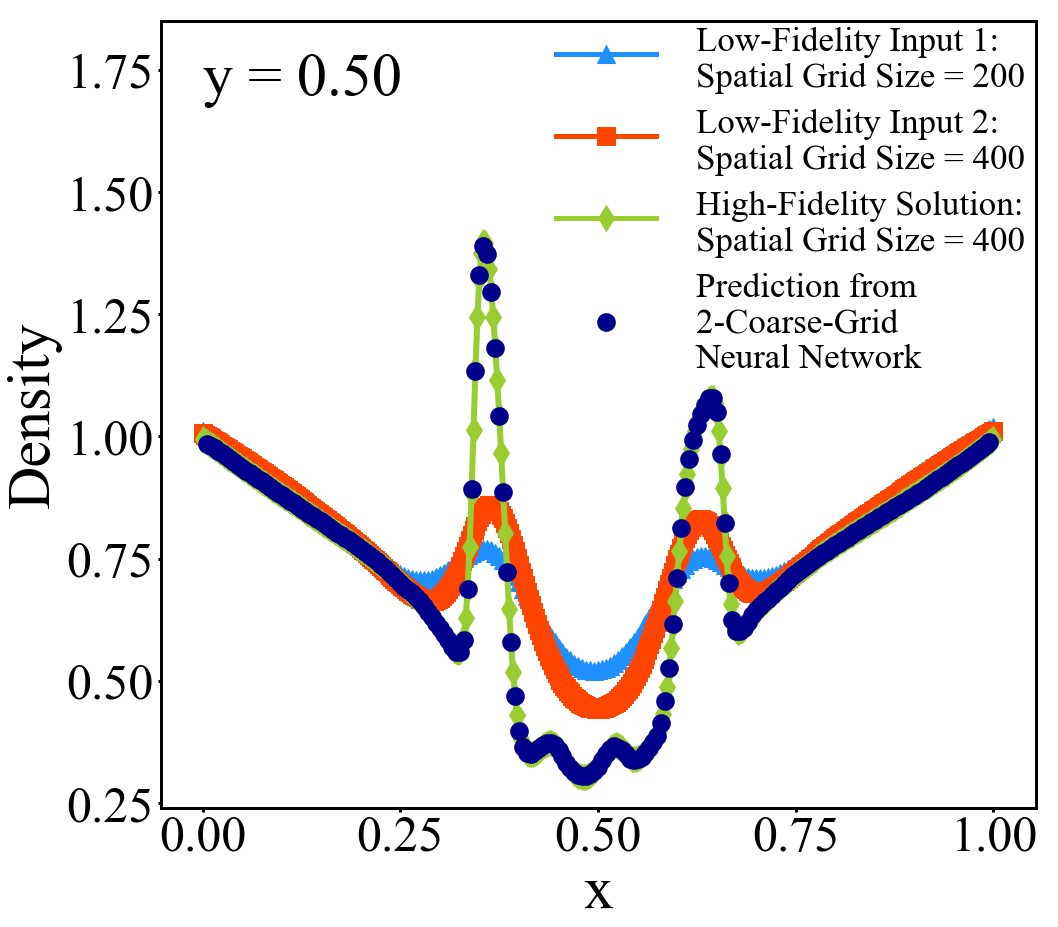}
    \caption{}
\end{subfigure}
\begin{subfigure}[b]{0.323\textwidth}
    \centering
    \includegraphics[width=1.0\linewidth]{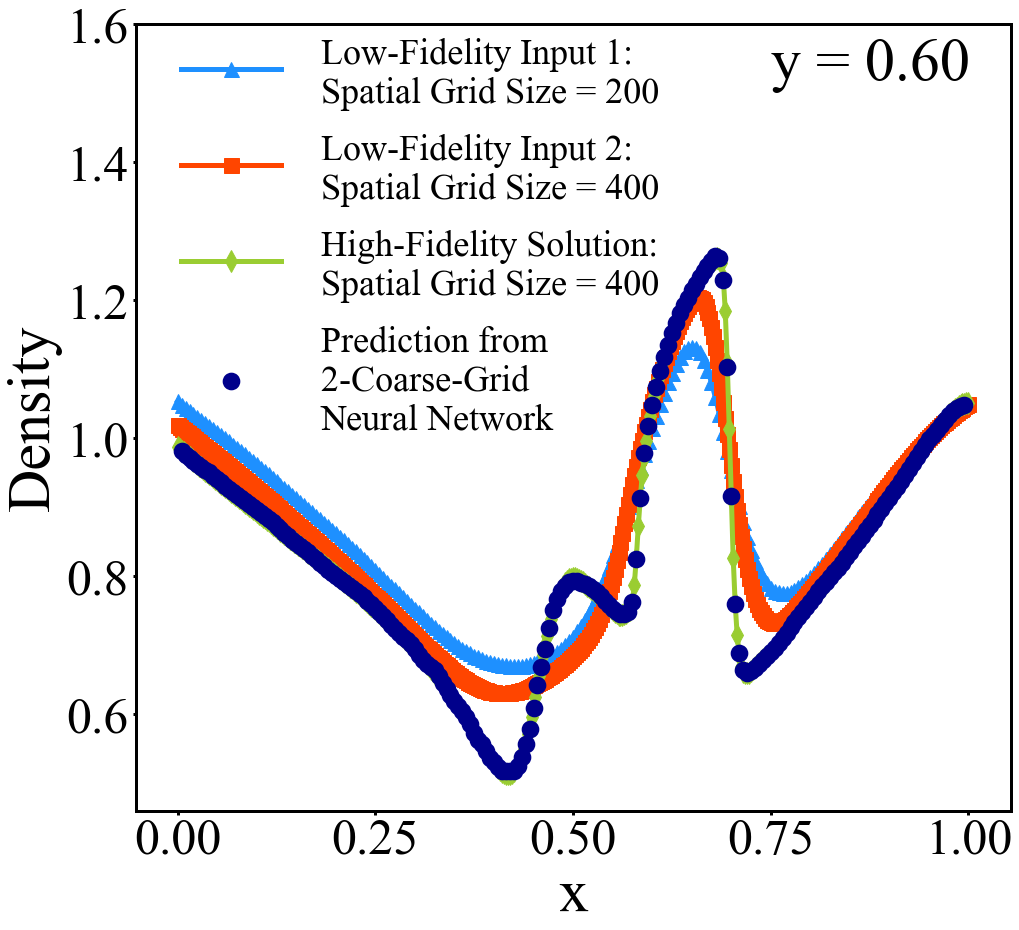}
    \caption{}
\end{subfigure}
\begin{subfigure}[b]{0.33\textwidth}
    \centering
    \includegraphics[width=1.0\linewidth]{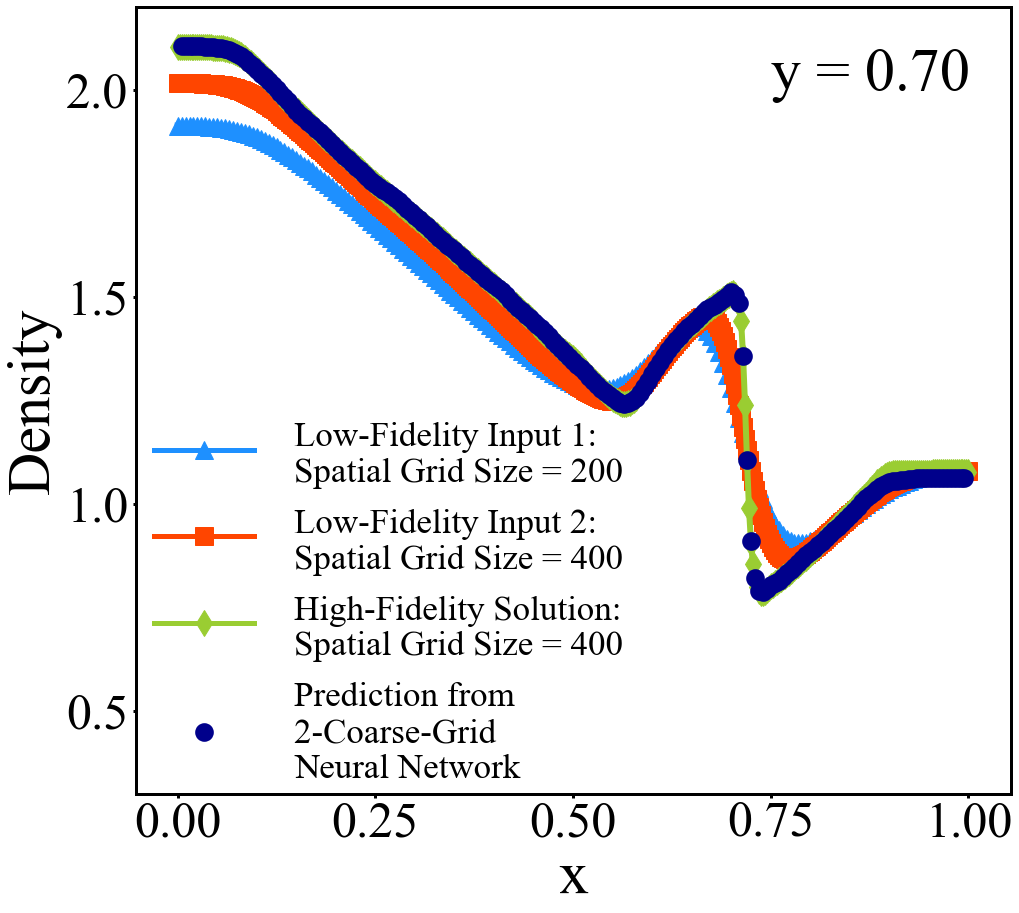}
    \caption{}
\end{subfigure}
\begin{subfigure}[b]{0.33\textwidth}
    \centering
    \includegraphics[width=1.0\linewidth]{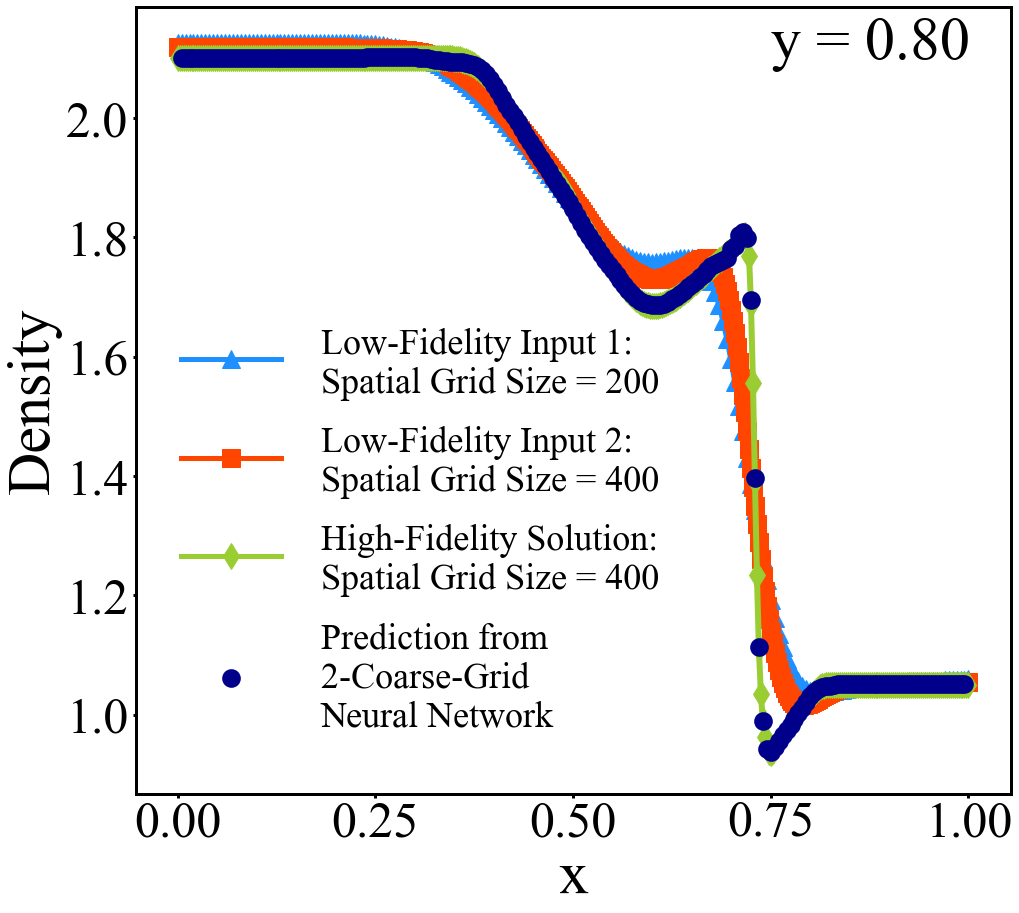}
    \caption{}
\end{subfigure}
\caption{Density cross-section profiles (along (a) \textit{y =} 0.34, (b) \textit{y =} 0.50, (c) \textit{y =} 0.60, (d) \textit{y =} 0.70, (e) \textit{y =} 0.80) of the 2CGNN prediction of the final-time $(t = 0.3)$ solution of the 2D Euler system, with {\bf initial value $+5\%$ perturbation of that of Configuration 6 in \cite{PDLax98}}. Predicted density (dark blue), low-fidelity input solutions (blue and red) by leapfrog and diffusion splitting scheme (\ref{leapfrog-diffusion-splitting, 2D}) on  $200\times 200$ and $400\times 400$ grids respectively, and ``exact'' (reference) solution (green).}
\label{2CGNN: Final time cross-section of config. 6, +5}
\end{figure}

We show below contour plots of predicted density for other the 2D Riemann problems from~\cite{PDLax98}, {\it i.e.}, Configurations 1, 2, 3, 4 and 8, and a cross-section profile perpendicular to y-axis for each problem.

Fig.~\ref{2CGNN: Final time of config. 1, original} shows the prediction of the final-time density solution of configuration 1 in \cite{PDLax98}, and the cross-section profile along $y = 0.50$.  
Fig.~\ref{2CGNN: Final time of config. 1, +5} shows the prediction of the final-time density solution of the 2D Euler system with initial value $+5\%$ perturbation of that of Configuration 1. The spatial computational domain is $[0.25, 0.95]\times [0.25, 0.95]$.

\begin{figure}[H] \centering
\begin{subfigure}[b]{.48\textwidth}
    \centering
    \includegraphics[width=1.0\linewidth]{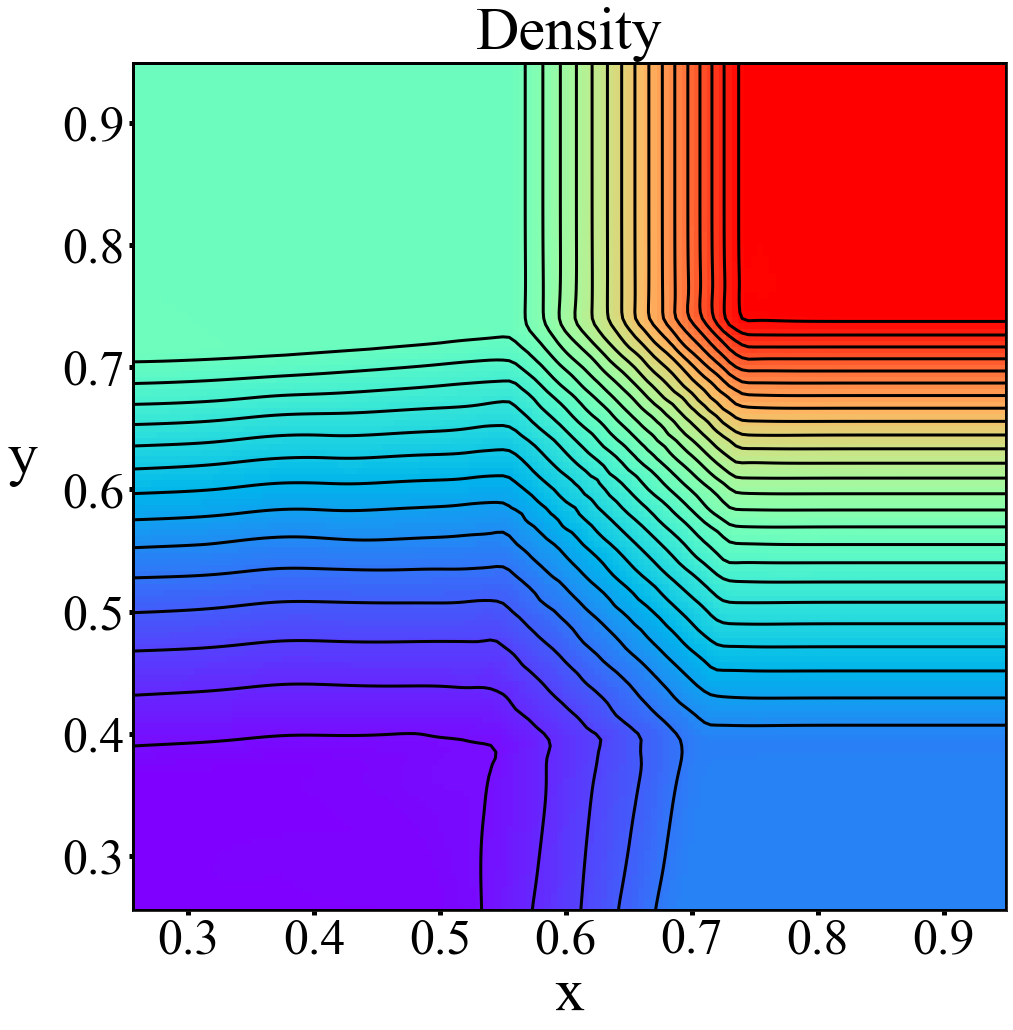}
    \caption{}
\end{subfigure}
\begin{subfigure}[b]{.48\textwidth}
    \centering
    \includegraphics[width=1.0\linewidth]{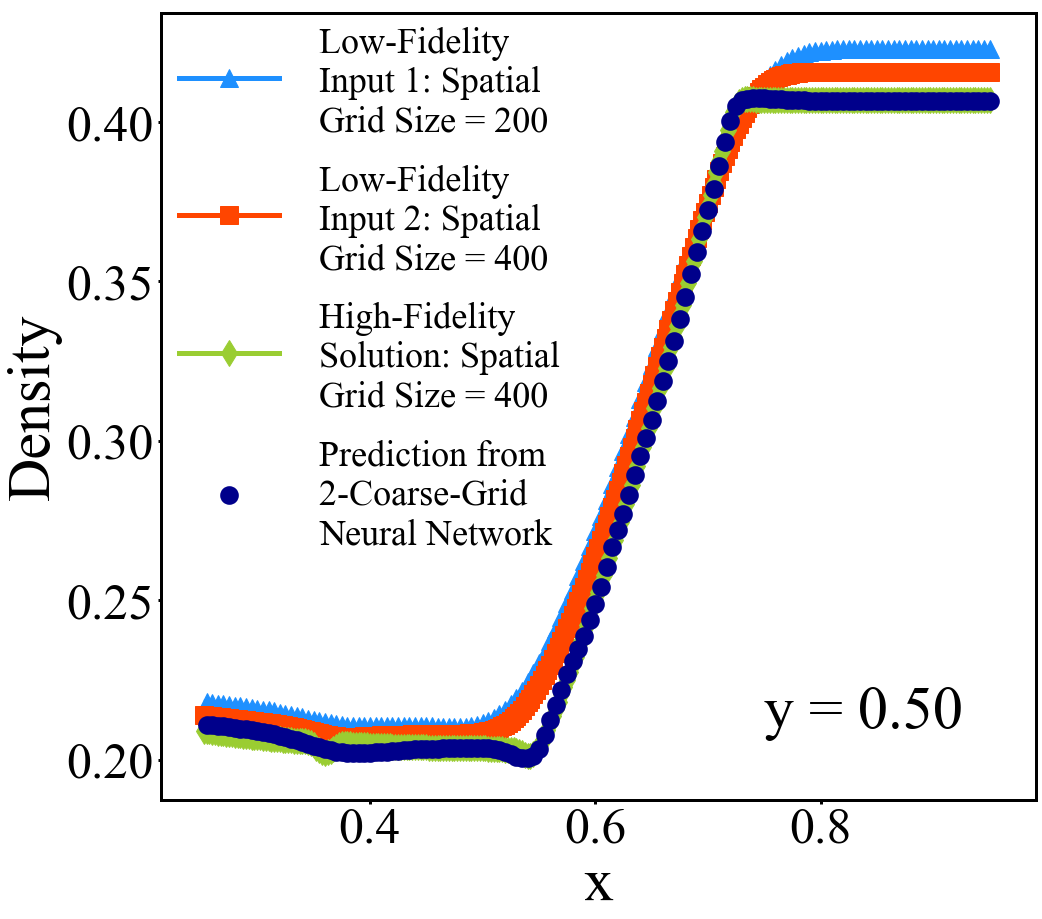}
    \caption{}
\end{subfigure}
\caption{2CGNN prediction of (a) the final-time $(t=0.2)$ density solution of \textbf{Configuration 1} in ~\cite{PDLax98}, and (b) its cross-section profile (dark blue) along \textit{y=} 0.50, compared to low-fidelity input solutions (blue and red) by leapfrog and diffusion splitting scheme (\ref{leapfrog-diffusion-splitting, 2D}) on $200\times 200$ and $400\times 400$ grids, respectively, and ``exact'' (reference) solution (green).}
\label{2CGNN: Final time of config. 1, original}
\end{figure}

\begin{figure}[H]
\centering
\begin{subfigure}[b]{.48\textwidth}
    \centering
    \includegraphics[width=1.0\linewidth]{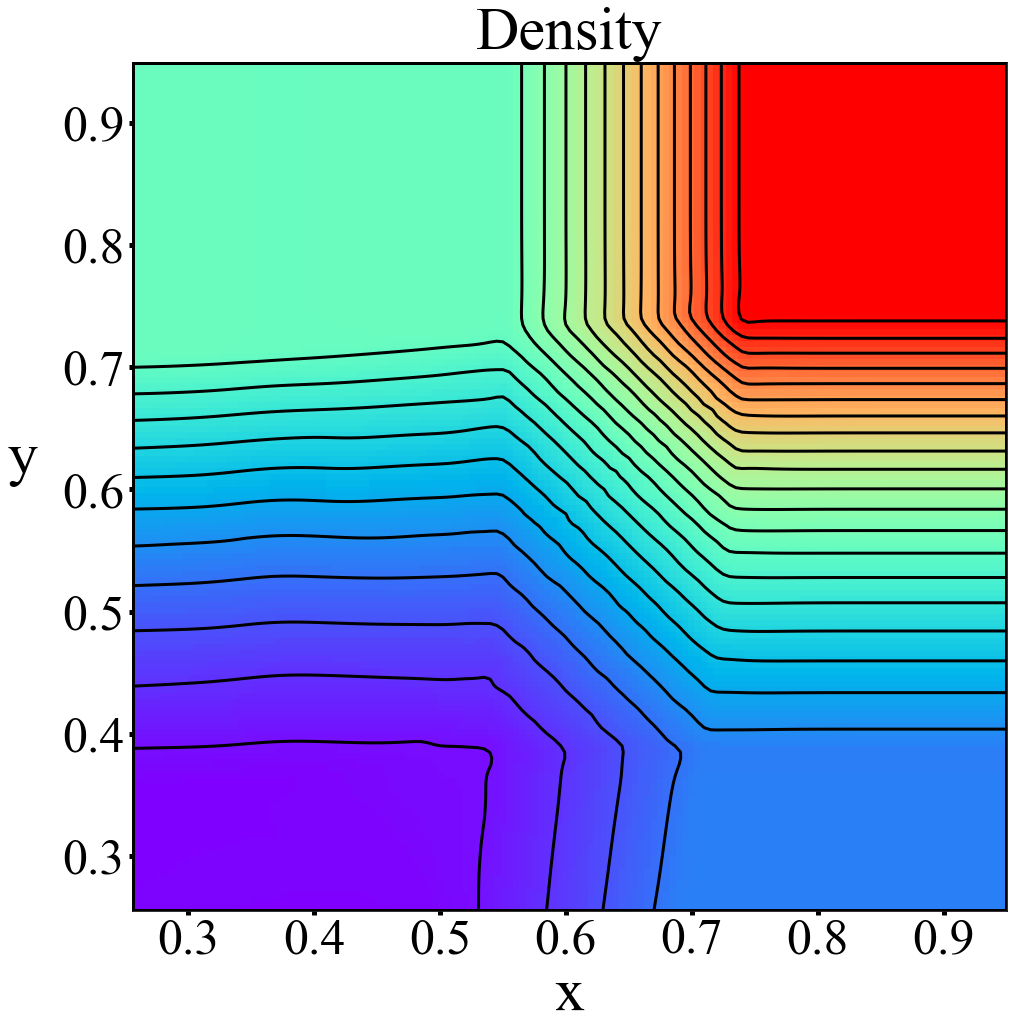}
    \caption{}
\end{subfigure}
\begin{subfigure}[b]{.48\textwidth}
    \centering
    \includegraphics[width=1.0\linewidth]{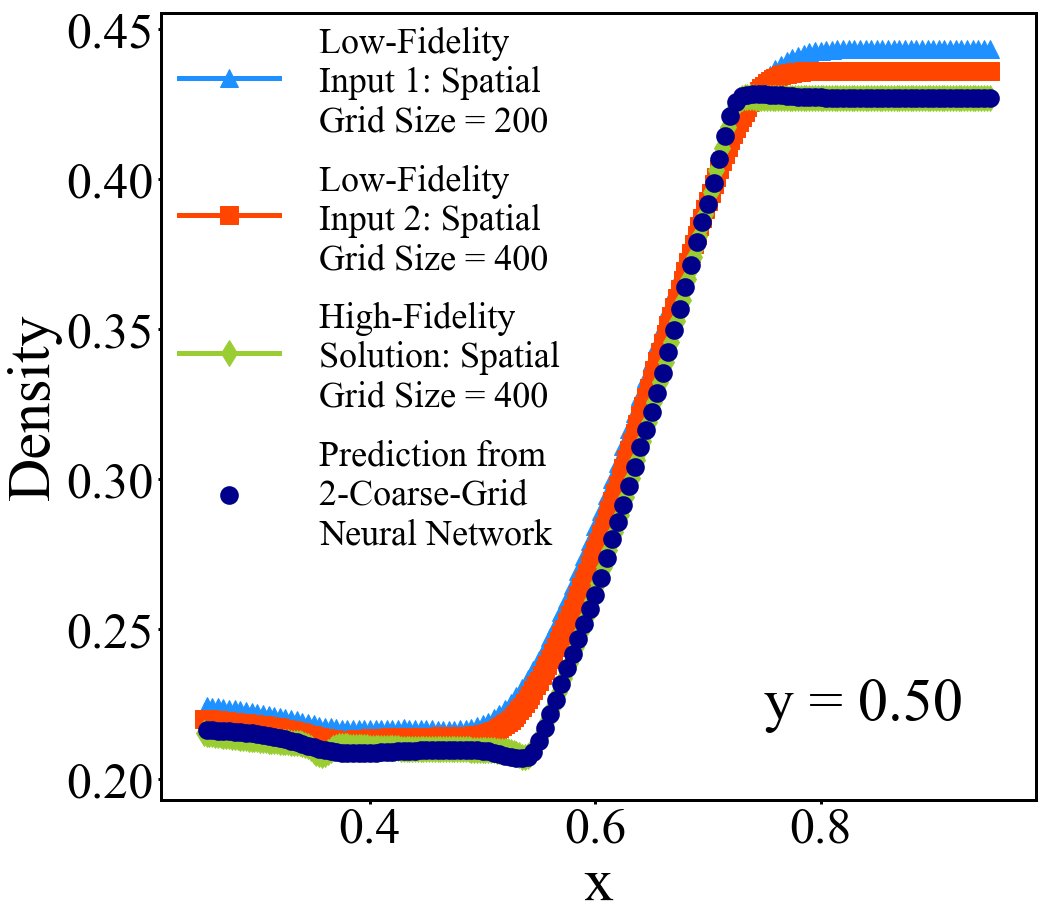}
    \caption{}
\end{subfigure}
\caption{2CGNN prediction of (a) final-time $(t=0.2)$ density solution of the 2D Euler system with {\bf initial value $+5\%$ perturbation of that of Configuration 1}, and (b) its cross-section profile (dark blue) at \textit{y=} 0.50, compared to low-fidelity input solutions (blue and red) by leapfrog and diffusion splitting scheme (\ref{leapfrog-diffusion-splitting, 2D}) on $200\times 200$ and $400\times 400$ grids, respectively, and ``exact'' (reference) solution (green).}
\label{2CGNN: Final time of config. 1, +5}
\end{figure}

Fig.~\ref{2CGNN: Final time of config. 2, original} shows the predicted final-time density solution of Configuration 2.  
Fig.~\ref{2CGNN: Final time of config. 2, +5} shows the predicted final-time density solution of the 2D Euler system with initial value $+5\%$ perturbation of that of Configuration 2. The spatial computational domain is $x,y\in [0, 0.85]$ for this configuration.

\begin{figure}[H] \centering
\begin{subfigure}[b]{.48\textwidth}
    \centering
    \includegraphics[width=1.0\linewidth]{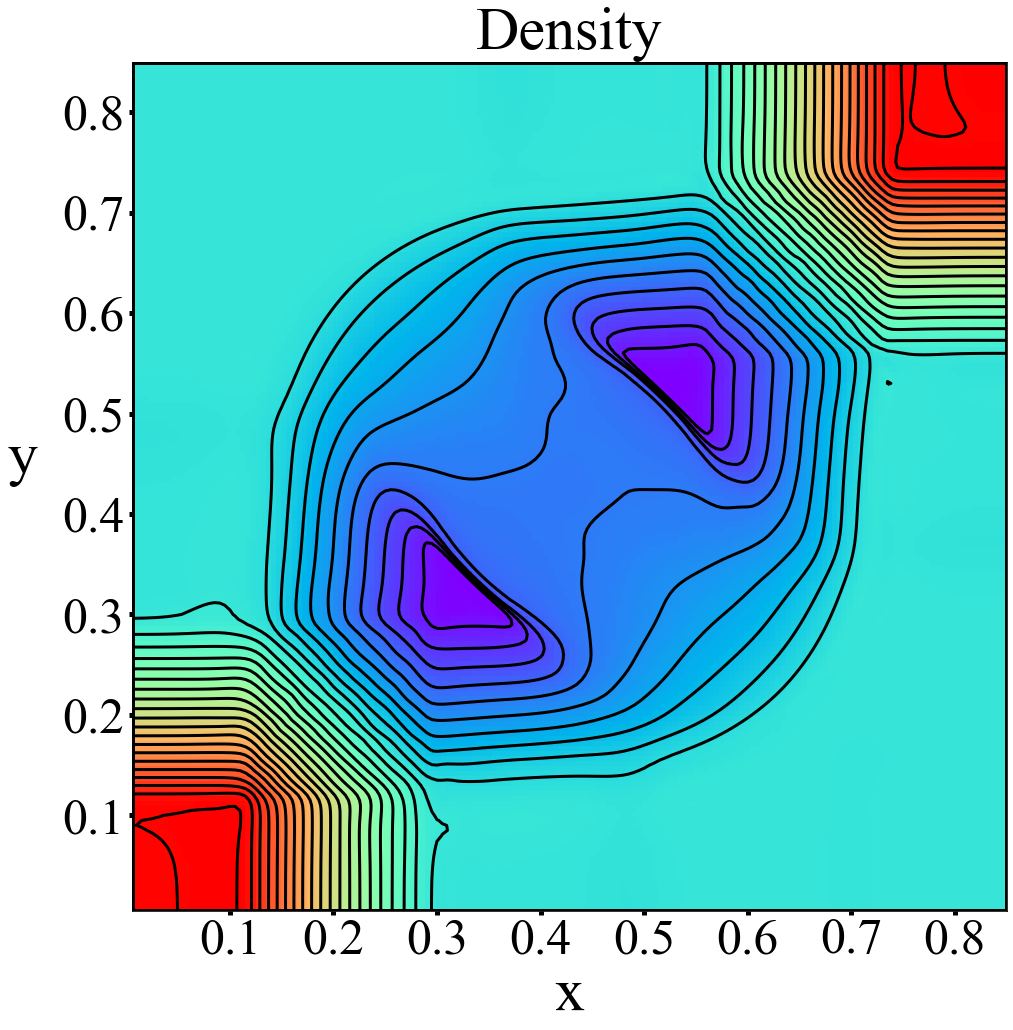}
    \caption{}
\end{subfigure}
\begin{subfigure}[b]{.48\textwidth}
    \centering
    \includegraphics[width=1.0\linewidth]{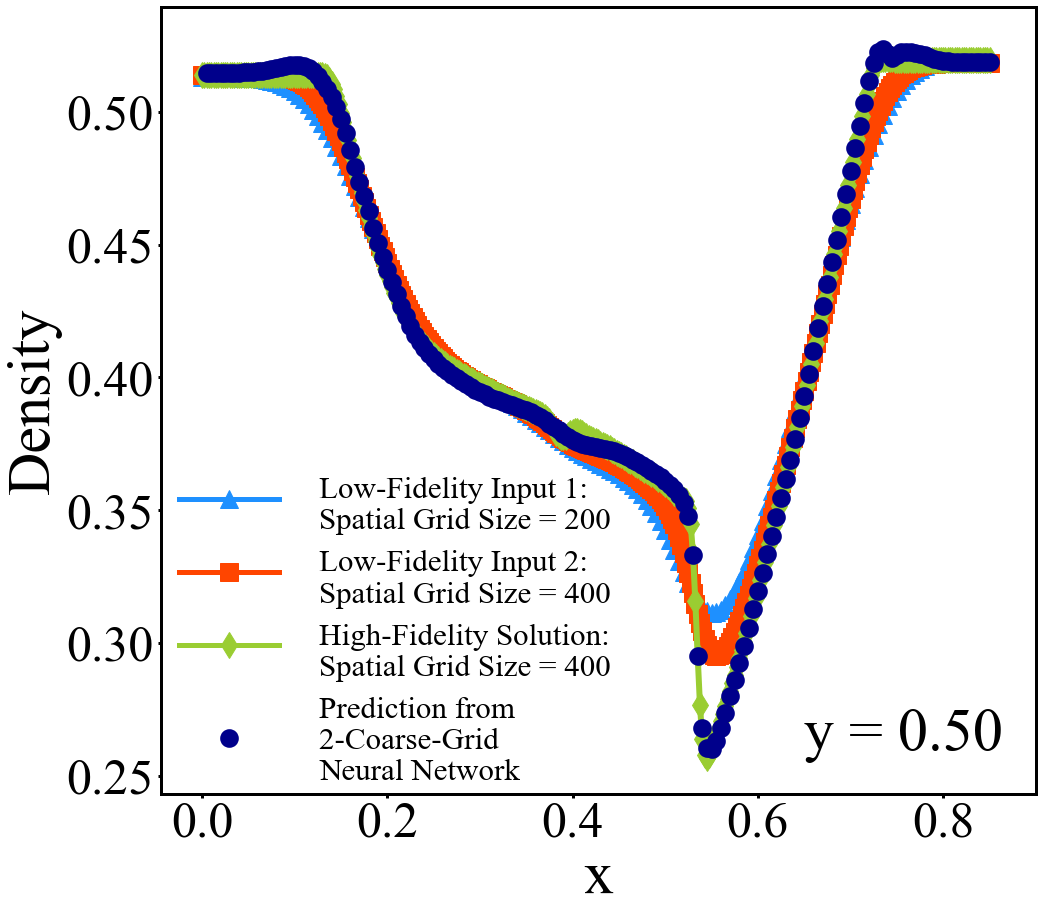}
    \caption{}
\end{subfigure}
\caption{2CGNN prediction of (a) the final-time $(t=0.2)$ density solution of \textbf{Configuration 2}, and (b) its cross-section profile (dark blue) along \textit{y=}0.50, compared to low-fidelity input solutions (blue and red) by leapfrog and diffusion splitting scheme (\ref{leapfrog-diffusion-splitting, 2D}) on $200\times 200$ and $400\times 400$ grids, respectively, and ``exact'' (reference) solution (green).}
\label{2CGNN: Final time of config. 2, original}
\end{figure}

\begin{figure}[H]
\centering
\begin{subfigure}[b]{.48\textwidth}
    \centering
    \includegraphics[width=1.0\linewidth]{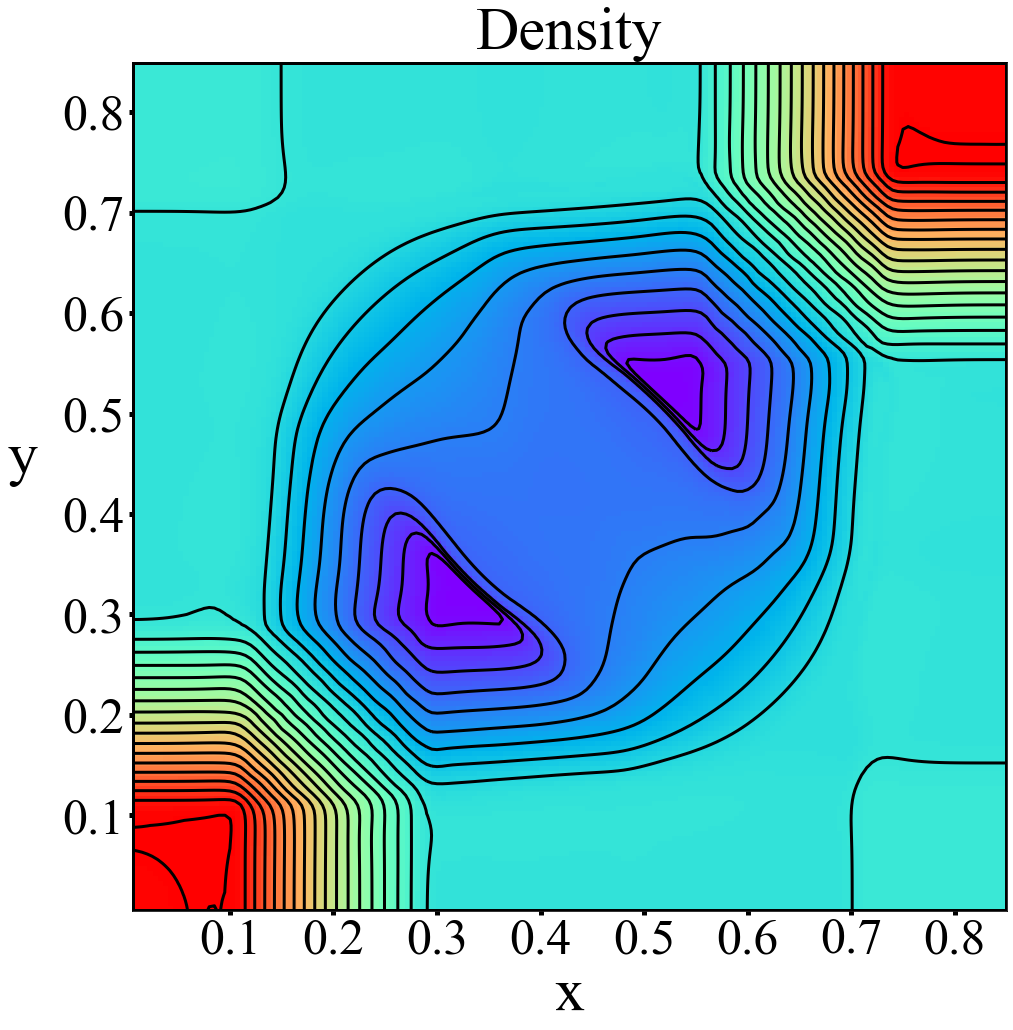}
    \caption{}
\end{subfigure}
\begin{subfigure}[b]{.48\textwidth}
    \centering
    \includegraphics[width=1.0\linewidth]{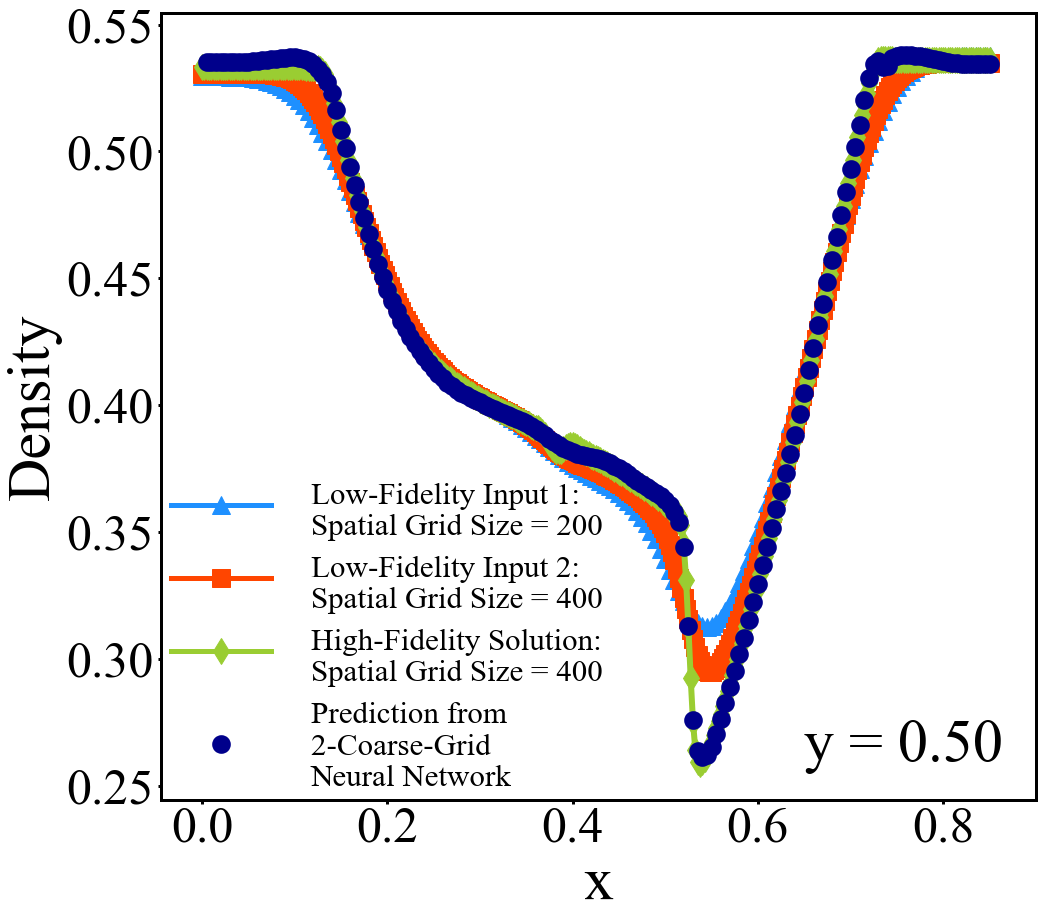}
    \caption{}
\end{subfigure}
\caption{2CGNN prediction of (a) the final-time $(t=0.2)$ density solution of the 2-D Euler system with {\bf initial value $+5\%$ perturbation of that of Configuration 2}, and (b) its cross-section profile (dark blue) at \textit{y=}0.50, compared to low-fidelity input solutions (blue and red) by leapfrog and diffusion splitting scheme (\ref{leapfrog-diffusion-splitting, 2D}) on $200\times 200$ and $400\times 400$ grids, respectively, and ``exact'' (reference) solution (green).}
\label{2CGNN: Final time of config. 2, +5}
\end{figure}

Fig.~\ref{2CGNN: Final time of config. 3, original} and ~\ref{2CGNN: Final time of config. 3, +5} show the predicted final-time solutions of the 2-D Euler system with initial values the original initial values and $+5\%$ perturbation of that of Configuration 3, respectively. The spatial computational domain is $x,y\in [0, 0.9]$.  

\begin{figure}[H] \centering
\begin{subfigure}[b]{.48\textwidth}
    \centering
    \includegraphics[width=1.0\linewidth]{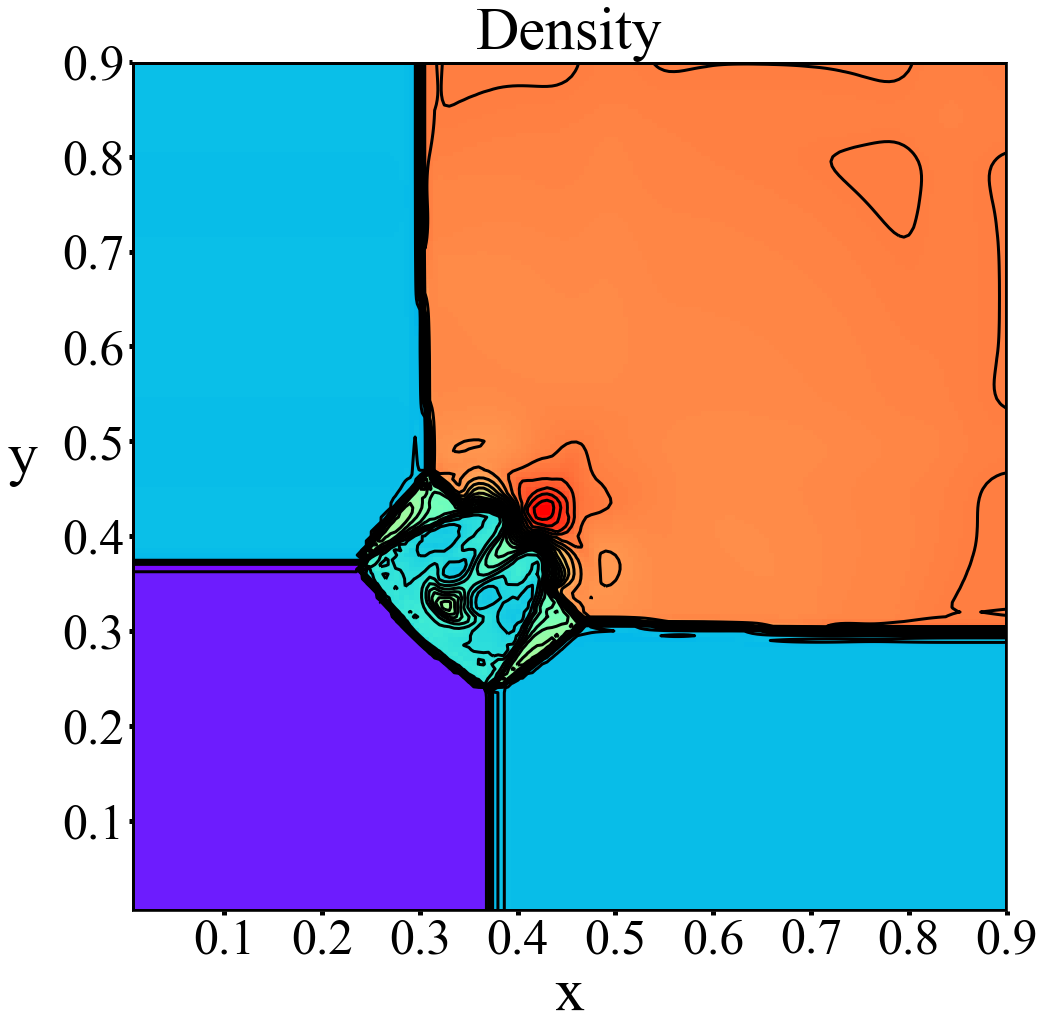}
    \caption{}
\end{subfigure}
\begin{subfigure}[b]{.48\textwidth}
    \centering
    \includegraphics[width=1.0\linewidth]{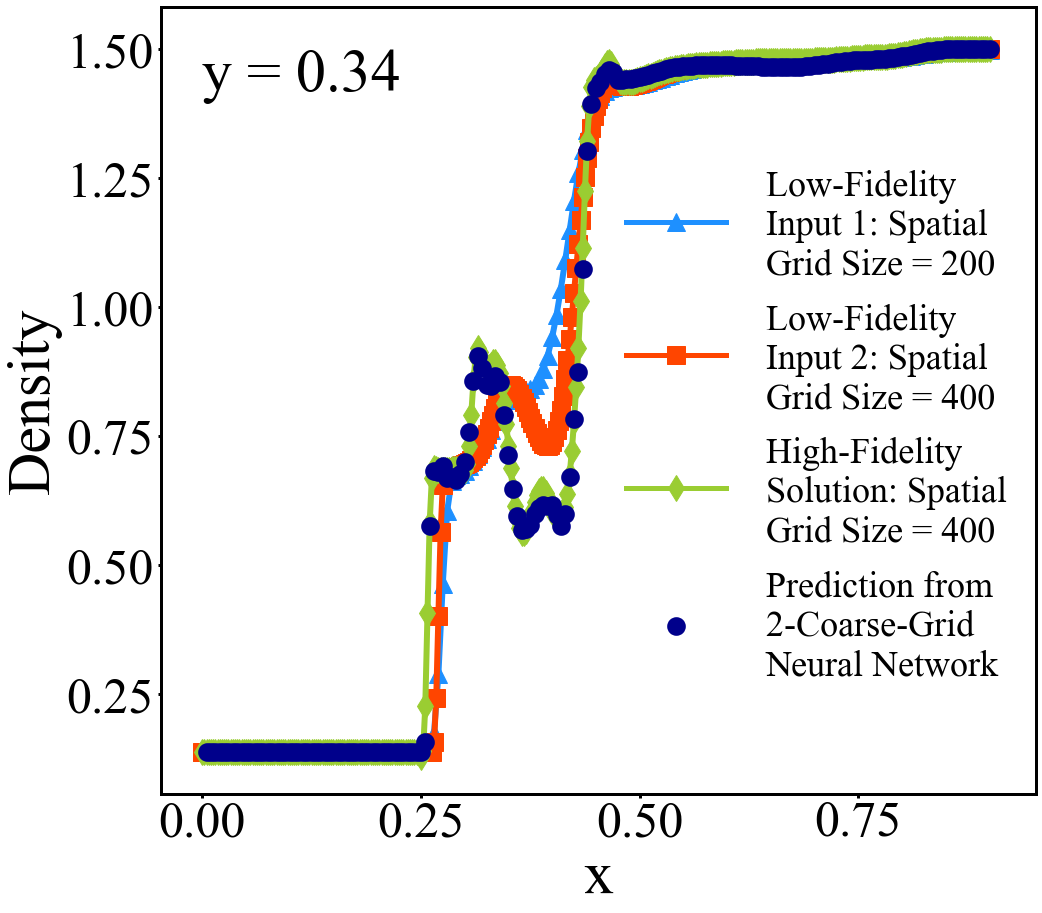}
    \caption{}
\end{subfigure}
\caption{2CGNN prediction of (a) the final-time $(t=0.3)$ density solution of \textbf{Configuration 3}, and (b) its cross-section profile (dark blue) along \textit{y=}0.34, compared to low-fidelity input solutions (blue and red) by leapfrog and diffusion splitting scheme (\ref{leapfrog-diffusion-splitting, 2D}) on $200\times 200$ and $400\times 400$ grids, respectively, and ``exact'' (reference) solution (green).}
\label{2CGNN: Final time of config. 3, original}
\end{figure}

\begin{figure}[H]
\centering
\begin{subfigure}[b]{.48\textwidth}
    \centering
    \includegraphics[width=1.0\linewidth]{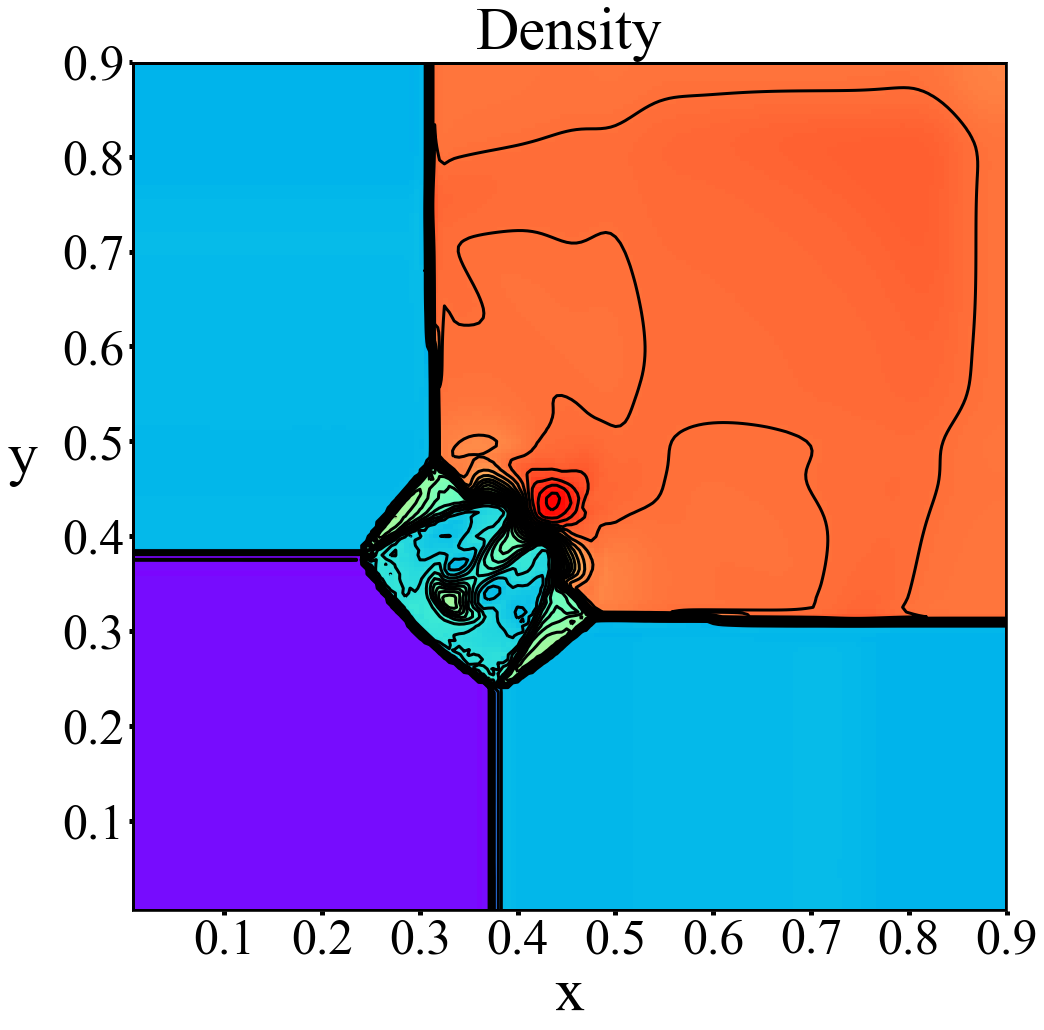}
    \caption{}
\end{subfigure}
\begin{subfigure}[b]{.48\textwidth}
    \centering
    \includegraphics[width=1.0\linewidth]{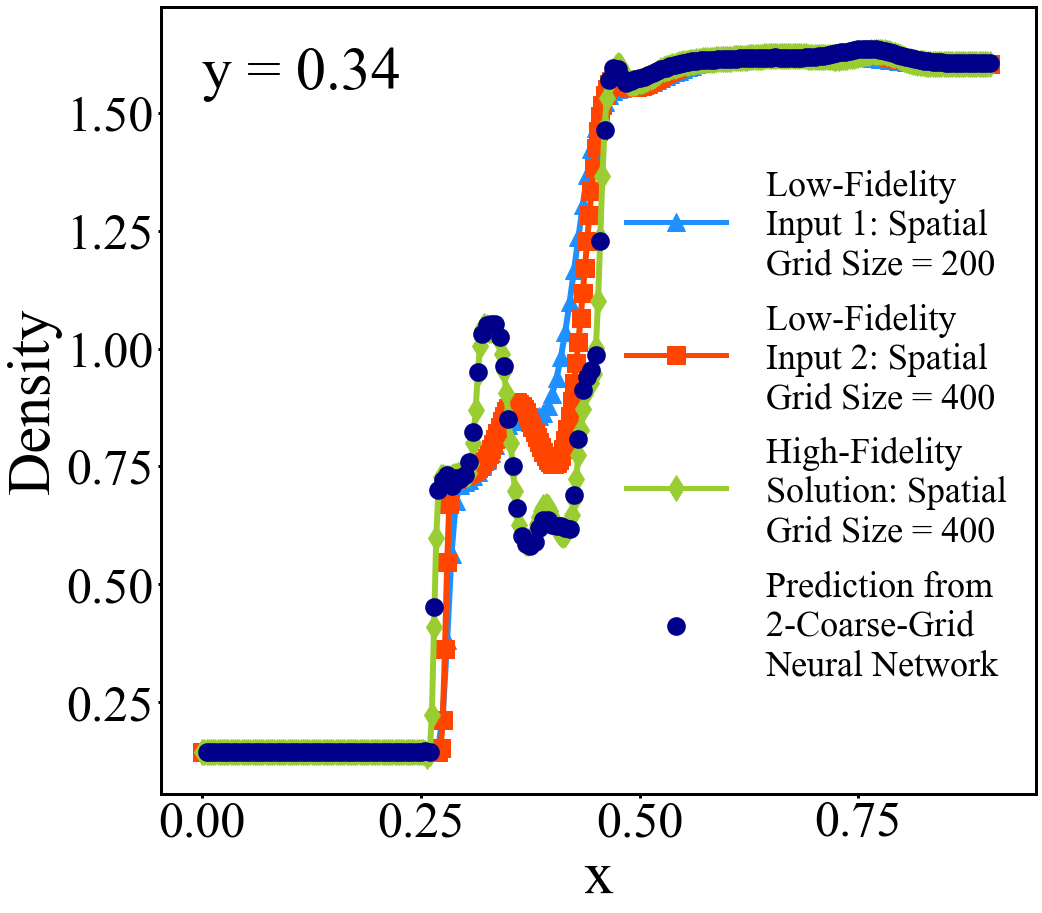}
    \caption{}
\end{subfigure}
\caption{2CGNN prediction of (a) the final-time $(t=0.3)$ density solution of the 2D Euler system with {\bf initial value $+5\%$ perturbation of that of Configuration 3}, and (b) its cross-section profile (dark blue) along \textit{y=}0.34, compared to low-fidelity input solutions (blue and red) by leapfrog and diffusion splitting scheme (\ref{leapfrog-diffusion-splitting, 2D}) on $200\times 200$ and $400\times 400$ grids, respectively, and ``exact'' (reference) solution (green).}
\label{2CGNN: Final time of config. 3, +5}
\end{figure}

The low-cost input solutions of Configuration 3 from the first order scheme (\ref{leapfrog-diffusion-splitting}) may not be qualitatively similar to the high-fidelity solution in certain areas for 2CGNN to make an accurate prediction. To improve the quality of inputs, we compute two input solutions on $100\times 100$ and $200\times 200$
grids by the same $4$th order scheme used for computing reference solutions. Fig. ~\ref{2CGNN: Final time of config. 3 src input, original} and ~\ref{2CGNN: Final time of config. 3 src input, +5} show the improved predictions.  Near the center of the region, the predictions capture fine details of the solution better than those in Fig.~\ref{2CGNN: Final time of config. 3, original} and \ref{2CGNN: Final time of config. 3, +5}.    

\begin{figure}[H] \centering
\begin{subfigure}[b]{.48\textwidth}
    \centering
    \includegraphics[width=1.0\linewidth]{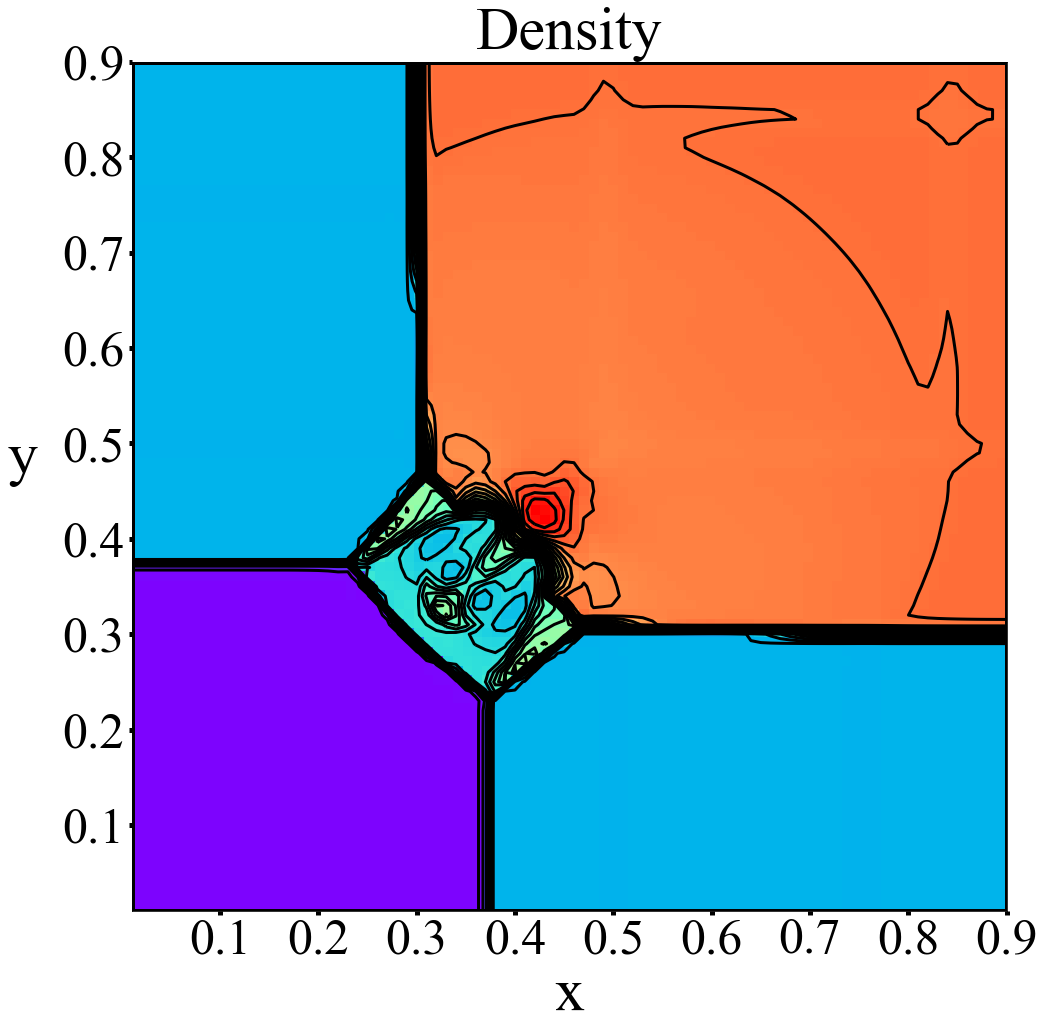}
    \caption{}
\end{subfigure}
\begin{subfigure}[b]{.48\textwidth}
    \centering
    \includegraphics[width=1.0\linewidth]{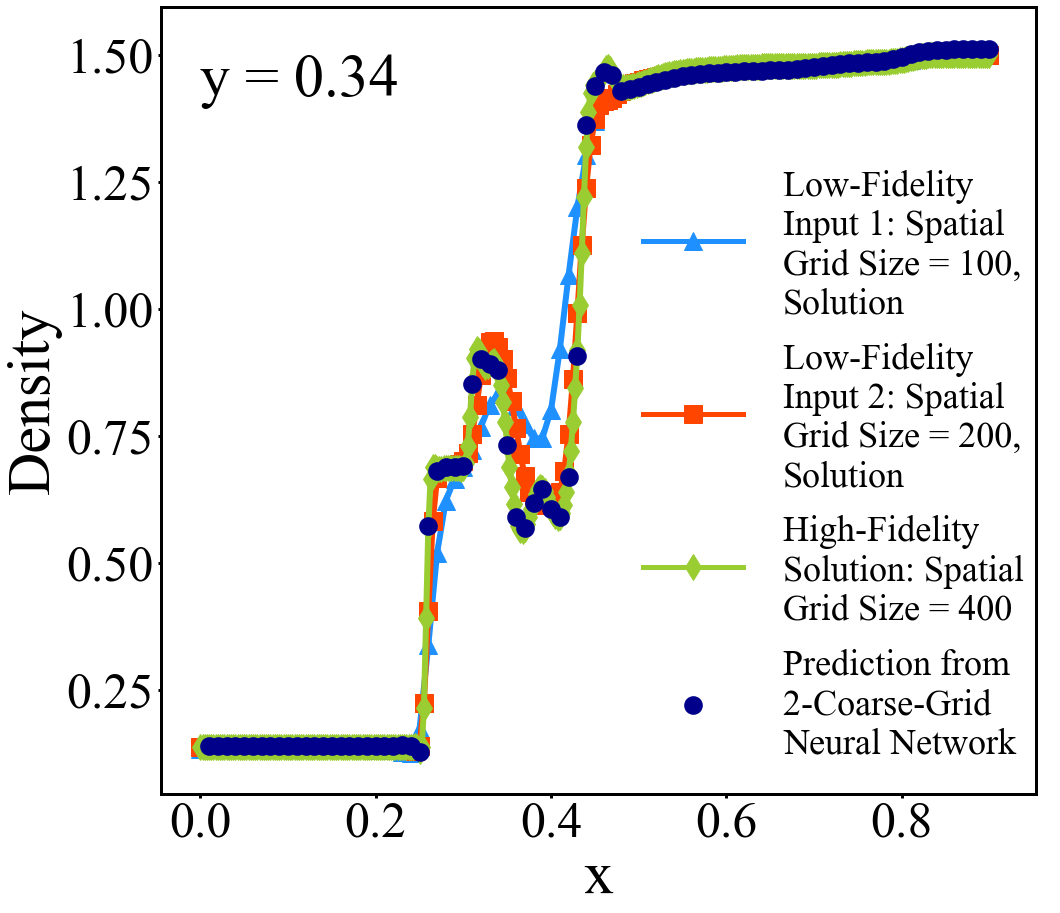}
    \caption{}
\end{subfigure}
\caption{2CGNN prediction of (a) the final-time $(t=0.3)$ density solution of \textbf{Configuration 3}, and (b) its cross-section profile (dark blue) along \textit{y=}0.34, compared to low-fidelity input solutions (blue and red) by a $4$th order scheme on $100\times 100$ and $200\times 200$ grids, respectively, and ``exact'' (reference) solution (green).}
\label{2CGNN: Final time of config. 3 src input, original}
\end{figure}

\begin{figure}[H]
\centering
\begin{subfigure}[b]{.48\textwidth}
    \centering
    \includegraphics[width=1.0\linewidth]{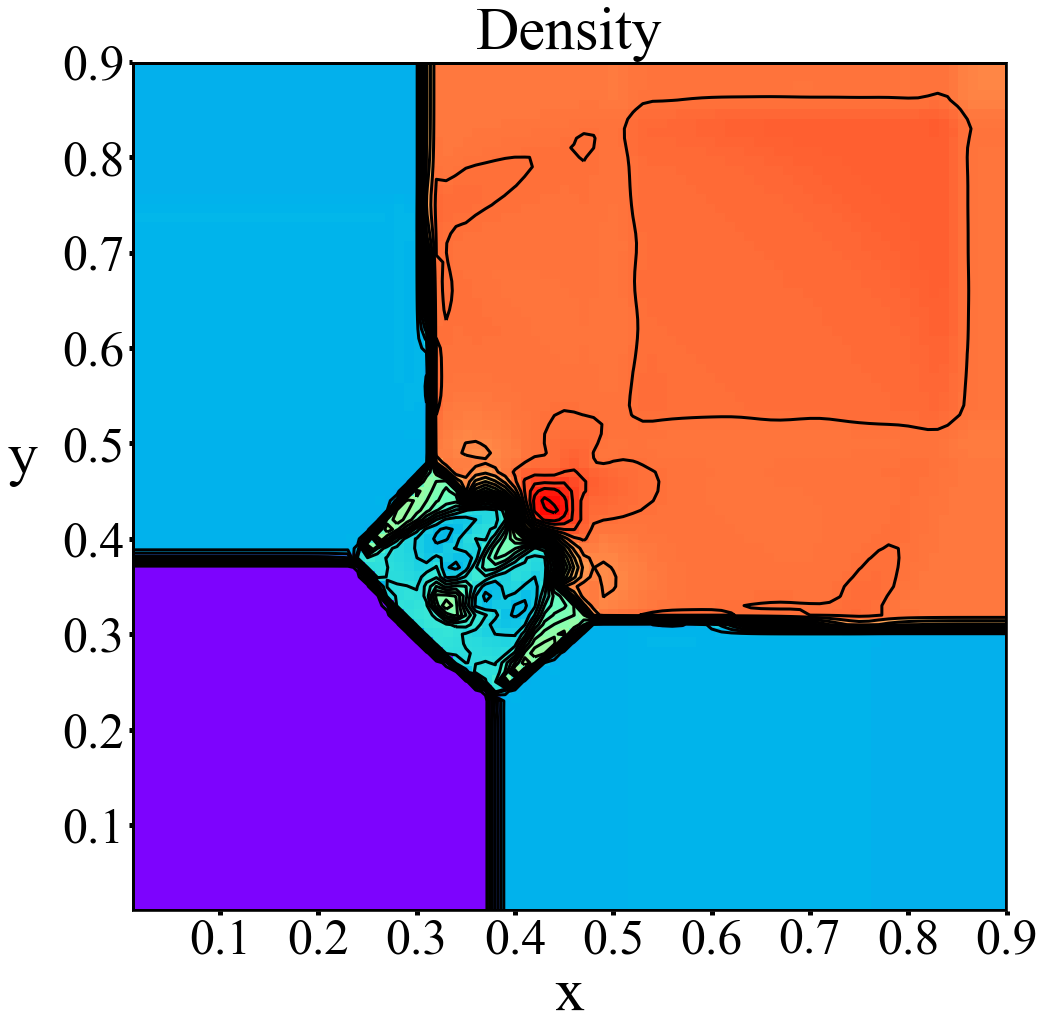}
    \caption{}
\end{subfigure}
\begin{subfigure}[b]{.48\textwidth}
    \centering
    \includegraphics[width=1.0\linewidth]{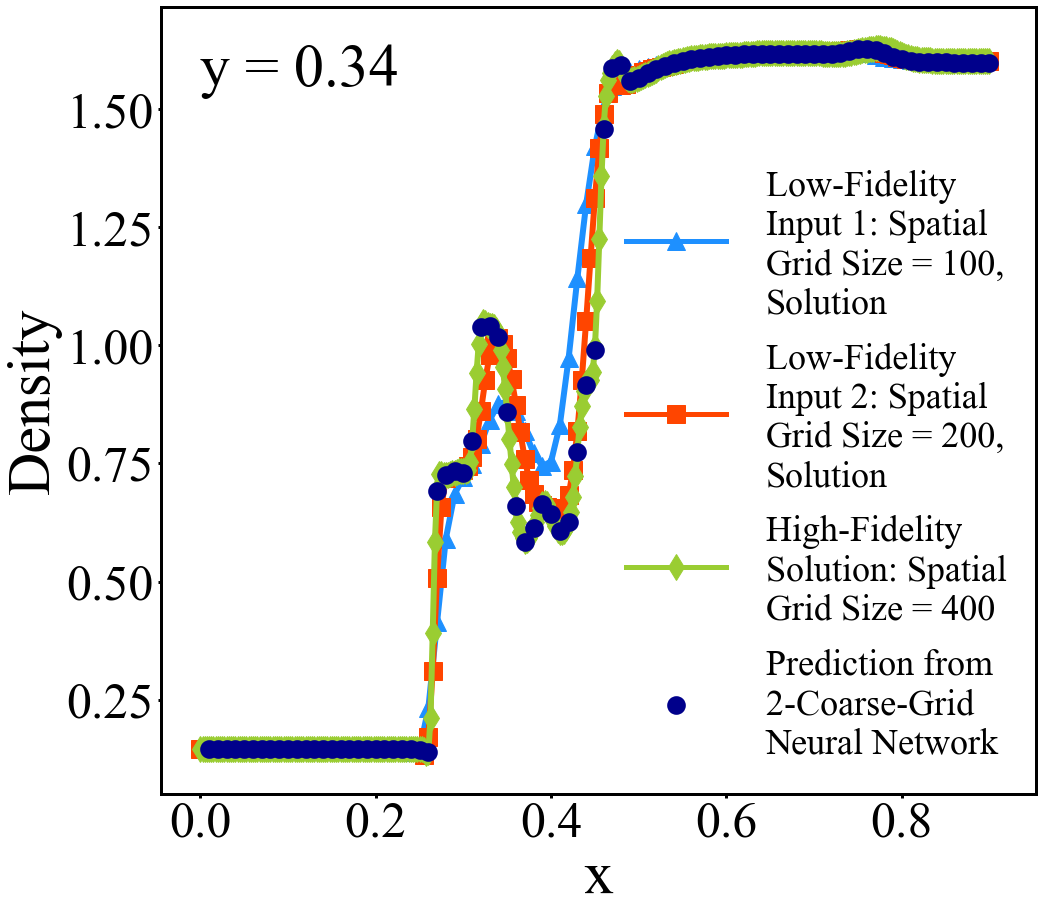}
    \caption{}
\end{subfigure}
\caption{2CGNN prediction of (a) the final-time $(t=0.3)$ density solution of the 2D Euler system with {\bf initial value $+5\%$ perturbation of that of Configuration 3}, and (b) its cross-section profile (dark blue) at \textit{y=}0.34, compared to low-fidelity input solutions (blue and red) by a $4$th order scheme on $100\times 100$ and $200\times 200$ grids, respectively, and ``exact'' (reference) solution (green).}
\label{2CGNN: Final time of config. 3 src input, +5}
\end{figure}

Fig.~\ref{2CGNN: Final time of config. 4, original} and ~\ref{2CGNN: Final time of config. 4, +5} show the predicted final-time solutions of the 2-D Euler system, with initial values the original initial value and $+5\%$ perturbation of that of Configuration 4, respectively. The spatial computational domain is $x,y\in [0.22, 0.98]$.

\begin{figure}[H] \centering
\begin{subfigure}[b]{.48\textwidth}
    \centering
    \includegraphics[width=1.0\linewidth]{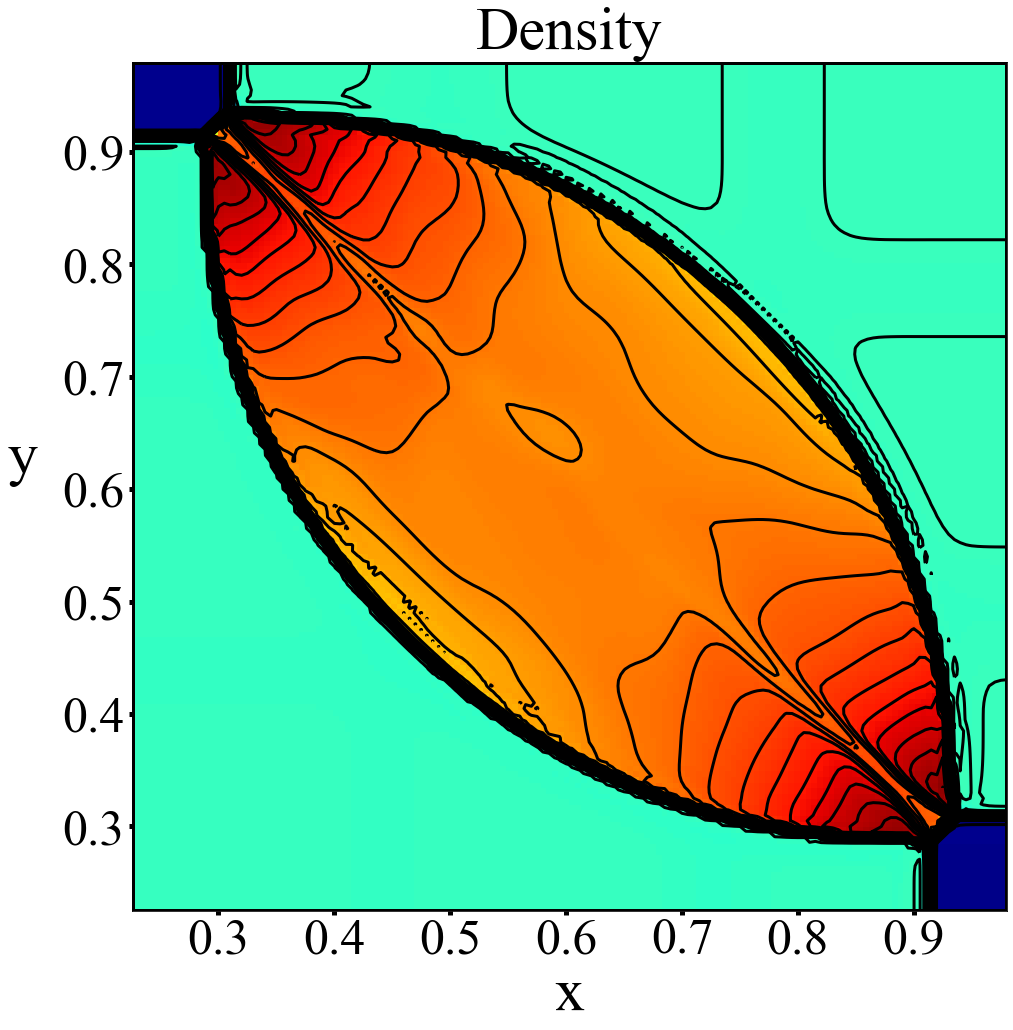}
    \caption{}
\end{subfigure}
\begin{subfigure}[b]{.48\textwidth}
    \centering
    \includegraphics[width=1.0\linewidth]{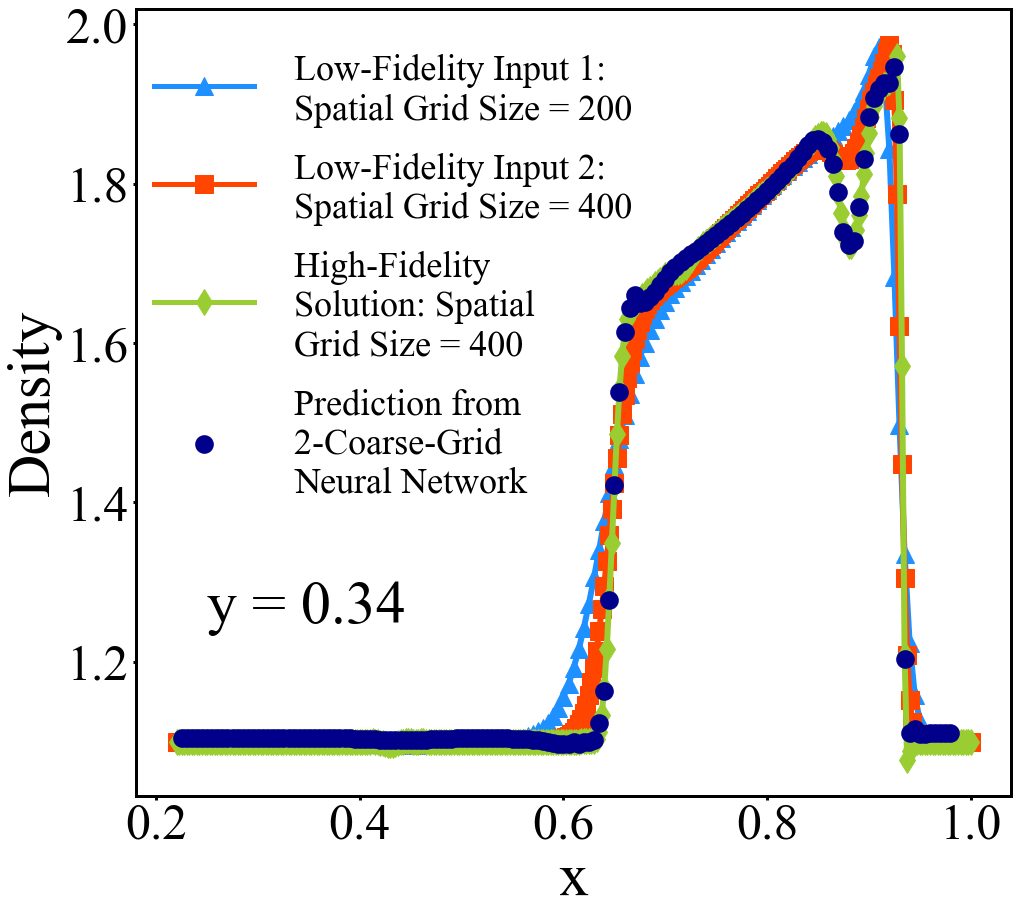}
    \caption{}
\end{subfigure}
\caption{2CGNN prediction of (a) the final-time $(t=0.25)$ density solution of \textbf{Configuration 4}, and (b) its cross-section profile (dark blue) along \textit{y=}0.34, compared to low-fidelity input solutions (blue and red) by  leapfrog and diffusion splitting scheme (\ref{leapfrog-diffusion-splitting, 2D}) on $200\times 200$ and $400\times 400$ grids, respectively, and ``exact'' (reference) solution (green).}
\label{2CGNN: Final time of config. 4, original}
\end{figure}

\begin{figure}[H]
\centering
\begin{subfigure}[b]{.48\textwidth}
    \centering
    \includegraphics[width=1.0\linewidth]{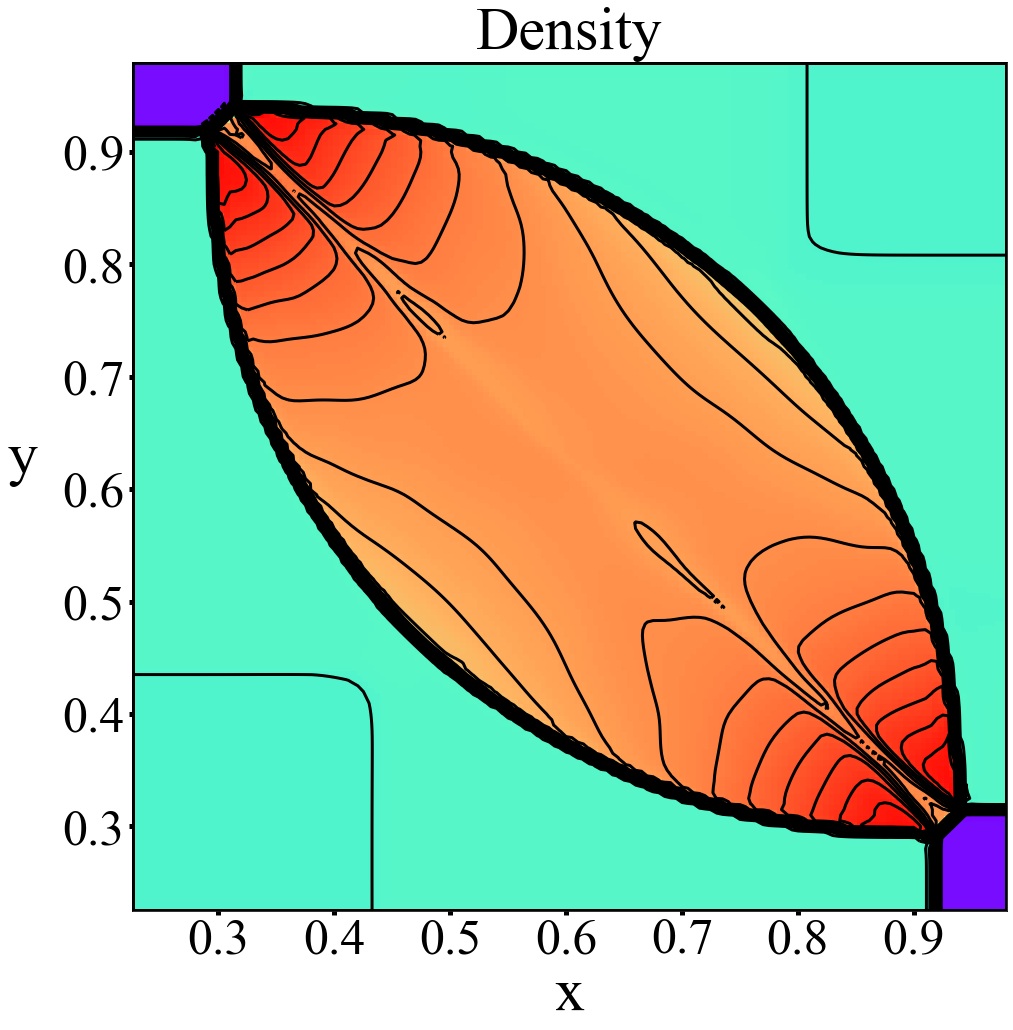}
    \caption{}
\end{subfigure}
\begin{subfigure}[b]{.48\textwidth}
    \centering
    \includegraphics[width=1.0\linewidth]{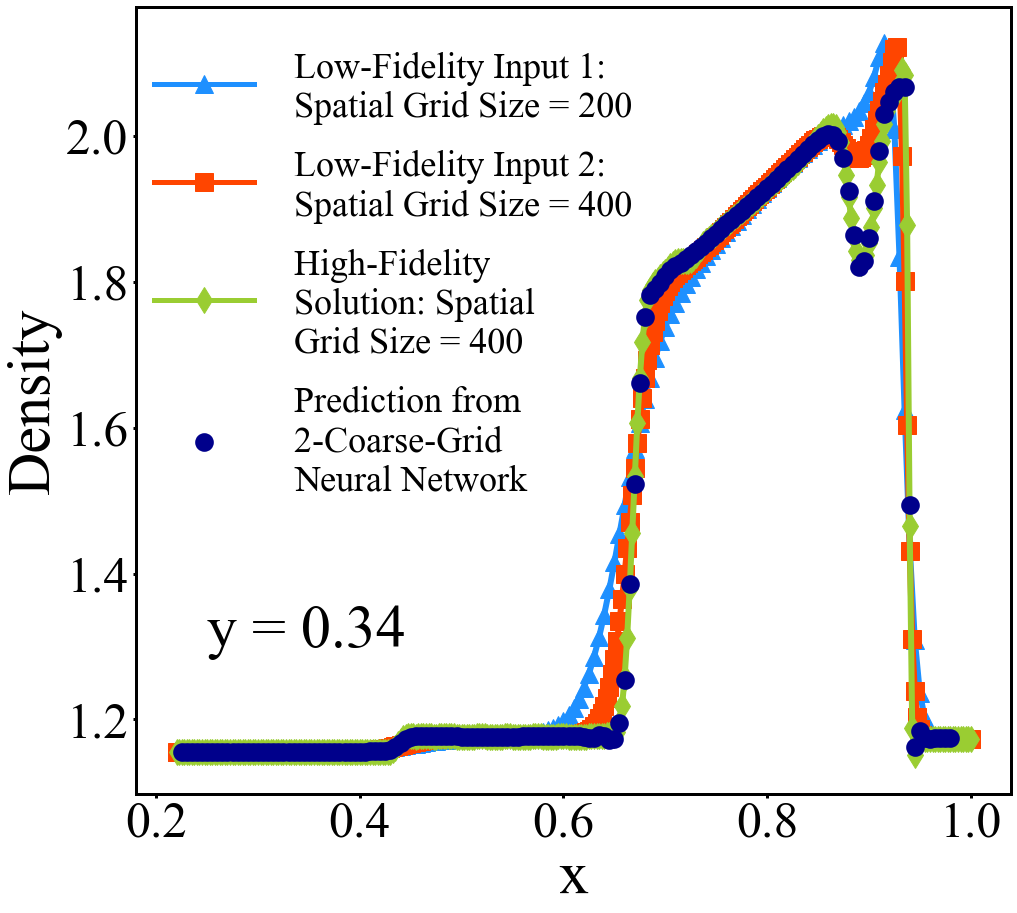}
    \caption{}
\end{subfigure}
\caption{2CGNN prediction of (a) the final-time $(t=0.25)$ density solution of the 2-D Euler system with {\bf initial value $+5\%$ perturbation of that of Configuration 4}, and (b) its cross-section profile (dark blue) along \textit{y=}0.34, compared to low-fidelity input solutions (blue and red) by leapfrog and diffusion splitting scheme (\ref{leapfrog-diffusion-splitting, 2D}) on $200\times 200$ and $400\times 400$ grids, respectively, and ``exact'' (reference) solution (green).}
\label{2CGNN: Final time of config. 4, +5}
\end{figure}

Fig.~\ref{2CGNN: Final time of config. 8, original} and ~\ref{2CGNN: Final time of config. 8, +5} show the predicted final-time solutions of the 2-D Euler system, with initial values the original initial value and $+5\%$ perturbation of that of Configuration 8, respectively. The spatial computational domain is $x,y\in [0, 0.9]$.

\begin{figure}[H] \centering
\begin{subfigure}[b]{.48\textwidth}
    \centering
    \includegraphics[width=1.0\linewidth]{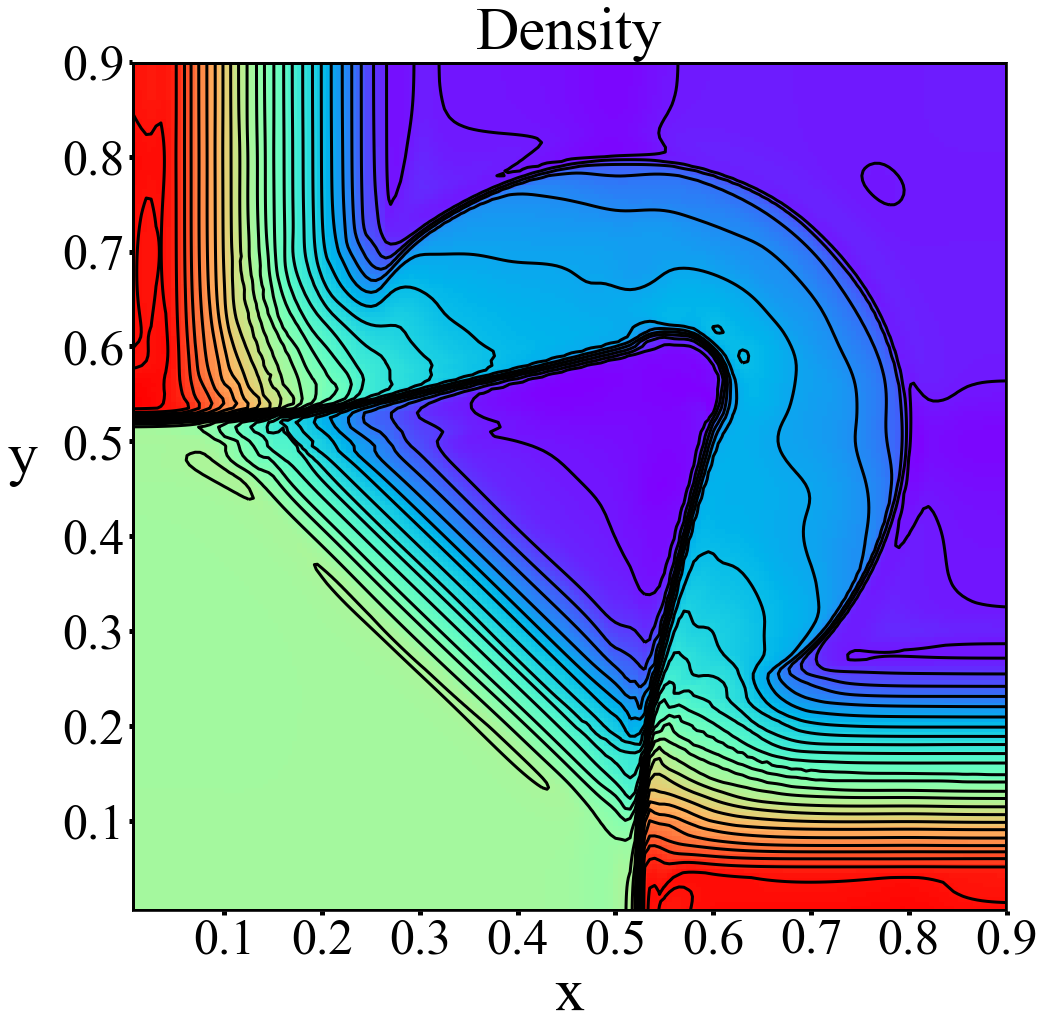}
    \caption{}
\end{subfigure}
\begin{subfigure}[b]{.48\textwidth}
    \centering
    \includegraphics[width=1.0\linewidth]{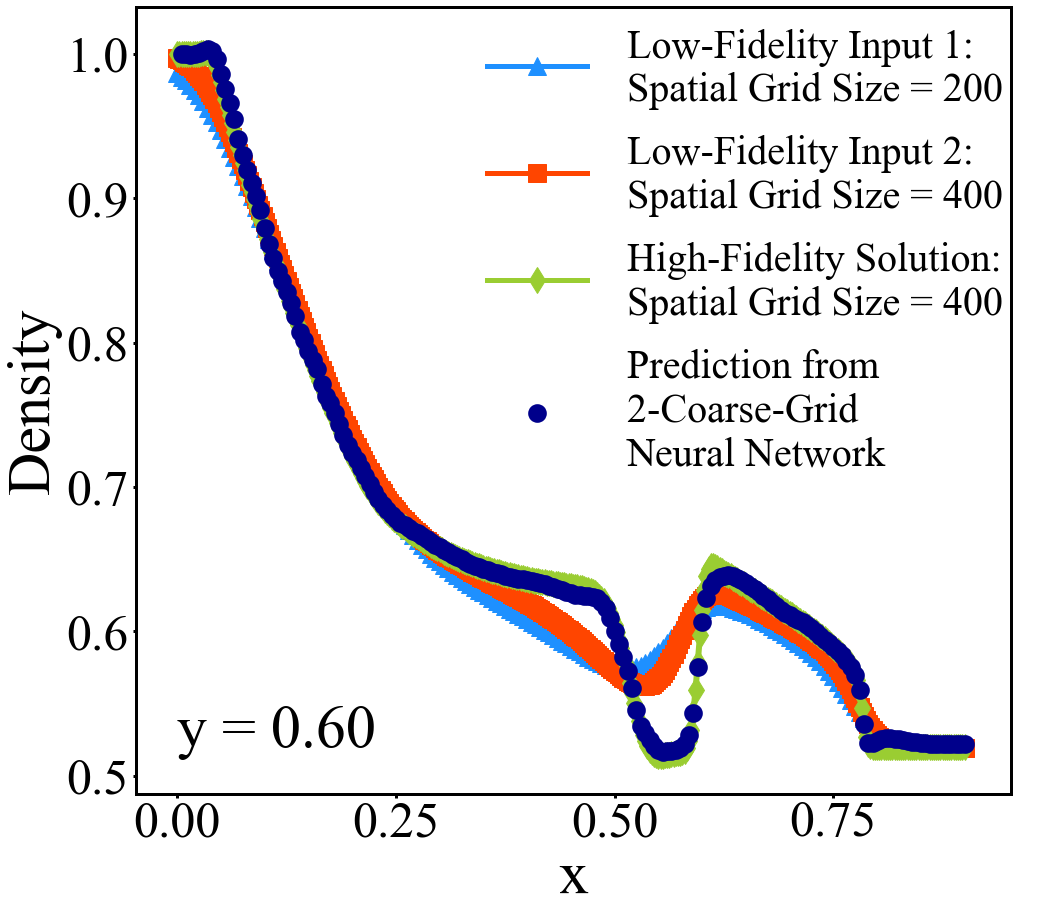}
    \caption{}
\end{subfigure}
\caption{2CGNN prediction of (a) the final-time $(t=0.25)$ density solution of \textbf{Configuration 8}, and (b) its cross-section profile (dark blue) along \textit{y=}0.60, compared to low-fidelity input solutions (blue and red) by leapfrog and diffusion splitting scheme (\ref{leapfrog-diffusion-splitting, 2D}) on $200\times 200$ and $400\times 400$ grids, respectively, and ``exact'' (reference) solution (green).}
\label{2CGNN: Final time of config. 8, original}
\end{figure}

\begin{figure}[H]
\centering
\begin{subfigure}[b]{.48\textwidth}
    \centering
    \includegraphics[width=1.0\linewidth]{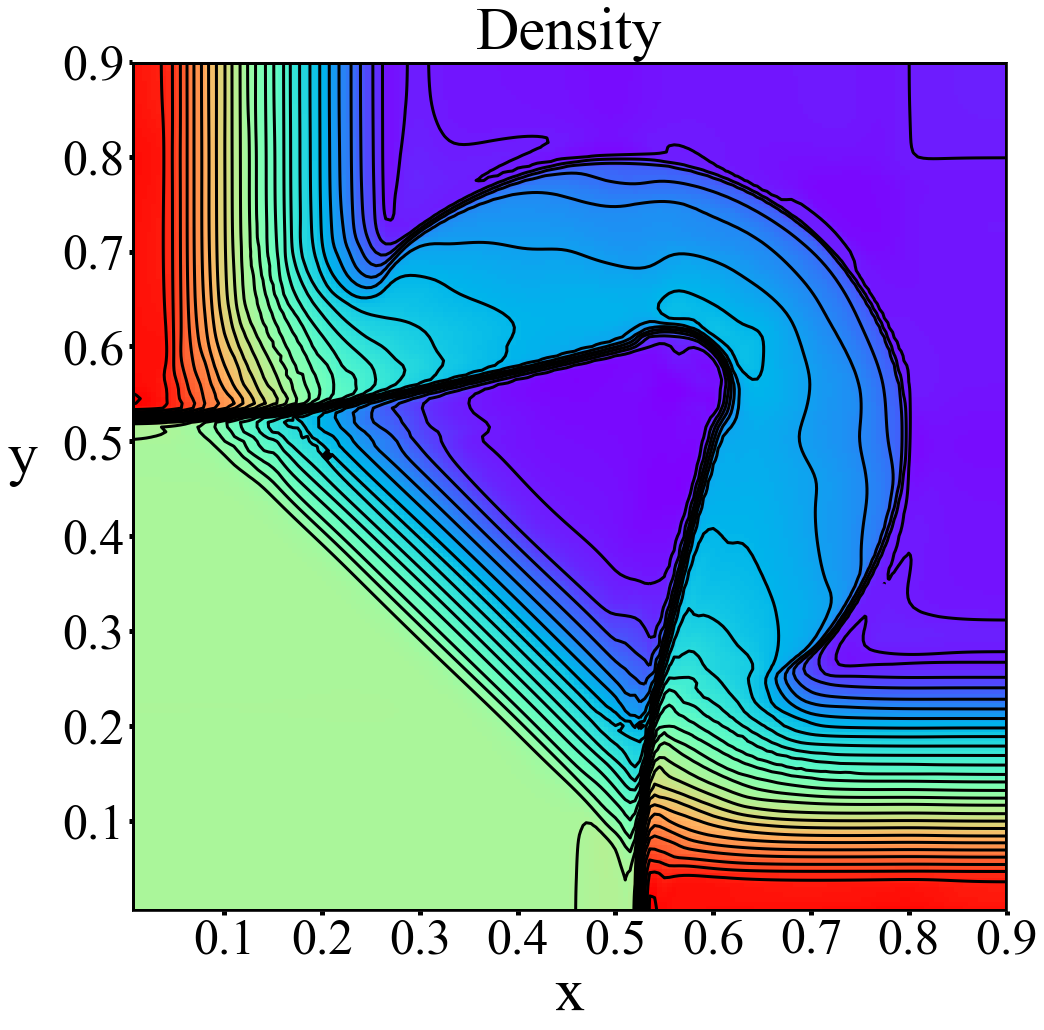}
    \caption{}
\end{subfigure}
\begin{subfigure}[b]{.48\textwidth}
    \centering
    \includegraphics[width=1.0\linewidth]{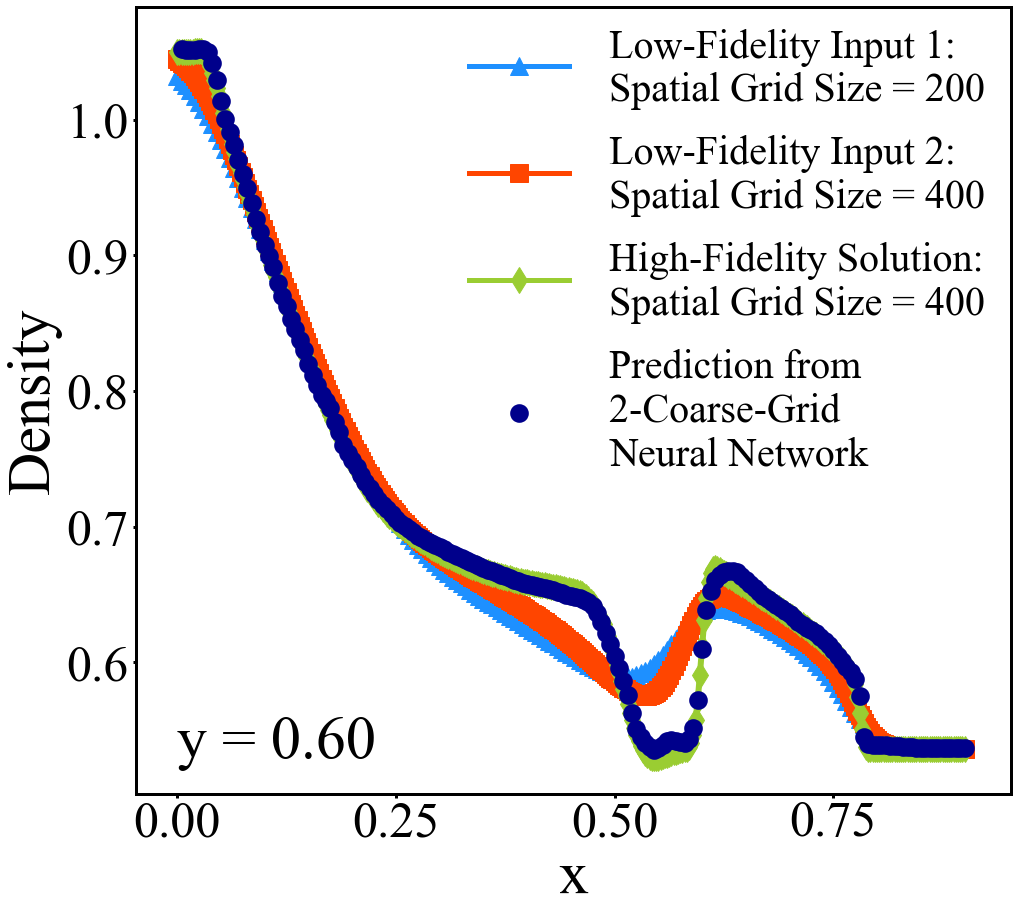}
    \caption{}
\end{subfigure}
\caption{2CGNN prediction of (a) the final-time $(t=0.25)$ density solution of the 2D Euler system with {\bf initial value $+5\%$ perturbation of that of Configuration 8}, and (b) its cross-section profile (dark blue) along \textit{y=}0.60, compared to low-fidelity input solutions (blue and red) by leapfrog and diffusion splitting scheme (\ref{leapfrog-diffusion-splitting, 2D}) on $200\times 200$ and $400\times 400$ grids, respectively, and ``exact'' (reference) solution (green).}
\label{2CGNN: Final time of config. 8, +5}
\end{figure}

\section{2-Diffusion-Coefficient Neural Network for 2D Problems}
\label{Sec: 2DCNN-2D}

Instead of computing for the input on two different grids, the input can be generated on a single grid using two equations perturbed with two different diffusion coefficients, using the vanishing viscosity approach~\cite{Neumann_Richtmyer50, PDLax57,Bianchini05}. This is referred to as a 2-Diffusion-Coefficient Neural Network (2DCNN) in \cite{2CGNN_1D_RiemannProb_21}. 
 
For example, we can approximate (\ref{cons-law}) using the leapfrog and diffusion splitting scheme (\ref{leapfrog-diffusion-splitting, 2D}) with $\alpha=\Delta x$ and $c\Delta x$ as follows
\begin{equation}
\label{leapfrog-diffusion-splitting-diff-coef-1}
\left \{
\begin{array}{l}
    \frac{\tilde{U}_{i,j} - U^{n-1}_{i,j}}{2\Delta t} +\frac{{f(U)}|^{n}_{i+1,j}-{f(U)}|^{n}_{i-1,j}}{2\Delta x} +\frac{{g(U)}|^{n}_{i,j+1}-{g(U)}|^{n}_{i,j-1}}{2\Delta y}=0~, \\
    \frac{U^{n+1}_{i,j} - \tilde{U}_{i,j}}{\Delta t} - \Delta x [\frac{\tilde{U}_{i+1,j} - 2\cdot \tilde{U}_{i,j} + \tilde{U}_{i-1,j}}{\Delta x^2} +  \frac{\tilde{U}_{i,j+1} - 2\cdot \tilde{U}_{i,j} + \tilde{U}_{i,j-1}}{\Delta y^2} ]=0~,
\end{array}
\right .
\end{equation}
and 
\begin{equation}
\label{leapfrog-diffusion-splitting-diff-coef-2}
\left \{
\begin{array}{l}
    \frac{\tilde{V}_{i,j} - V^{n-1}_{i,j}}{2\Delta t} +\frac{{f(V)}|^{n}_{i+1,j}-{f(V)}|^{n}_{i-1,j}}{2\Delta x} +\frac{{g(V)}|^{n}_{i,j+1}-{g(V)}|^{n}_{i,j-1}}{2\Delta y}=0~, \\
    \frac{V^{n+1}_{i,j} - \tilde{V}_{i,j}}{\Delta t} - c\Delta x [\frac{\tilde{V}_{i+1,j} - 2\cdot \tilde{V}_{i,j} + \tilde{V}_{i-1,j}}{\Delta x^2} + \frac{\tilde{V}_{i,j+1} - 2\cdot \tilde{V}_{i,j} + \tilde{V}_{i,j-1}}{\Delta y^2} ]=0~.
\end{array}
\right .
\end{equation}

Then the input computed on a single uniform grid can be written as

\begin{equation}
\begin{array}{cc}
\label{2DCNN-input}
    \{U^{n-1}_{i-1,j-1}, U^{n-1}_{i-1,j}, U^{n-1}_{i-1,j+1}, U^{n-1}_{i,j-1}, U^{n-1}_{i,j}, U^{n-1}_{i,j+1}, U^{n-1}_{i+1,j-1}, U^{n-1}_{i+1,j}, U^{n-1}_{i+1,j+1}, U^{n}_{i,j},\\ V^{n-1}_{i-1,j-1}, V^{n-1}_{i-1,j}, V^{n-1}_{i-1,j+1}, V^{n-1}_{i,j-1}, V^{n-1}_{i,j}, V^{n-1}_{i,j+1}, V^{n-1}_{i+1,j-1}, V^{n-1}_{i+1,j}, V^{n-1}_{i+1,j+1}, V^{n}_{i,j}\}~,
\end{array}
\end{equation}
as for the Euler system.
See Fig.~\ref{2DCNN: Final time of config. 8, dx, 4dx, original} and \ref{2DCNN: Final time of config. 8, dx, 4dx, +5} for the predictions by 2DCNN, and their cross-section profiles, with the input computed on a $400\times 400$ uniform grid and $c=4$. Note that for the second equation
(\ref{leapfrog-diffusion-splitting-diff-coef-2}), the larger diffusion coefficient may require a smaller time step size than that of the first equation of (\ref{leapfrog-diffusion-splitting-diff-coef-1}). In that case, the second equation of (\ref{leapfrog-diffusion-splitting-diff-coef-2}) will be computed in two time steps, with the time step size $\frac12\Delta t$ for each time step. 

\begin{figure}[H] \centering
\begin{subfigure}[b]{.48\textwidth}
    \centering
    \includegraphics[width=1.0\linewidth]{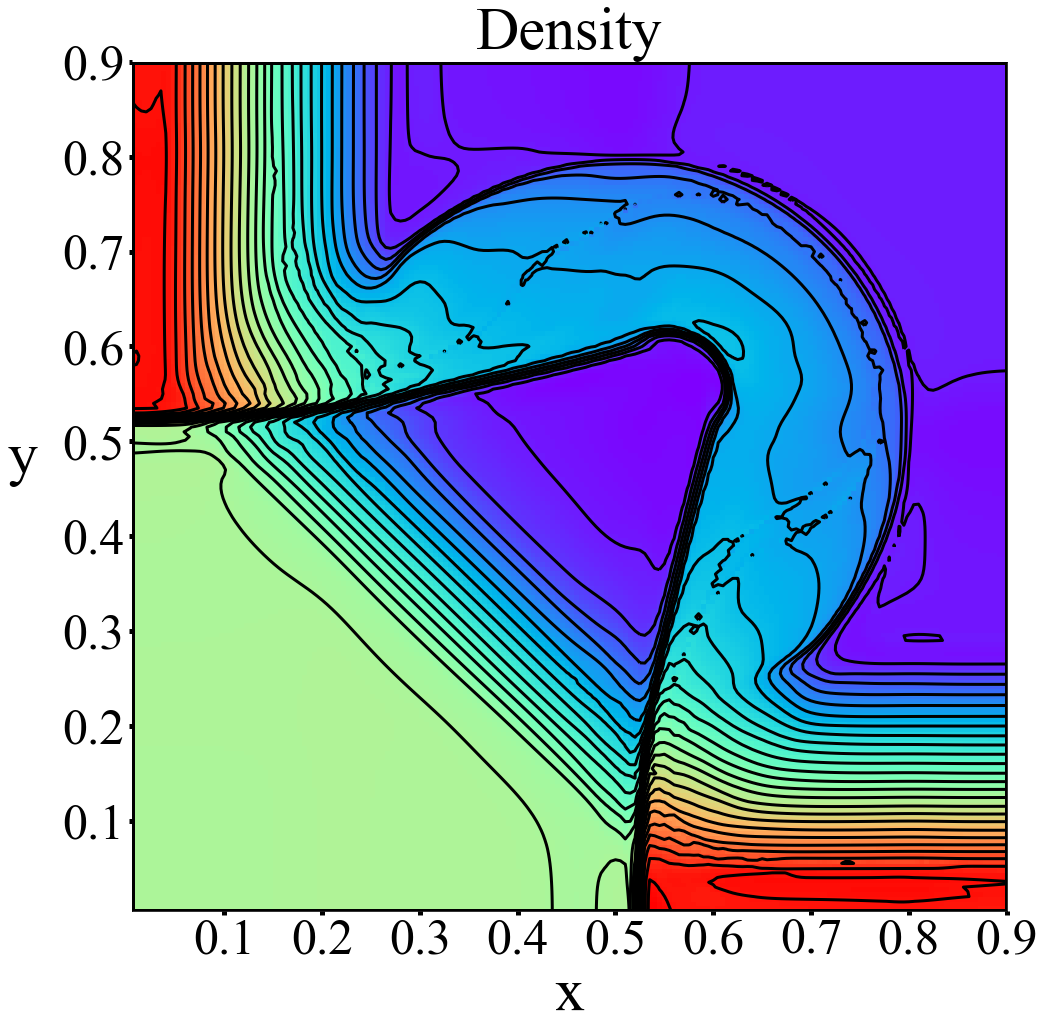}
    \caption{}
\end{subfigure}
\begin{subfigure}[b]{.48\textwidth}
    \centering
    \includegraphics[width=1.0\linewidth]{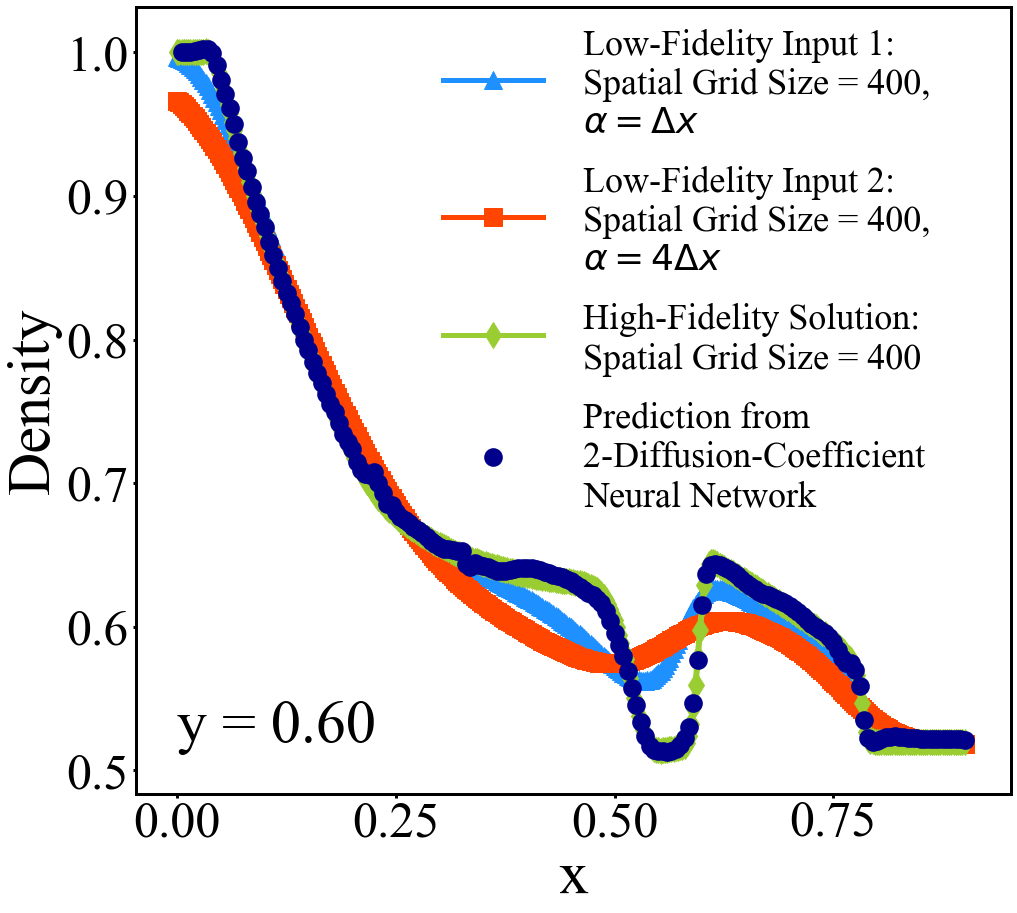}
    \caption{}
\end{subfigure}
\caption{2DCNN prediction of (a) the final-time $(t=0.25)$ density solution of \textbf{Configuration 8}, and (b) its cross-section profile (dark blue) along \textit{y=}0.60, compared to low-fidelity input solutions (blue and red) by leapfrog and diffusion splitting schemes (\ref{leapfrog-diffusion-splitting-diff-coef-1}) and (\ref{leapfrog-diffusion-splitting-diff-coef-2}) ($c=4$), respectively on the $400 \times 400$ grid, and ``exact'' (reference) solution (green).}
\label{2DCNN: Final time of config. 8, dx, 4dx, original}
\end{figure}

\begin{figure}[H]
\centering
\begin{subfigure}[b]{.48\textwidth}
    \centering
    \includegraphics[width=1.0\linewidth]{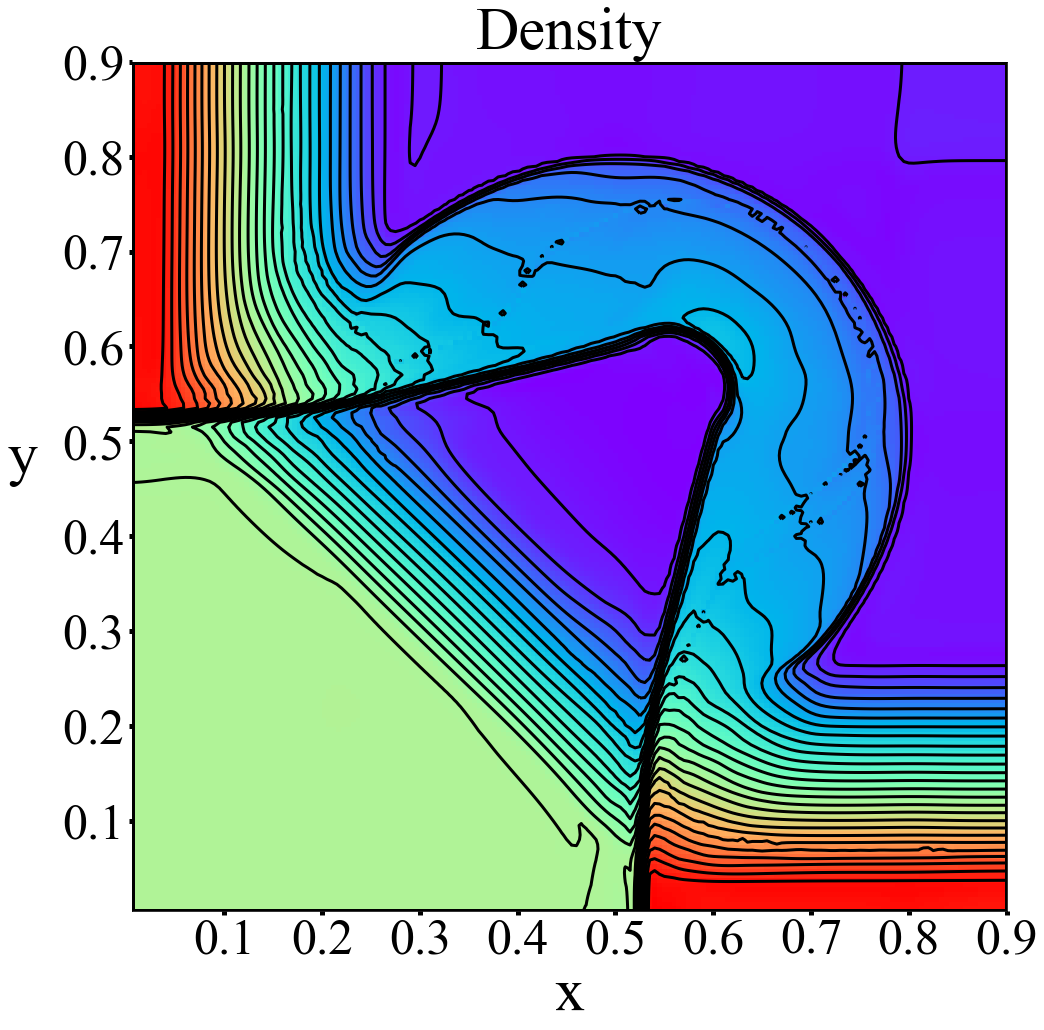}
    \caption{}
\end{subfigure}
\begin{subfigure}[b]{.48\textwidth}
    \centering
    \includegraphics[width=1.0\linewidth]{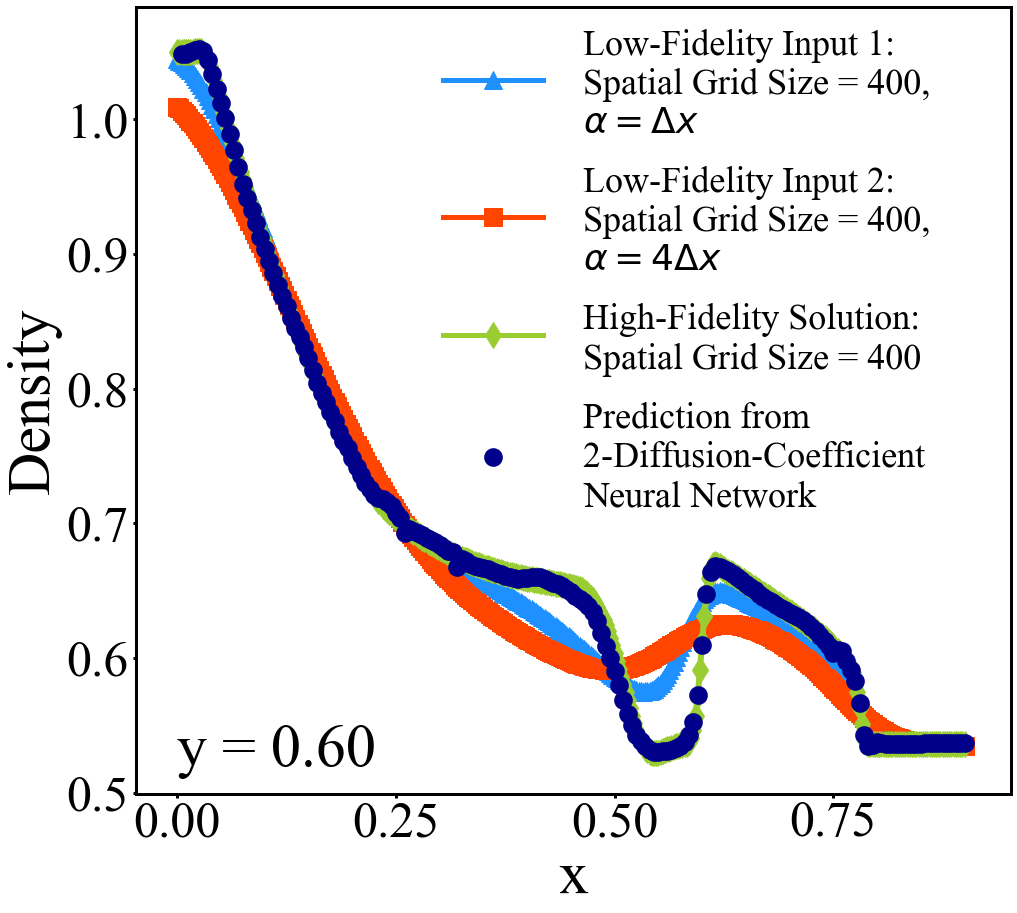}
    \caption{}
\end{subfigure}
\caption{2DCNN prediction of (a) the final-time $(t=0.25)$ density solution of the 2D Euler system with {\bf initial value $+5\%$ perturbation of that of Configuration 8}, and (b) its cross-section profile (dark blue) along \textit{y=}0.60, compared to low-fidelity input solutions (blue and red) by leapfrog and diffusion splitting schemes (\ref{leapfrog-diffusion-splitting-diff-coef-1}) and (\ref{leapfrog-diffusion-splitting-diff-coef-2}) ($c=4$), respectively on the $400 \times 400$ grid, and ``exact'' (reference) solution (green).}
\label{2DCNN: Final time of config. 8, dx, 4dx, +5}
\end{figure}

 We can also compute the two parts of the input by a first order scheme and a high order scheme, respectively, as in Sec.~\ref{low-high-2DCNN-1D}. For example, the first part can be computed by the leapfrog and diffusion splitting scheme (\ref{leapfrog-diffusion-splitting, 2D}) on $400\times 400$ grid with $\alpha=4\Delta x$, the second part by the $4$th order scheme used for computing reference solutions on $200\times 200$ grid, and the input is formatted on the $200\times 200$ grid. Figures \ref{2DCNN: Final time of config. 8, 4dx, src_input, original} and \ref{2DCNN: Final time of config. 8, 4dx, src_input, +5} show the prediction results from this approach.

\begin{figure}[H] \centering
\begin{subfigure}[b]{.48\textwidth}
    \centering
    \includegraphics[width=1.0\linewidth]{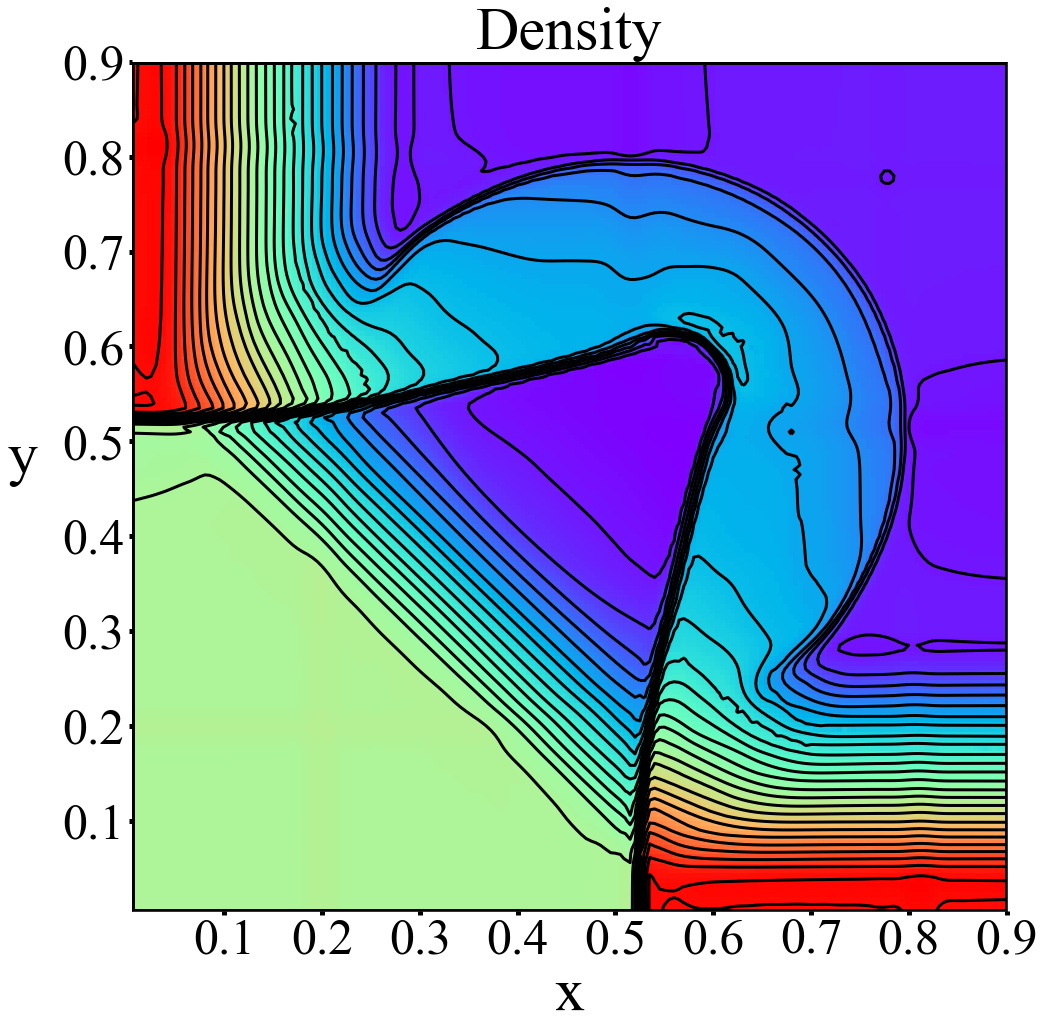}
    \caption{}
\end{subfigure}
\begin{subfigure}[b]{.48\textwidth}
    \centering
    \includegraphics[width=1.0\linewidth]{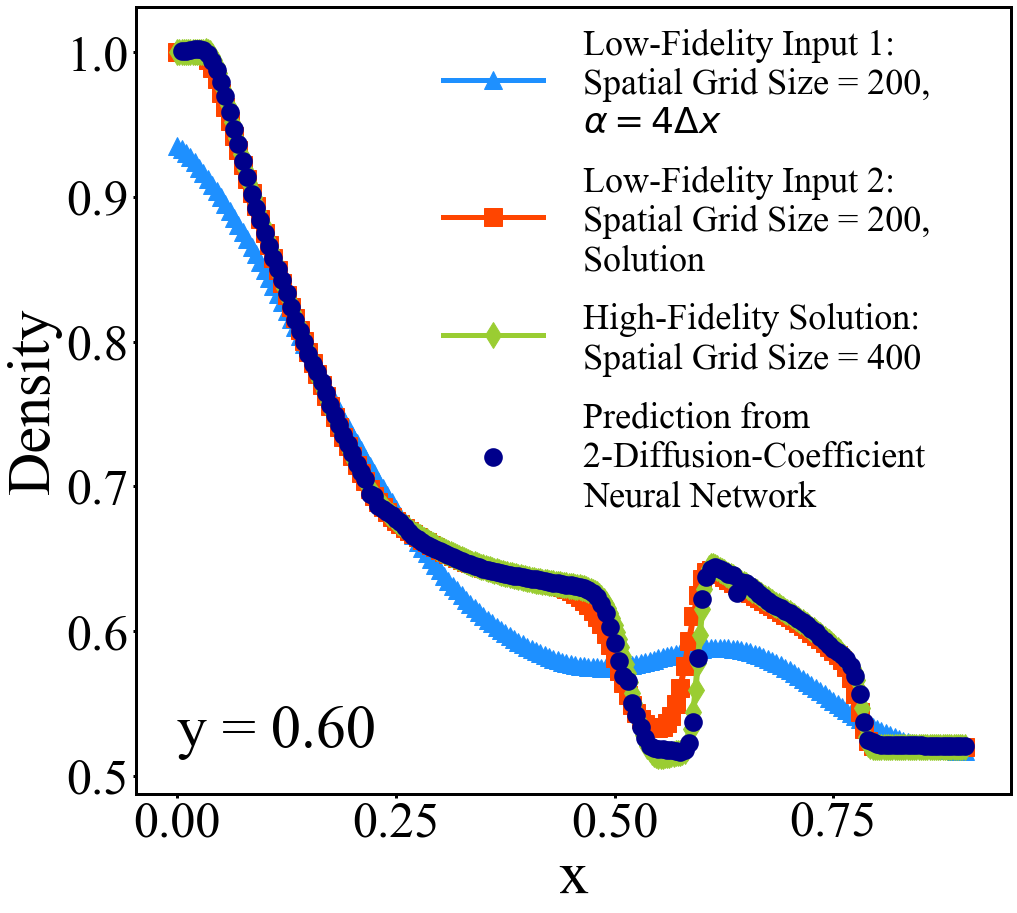}
    \caption{}
\end{subfigure}
\caption{2DCNN prediction of (a) final-time $(t=0.25)$ density solution of \textbf{Configuration 8}, and (b) its cross-section profile (dark blue) along \textit{y=}0.60, compared to low-fidelity input solutions (blue and red) by leapfrog and diffusion splitting scheme (\ref{leapfrog-diffusion-splitting, 2D}) ($\alpha = 4\Delta x$) on $400 \times 400$ grid, and a $4$th order scheme on $200 \times 200$ grid, respectively, and ``exact'' (reference) solution (green).}
\label{2DCNN: Final time of config. 8, 4dx, src_input, original}
\end{figure}

\begin{figure}[H]
\centering
\begin{subfigure}[b]{.48\textwidth}
    \centering
    \includegraphics[width=1.0\linewidth]{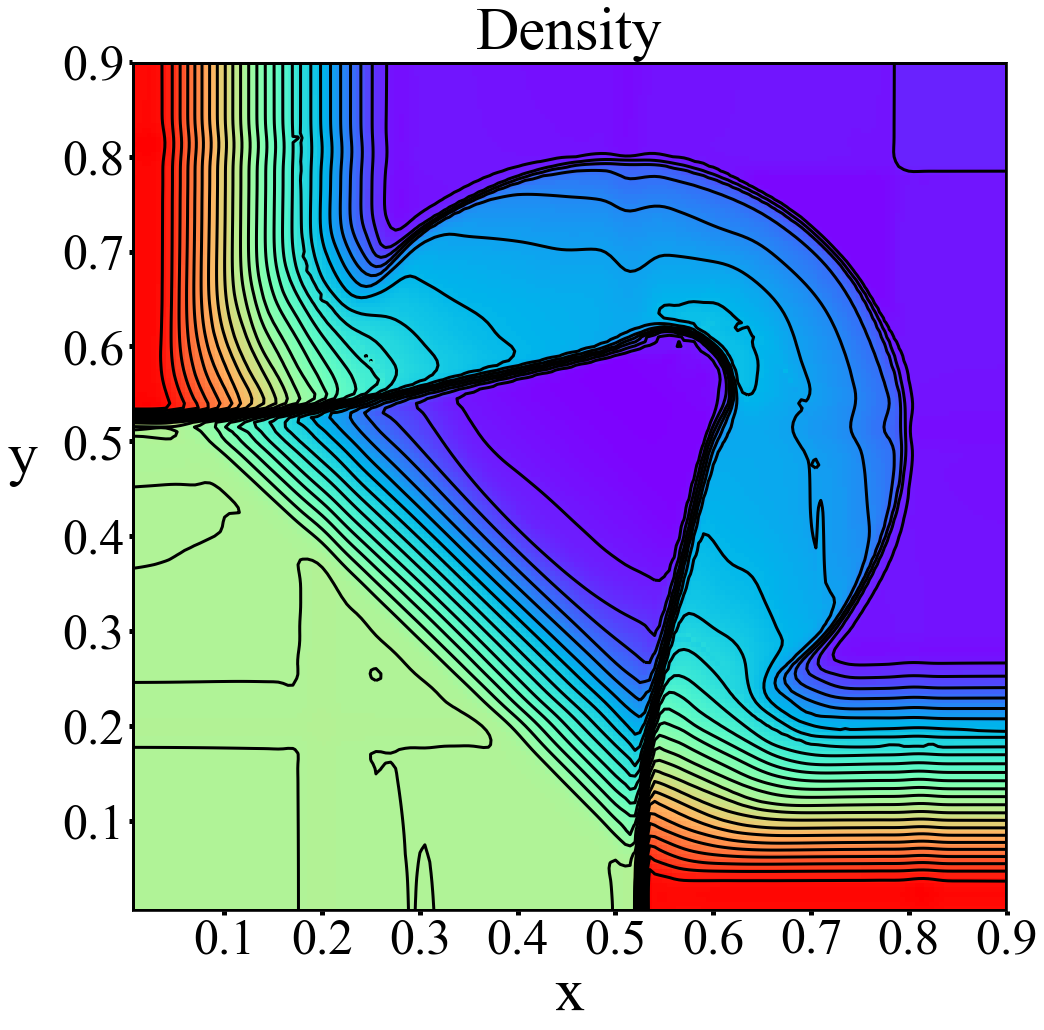}
    \caption{}
\end{subfigure}
\begin{subfigure}[b]{.48\textwidth}
    \centering
    \includegraphics[width=1.0\linewidth]{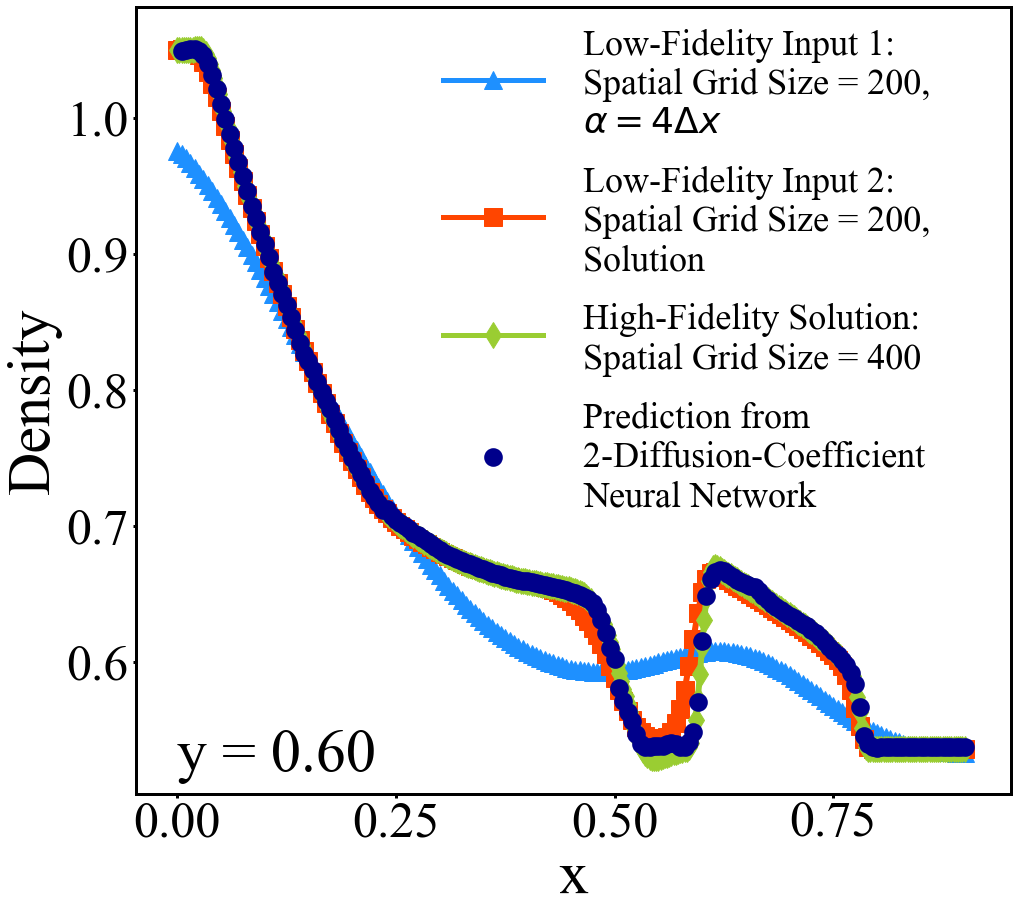}
    \caption{}
\end{subfigure}
\caption{2DCNN prediction of (a) the final-time $(t=0.25)$ density solution of the 2D Euler system with {\bf initial value $+5\%$ perturbation of that of Configuration 8}, and (b) its cross-section profile (dark blue) along \textit{y=}0.60, compared to low-fidelity input solutions (blue and red) by leapfrog and diffusion splitting scheme (\ref{leapfrog-diffusion-splitting, 2D}) ($\alpha = 4\Delta x$) on $400 \times 400$ grid, and a $4$th order scheme on $200 \times 200$ grid, respectively, and ``exact'' (reference) solution (green).}
\label{2DCNN: Final time of config. 8, 4dx, src_input, +5}
\end{figure}

We can also improve the first part of the input slightly by computing it by the leapfrog and diffusion splitting scheme (\ref{leapfrog-diffusion-splitting, 2D}) on $400\times 400$ grid with $\alpha=\Delta x$.
Figures \ref{2DCNN: Final time of config. 8, dx, src_input, original} and \ref{2DCNN: Final time of config. 8, dx, src_input, +5} show improved predictions from this minor adjustment.

\begin{figure}[H] \centering
\begin{subfigure}[b]{.48\textwidth}
    \centering
    \includegraphics[width=1.0\linewidth]{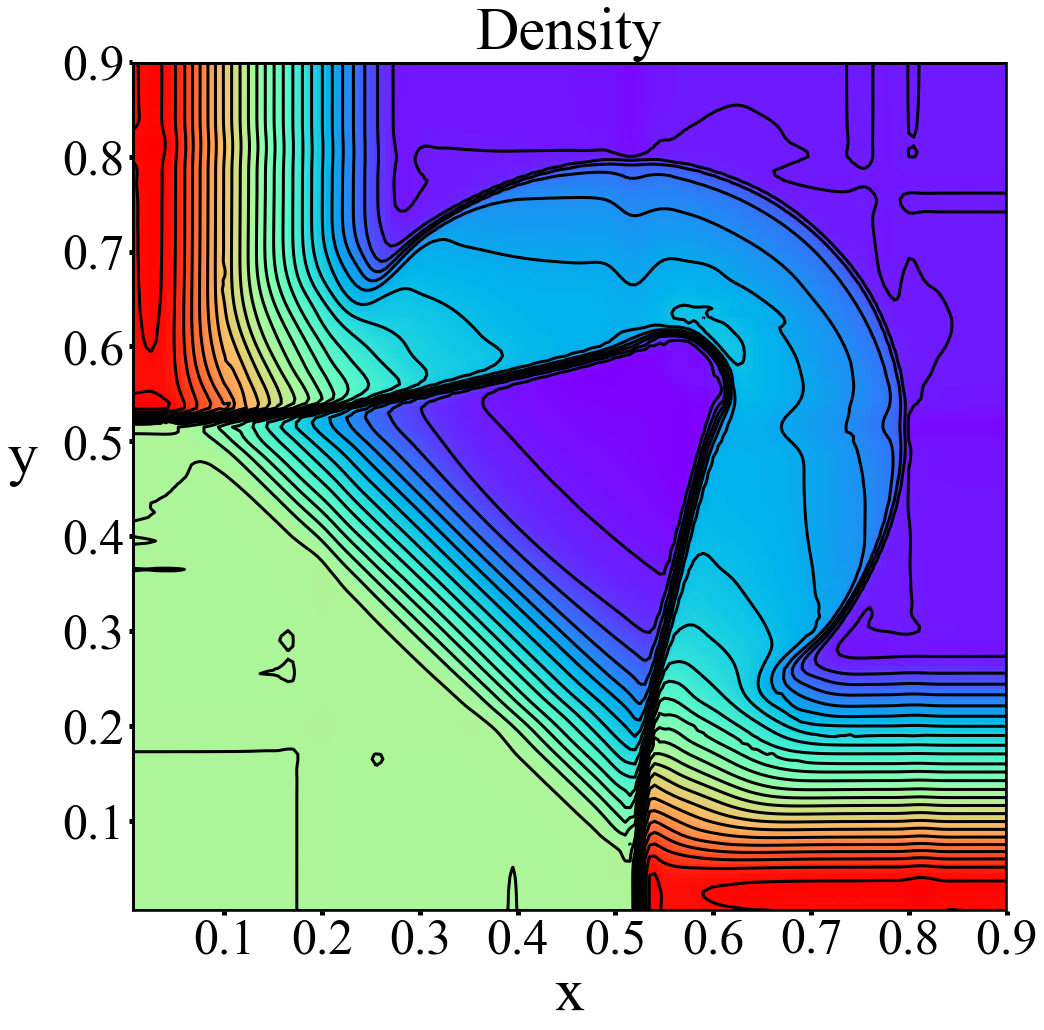}
    \caption{}
\end{subfigure}
\begin{subfigure}[b]{.48\textwidth}
    \centering
    \includegraphics[width=1.0\linewidth]{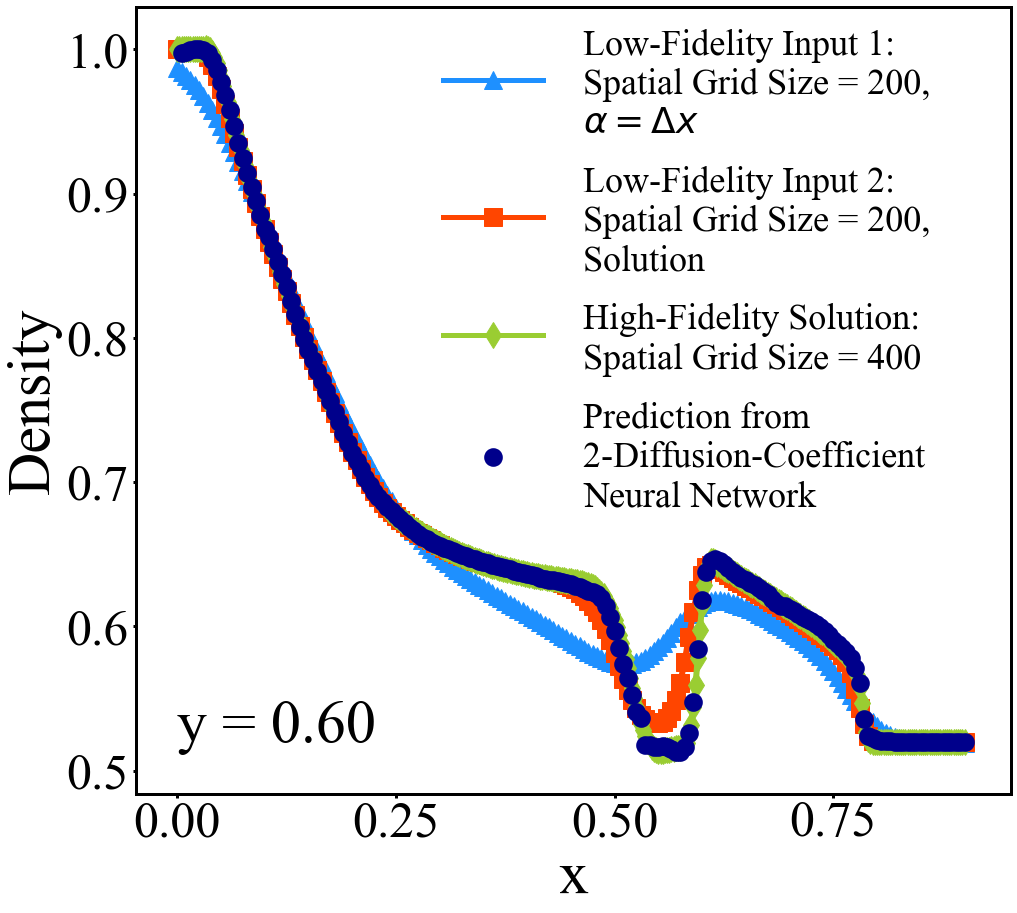}
    \caption{}
\end{subfigure}
\caption{2DCNN prediction of (a) final-time $(t=0.25)$ density solution of \textbf{Configuration 8}, and (b) its cross-section profile (dark blue) along \textit{y=}0.60, compared to low-fidelity input solutions (blue and red) by leapfrog and diffusion splitting scheme (\ref{leapfrog-diffusion-splitting, 2D}) ($\alpha = \Delta x$) on $400 \times 400$ grid, and a $4$th order scheme on $200 \times 200$ grid, respectively, and ``exact'' (reference) solution (green).}
\label{2DCNN: Final time of config. 8, dx, src_input, original}
\end{figure}

\begin{figure}[H]
\centering
\begin{subfigure}[b]{.48\textwidth}
    \centering
    \includegraphics[width=1.0\linewidth]{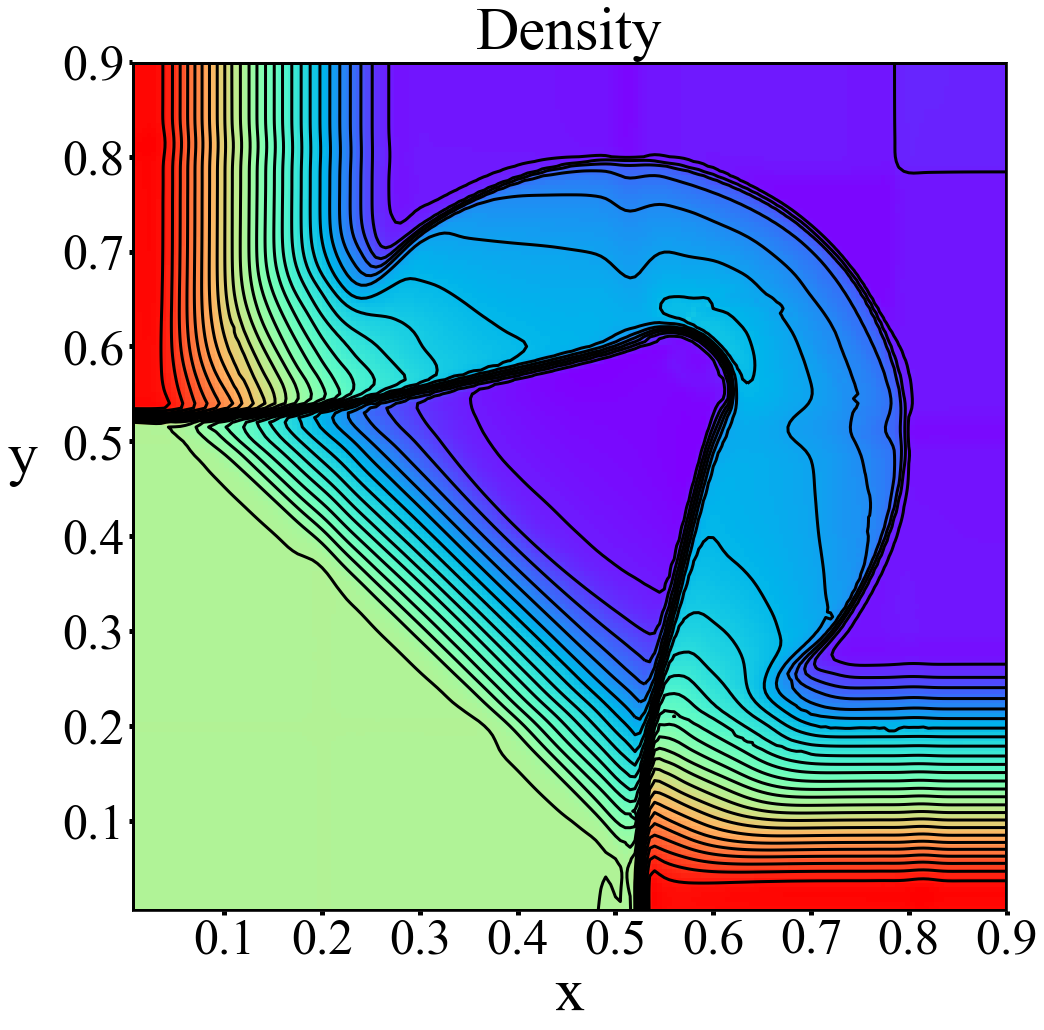}
    \caption{}
\end{subfigure}
\begin{subfigure}[b]{.48\textwidth}
    \centering
    \includegraphics[width=1.0\linewidth]{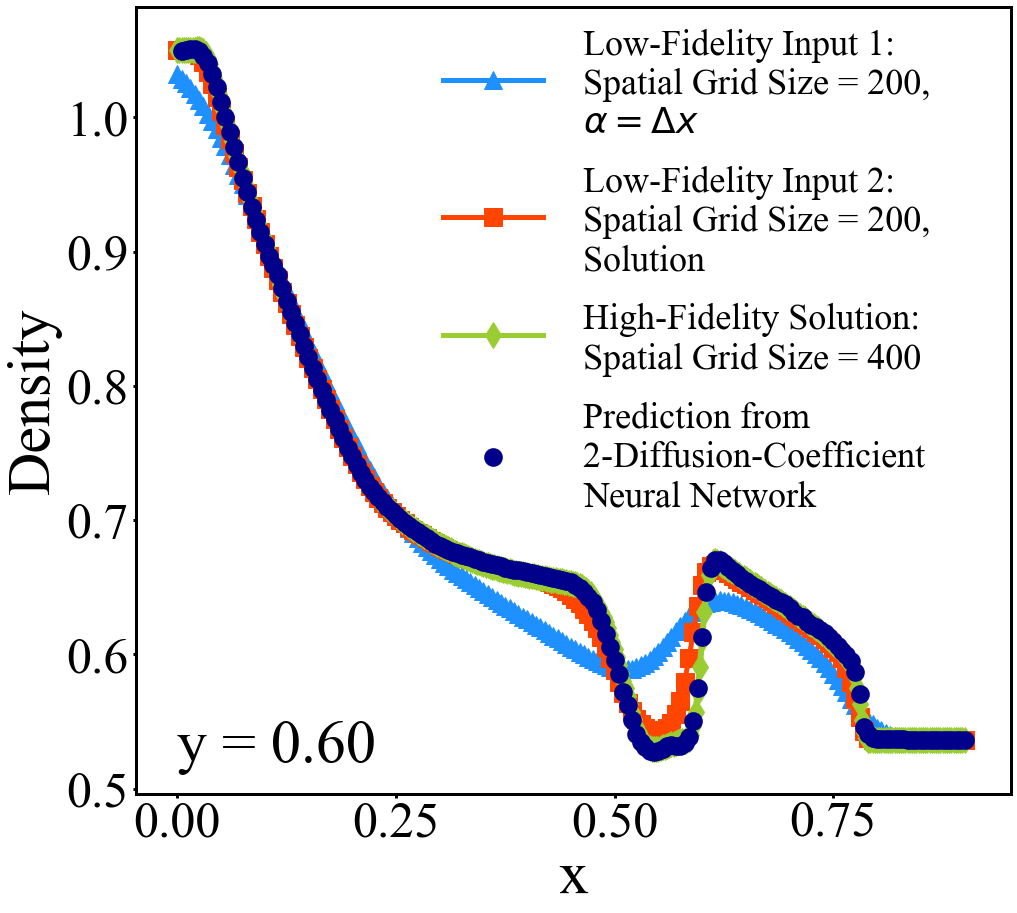}
    \caption{}
\end{subfigure}
\caption{2DCNN prediction of (a) the final-time $(t=0.25)$ density solution of the 2D Euler system with {\bf initial value $+5\%$ perturbation of that of Configuration 8}, and (b) its cross-section profile (dark blue) along \textit{y=}0.60, compared to low-fidelity input solutions (blue and red) by leapfrog and diffusion splitting scheme (\ref{leapfrog-diffusion-splitting, 2D}) ($\alpha = \Delta x$) on $400 \times 400$ grid, and a $4$th order scheme on $200 \times 200$ grid, respectively, and ``exact'' (reference) solution (green).}
\label{2DCNN: Final time of config. 8, dx, src_input, +5}
\end{figure}

\section{Summary}
\label{Sec: summary}
The tables below summarize the results of the methods tested. 

\begin{table}[H]\centering
\caption{Comparison of cross-training in Sec.~\ref{Cross-training} with 3 different types of training or inputs: Relative $l_2$ norm errors between the predicted and ``exact'' (reference) solutions for the Lax problem.}\label{tb: Summary of cross-training errors, Lax Prob. }

\begin{tabular}{|c|c|c|c|} \hline
& \multicolumn{3}{|c|}{\textbf{Lax Problem}}\\ \hline
& & $\Delta t$ in Input and & \\ 
& $\Delta t$ in Input & 2 Intermediate & $\Delta t$ not in Input \\ 
&  & $\Delta t$ values in training & \\ \hline
Initial Value & Relative $l_2$ Errors & Relative $l_2$ Errors & Relative $l_2$ Errors \\ \hline
Original & $9.36e-03$ & $6.95e-03$ & $1.22e-03$ \\ \hline
$+3\%$ & $4.29e-03$ & $3.65e-03$ & $1.23e-03$ \\ \hline
$-3\%$ & $8.89e-03$ & $6.01e-03$ & $1.49e-03$ \\ \hline
$+5\%$ & $8.46e-03$ & $4.69e-03$ & $1.29e-03$ \\ \hline
$-5\%$ & $9.35e-03$ & $8.13e-03$ & $1.40e-03$ \\ \hline
$+7\%$ & $9.04e-03$ & $7.01e-03$ & $1.27e-03$ \\ \hline
$-7\%$ & $1.11e-02$ & $8.81e-03$ & $1.21e-03$ \\ \hline
\end{tabular}
\end{table}

\begin{table}[H]\centering
\caption{Comparison of cross-training in Sec.~\ref{Cross-training} with 3 different types of training or inputs: Relative $l_2$ norm errors between the predicted and ``exact'' (reference) solutions for the Sod problem.}\label{tb: Summary of cross-training errors, Sod Prob. }
\begin{tabular}{|c|c|c|c|} \hline
& \multicolumn{3}{|c|}{\textbf{Sod Problem}}\\ \hline
& & $\Delta t$ in Input and & \\ 
& $\Delta t$ in Input & 2 Intermediate & $\Delta t$ not in Input \\ 
&  & $\Delta t$ values in training & \\ \hline
Initial Value & Relative $l_2$ Errors & Relative $l_2$ Errors & Relative $l_2$ Errors \\ \hline
Original & $3.31e-03$ & $5.36e-03$ & $6.45e-03$ \\ \hline
$+3\%$ & $3.03e-03$ & $4.93e-03$ & $6.39e-03$ \\ \hline
$-3\%$ & $4.63e-03$ & $7.35e-03$ & $6.54e-03$ \\ \hline
$+5\%$ & $3.45e-03$ & $6.43e-03$ & $6.38e-03$ \\ \hline
$-5\%$ & $5.59e-03$ & $8.21e-03$ & $6.63e-03$ \\ \hline
$+7\%$ & $4.19e-03$ & $6.59e-03$ & $6.42e-03$ \\ \hline
$-7\%$ & $6.62e-03$ & $9.29e-03$ & $6.75e-03$ \\ \hline
\end{tabular}
\end{table}

\begin{table}[H]
\centering
\caption{Relative $l_2$ errors of 2CGNN predictions of solutions of the Euler system with initial values perturbed from those of the W-C problem,  with low-fidelity input solutions computed by a $3$rd order scheme on $200\times 200$ and $400\times 400$ grids (Sec.~\ref{2CGNN-4-WC}.)}\label{tb: Summary of total errors from predictions, 2CGNN WC prob.}
\begin{tabular}{|c|c|} \hline
&\textbf{W-C Problem} \\ \hline
Initial Value &Relative $l_2$ Error \\ \hline
Original &$1.29e-04$  \\ \hline
$+3\%$ &$8.22e-05$ \\ \hline
$-3\%$ &$8.38e-05$ \\ \hline
$+5\%$ &$9.43e-05$ \\ \hline
$-5\%$ &$8.63e-05$ \\ \hline
$+7\%$ &$1.23e-04$ \\ \hline
$-7\%$ &$9.10e-05$ \\ \hline
\end{tabular}
\label{2CGNN Src. input WC results}
\end{table}

\begin{table}[H]
\centering
\caption{Relative $l_2$ errors of 2DCNN predictions of solutions of the Euler system with initial values perturbed from those of the W-C problem,  with low-fidelity input solutions computed by the Rusanov scheme and a $3$rd order scheme on $400\times 400$ grid (Sec.~\ref{low-high-2DCNN-1D}.)}\label{tb: Summary of total errors from predictions, 2DCNN WC prob.}
\begin{tabular}{|c|c|} \hline
&\textbf{WC Problem} \\ \hline
Initial Value &Relative $l_2$ Error \\ \hline
Original &$2.56e-04$  \\ \hline
$+3\%$ &$1.41e-04$ \\ \hline
$-3\%$ &$1.45e-04$ \\ \hline
$+5\%$ &$1.33e-04$ \\ \hline
$-5\%$ &$1.44e-04$ \\ \hline
$+7\%$ &$1.32e-04$ \\ \hline
$-7\%$ &$1.44e-04$ \\ \hline
\end{tabular}
\label{2DCNN Rusanov Src. input WC results}
\end{table}

\begin{table}[H]\centering
\caption{Relative $l_2$ errors of 2CGNN predictions of solutions of the 2D Euler system with initial values perturbed from those of 2D Riemann problems, with low-fidelity input solutions computed by leapfrog and diffusion splitting scheme (\ref{leapfrog-diffusion-splitting, 2D}) on $200\times 200$ and $400\times 400$ grids respectively(Sec.~\ref{2d-2CGNN-leapfrog}.)}\label{tb: Summary of cross-training errors, 2D Riemann Prob. }

\begin{tabular}{|c|c|c|c|c|c|} \hline
& \multicolumn{5}{|c|}{\textbf{2-D Riemann Problems}}\\ \hline
& Config. 1 & Config. 2 & Config. 4 & Config. 6 & Config. 8\\ \hline
Initial Value & \multicolumn{5}{|c|}{Relative $l_2$ Errors}\\ \hline
Original & $1.78e-03$ & $3.27e-03$ & $5.71e-03$ & $6.56e-03$ & $4.03e-03$\\ \hline
$+3\%$ & $1.69e-03$ & $3.01e-03$ & $4.95e-03$ & $3.89e-03$ & $3.86e-03$\\ \hline
$-3\%$ & $1.67e-03$ & $3.07e-03$ & $3.88e-03$ & $3.71e-03$ & $4.07e-03$\\ \hline
$+5\%$ & $1.59e-03$ & $3.05e-03$ & $3.92e-03$ & $4.11e-03$ & $3.72e-03$\\ \hline
$-5\%$ & $1.63e-03$ & $3.05e-03$ & $4.06e-03$ & $4.02e-03$ & $4.94e-03$\\ \hline

\end{tabular}
\end{table}

\begin{table}[H]\centering
\caption{Relative $l_2$ errors of 2CGNN predictions of solutions of the 2D Euler system with initial values perturbed from those of Configuration 3, with 2 different types of inputs: (1) Low-fidelity input solutions computed by leapfrog and diffusion splitting scheme (\ref{leapfrog-diffusion-splitting, 2D}) on $200\times 200$ and $400\times 400$ grids respectively; (2) low-fidelity input solutions computed by a $4$th order scheme on $100\times 100$ and $200\times 200$ grids respectively (Sec.~\ref{2d-2CGNN-leapfrog}.)}\label{tb: Summary of 2D 2CGNN errors, config. 3 }

\begin{tabular}{|c|c|c|} \hline
& \multicolumn{2}{|c|}{\textbf{2-D Riemann Problem, Config. 3}}\\ \hline
& Input by leapfrog and & Input by a 4th \\ 
& diffusion splitting scheme & order scheme \\ \hline
Initial Value & \multicolumn{2}{|c|}{Relative $l_2$ Errors}\\ \hline
Original & $1.71e-02$ & $7.81e-03$ \\ \hline
$+3\%$ & $2.62e-02$ & $2.99e-03$ \\ \hline
$-3\%$ & $2.77e-02$ & $2.00e-03$ \\ \hline
$+5\%$ & $4.02e-02$ & $3.77e-03$ \\ \hline
$-5\%$ & $1.55e-02$ & $3.88e-03$ \\ \hline

\end{tabular}
\end{table}

\begin{table}[H]\centering
\caption{Relative $l_2$ errors of 2DCNN predictions of solutions of the 2D Euler system with initial values perturbed from those of Configuration 8, with low-fidelity inputs computed by 3 different methods: (1) leapfrog and diffusion splitting scheme (\ref{leapfrog-diffusion-splitting, 2D}) with diffusion coefficients $\Delta x$ and $4\Delta x$ respectively on $400 \times 400$ grid; (2) leapfrog and diffusion splitting scheme (\ref{leapfrog-diffusion-splitting, 2D}) with diffusion coefficient $4\Delta x$ on $400 \times 400$ grid and a $4$th order scheme on $200 \times 200$ grid respectively; (3) leapfrog and diffusion splitting scheme (\ref{leapfrog-diffusion-splitting, 2D}) with diffusion coefficient $\Delta x$ on $400 \times 400$ grid and a $4$th order scheme on $200 \times 200$ grid respectively (Sec.~\ref{Sec: 2DCNN-2D}.)}\label{tb: Summary of 2D 2DCNN errors, config. 8}

\begin{tabular}{|c|c|c|c|} \hline
& \multicolumn{3}{|c|}{\textbf{2-D Riemann Problem, Config. 8}}\\ \hline
& Input by leapfrog and & Input by leapfrog and  & Input by leapfrog and\\ 
& diffusion splitting & diffusion splitting & diffusion splitting \\ 
& with diffusion coeff. & with diffusion coeff. $4\Delta x$ & with diffusion coeff. $\Delta x$\\ 
&  $\Delta x$ and $4\Delta x$ & and a 4th order scheme & and a 4th order scheme \\ \hline
Initial Value & \multicolumn{3}{|c|}{Relative $l_2$ Errors}\\ \hline
Original & $3.95e-03$ & $3.05e-03$ & $2.02e-03$ \\ \hline
$+3\%$ & $3.16e-03$ & $1.97e-03$ & $1.99e-03$ \\ \hline
$-3\%$ & $2.97e-03$ & $2.57e-03$ & $1.93e-03$ \\ \hline
$+5\%$ & $2.95e-03$ & $1.80e-03$ & $1.92e-03$ \\ \hline
$-5\%$ & $3.10e-03$ & $2.17e-03$ & $2.02e-03$ \\ \hline

\end{tabular}
\end{table}

\section{Conclusion}
\label{Sec: conclusion}
We have extended the neural network method presented in previous work~\cite{2CGNN_1D_RiemannProb_21} to 2D and developed several variants of the method. In particular, in one variant we introduce the time step size into the input so as to allow the computation of inputs to use different time step sizes.
In another variant, we use a high order scheme to compute inputs on two coarse grids, which improves the prediction of solutions with very fine detail.
In other variants, we use a low order scheme and a high order scheme to compute inputs on coarse grids, and this can be viewed as a generalization of 2DCNN. Numerical experiments show that the method is of relatively low cost to train, robust and efficient.
 
\bibliography{mybibfile}

\end{document}